\documentclass[a4paper,10pt]{article}

\usepackage[title]{appendix}
\usepackage{setspace}
\usepackage[bottom]{footmisc}
\usepackage{color}
\usepackage{subfig}
\usepackage{svg}
\usepackage{hyperref}
\usepackage{multirow}
\usepackage{cite}
\usepackage{makecell}
\usepackage{amsmath} 
\usepackage{mathtools}
\usepackage{cleveref}
\usepackage{graphicx}
\usepackage{epsfig}
\usepackage{epstopdf}
\usepackage{float}
\usepackage{titling}
\usepackage[margin=1.0in]{geometry}
\usepackage{etoolbox}

\usepackage{standalone}     

\usepackage{tikz}
\usetikzlibrary{calc, patterns}
\usepackage{tkz-base}       
\usepackage{tkz-euclide}
\usetikzlibrary{shapes.geometric, shapes.arrows}
\usetikzlibrary{decorations.text}

\usepackage{color, colortbl}
\usepackage{enumitem}

\usepackage{siunitx}  
\usepackage{booktabs} 

\newcommand{\volume}{{\ooalign{\hfil$V$\hfil\cr\kern0.08em--\hfil\cr}}}

\usepackage{algorithm}
\usepackage{algpseudocode}

\title{A study of troubled-cell indicators applied to finite volume methods
using a novel monotonicity parameter}
\makeatletter
\renewcommand\@date{{%
  \vspace{-\baselineskip}%
  \large\centering
  \begin{tabular}{@{}c@{}}
    R. Shivananda Rao\textsuperscript{1} \\
    \normalsize nandu.rsr@gmail.com
  \end{tabular}%
  \quad and\quad
  \begin{tabular}{@{}c@{}}
    M Ramakrishna\textsuperscript{2} \\
    \normalsize krishna@ae.iitm.ac.in
  \end{tabular}

  \bigskip

  \textsuperscript{1}Dept. of Aerospace Engg., IIT Madras.\par
  \textsuperscript{2}Professor, Dept. of Aerospace Engg., IIT Madras.

  \bigskip

  \today
}}
\makeatother
\begin{document}

\maketitle

\begin{abstract}
 We adapt a troubled-cell indicator developed for discontinuous Galerkin (DG) methods \cite{FU2017} to the finite volume method (FVM) framework for solving hyperbolic conservation laws. This indicator depends solely on the cell-average data of the target cell and its immediate neighbours. Once the troubled-cells are identified, we apply the limiter only in these cells instead of applying in all computational cells. We introduce a novel technique to quantify the quality of the solution in the neighbourhood of the shock by defining a monotonicity parameter $\mu$. Numerical results from various two-dimensional simulations on the hyperbolic systems of Euler equations using a finite volume solver employing MUSCL reconstruction validate the performance of the troubled-cell indicator and the approach of limiting only in the troubled-cells. These results show that limiting only in the troubled-cells is preferable to limiting everywhere as it improves convergence without compromising on the solution accuracy.
 
\end{abstract}

\section{Introduction}
\label{sec:Intro}

In this paper, we adapt a troubled-cell indicator originally developed for discontinuous Galerkin (DG) methods to finite volume methods (FVM). We introduce a novel technique to quantify the quality of the solution in the neighbourhood of the shock. We show that limiting only in troubled-cells is preferable to limiting everywhere.

The Finite Volume Method (FVM) is a widely employed numerical technique for solving hyperbolic conservation laws arising in various scientific and engineering disciplines, particularly compressible fluid dynamics \cite{LeVeque2002} due to its robustness and adaptability across diverse applications. The Monotone Upstream-centered Schemes for Conservation Laws (MUSCL) \cite{VANLEER1979} reconstruction is widely used within the finite volume framework to attain high-order accuracy by reconstructing cell interface values from neighboring cell-average data. However, like other high-order numerical schemes, MUSCL face challenges in accurately capturing shocks and discontinuities without introducing non-physical oscillations, which can potentially lead to numerical instability.

Limiting techniques are often incorporated into the MUSCL reconstruction to prevent such oscillations. These limiters control spurious oscillations in the vicinity of discontinuities, thereby preserving monotonicity while retaining high-order accuracy of the numerical scheme in regions of smooth solution \cite{VANLEER1979,VANLEER1974}. In conventional implementations of the MUSCL reconstruction, the limiter function is typically applied to all computational cells. However, this function remains inactive in most of the cells, particularly, in regions of smooth solution and it often identifies regions near smooth extrema as requiring limiting \cite{BISWAS1994}. This approach guarantees monotonicity but it can lead to unnecessary computational effort and may adversely affect the convergence of the solution to a steady state \cite{WAN2022}. In the subsequent sections, we refer to this as the ``\textbf{limiting everywhere approach}''.

Krivodonova et al. \cite{KRIVODONOVA2004} proposed a discontinuity detector for Discontinuous Galerkin (DG) methods. By applying limiters exclusively in the computational cells suggested by discontinuity detector, they demonstrated a significant improvement over the conventional approach of applying limiters in all computational cells. Qiu and Shu \cite{SHU2005} proposed a weighted essentially non-oscillatory (WENO) finite volume methodology as limiters for DG methods. Their approach involves detecting ``troubled-cells" that potentially contain discontinuities and may require limiting. Within these troubled cells, the original high-order DG solution polynomials are replaced with WENO reconstructed polynomials to retain the original high-order accuracy of the DG method.

In the present paper, we investigate this approach of limiting only in the regions near discontinuities for an FVM solver employing MUSCL reconstruction. To detect such regions near discontinuities, where limiting is needed, we adapt the troubled-cell indicator developed by Fu and Shu \cite{FU2017} for DG methods to our FVM framework. To validate the effectiveness of this limiting approach, we employ a series of well-defined two-dimensional simulations on the hyperbolic systems of Euler equations involving discontinuities, particularly shocks. Limiting in troubled-cells alone shows better convergence without compromising the solution accuracy as compared to limiting in all computational cells. In the subsequent sections, we refer to the method of applying the limiter function only in troubled-cells as the ``\textbf{limiting restricted region approach}''.

We also present a new technique to quantify the quality of the solution in the neighbourhood of the shock by defining a monotonicity parameter $\mu$. Using the monotonicity parameter, we compare the solutions obtained from the limiting restricted region approach for different threshold constants of the troubled-cell indicator, as well as the solution obtained from the limiting everywhere approach for the given test problem.

The outline of the rest of the paper is as follows. In section \ref{sec:Methodology}, we briefly review the two-dimensional compressible Euler equations governing inviscid flow, along with the FVM employing the MUSCL reconstruction scheme to solve these equations. This section also presents the description of various two-dimensional test cases employed in this study. Section \ref{sec:TCI} introduces the troubled-cell indicator, adapted from DG methods to the FVM framework, and demonstrates its effectiveness across various test cases. Section \ref{sec:Limiting} presents numerical results for the test problems outlined in Section \ref{sec:Methodology}, by solving them using both limiting approaches and comparing the solutions through various metrics, including the monotonicity parameter. Finally, we conclude in section \ref{sec:Conclusion}.

\section{Numerical Methodology}
\label{sec:Methodology}
In this section, we provide a concise overview of the governing equations and the base solver employed to solve them. We describe the various two-dimensional test cases, both steady-state and unsteady problems, employed in this study.

The base solver utilizes a finite volume approach with MUSCL reconstruction and the Koren limiter \cite{HEMKER1988}. No other scheme is used to improve the quality of the solution or enhance convergence to steady-state so as to isolate the effect of the limiting strategy. As a consequence, there are instances where the base solver stalls and convergence to steady-state is not achieved.

\subsection*{Governing Equations and Base Solver:}
The governing equations for compressible inviscid flows can be written in integral form over a control volume with volume $\mathcal{V}$ and surface area $S$ as
\begin{equation}\label{eq:2d_euler}
\frac{d}{dt}\int_{\mathcal{V}}\textbf{Q} \,d\mathcal{V} + \int_{S} \textbf{H}(\textbf{Q}) \cdot \,\hat{\textbf{n}}\, dS= 0.
\end{equation}

In two dimensions, the vector of conservative variables $\textbf{Q} = [\rho, \rho u, \rho v, \rho E]^T$, the convective flux vector $\textbf{H} = (\textbf{F}, \textbf{G})$, $\textbf{F} = [\rho u, \rho u^2 + p, \rho uv, \rho uH]^T$, $\textbf{G} = [\rho v, \rho uv, \rho v^2 + p, \rho vH]^T$, pressure $p = (\gamma - 1)\rho\left[E - \frac{1}{2}(u^2 + v^2)\right]$, enthalpy $H = E + p/\rho$ and $\gamma = 1.4$.

The region of interest in the flow field is discretized and replaced by an array of non-overlapping cells (or volumes). For a given cell, using the finite volume method, equation (\ref{eq:2d_euler}) can be discretized in space to obtain the semi-discretized form
\begin{equation}\label{eq:semidiscrete}
 \Omega\,\frac{d\bar{\textbf{Q}}}{dt} + \textbf{R}(\bar{\textbf{Q}})= 0
\end{equation}
where, $\bar{\textbf{Q}}$ is the cell average of $\textbf{Q}$ and $\Omega$ is the volume of the cell. The residue $\textbf{R}(\bar{\textbf{Q}})$ is given as
\begin{equation}\label{eq:residue}
\textbf{R}(\bar{\textbf{Q}}) = \sum_{f} \tilde{\textbf{H}}_f(\bar{\textbf{Q}}) \cdot \,\hat{\textbf{n}}_f \,s_f
\end{equation}
where, $f$ is an index over the faces of the cell, $\tilde{\textbf{H}}_f$ is the numerical inviscid flux vector at the face, $\hat{\textbf{n}}_f$ is the unit normal vector at the face and $s_f$ is the area of the cell face.

In this investigation, we use the Advection Upstream Splitting Method (AUSM) family scheme AUSM+ \cite{LIOU1996} to evaluate the numerical inviscid flux $\tilde{\textbf{H}}_f$. This requires the state to the left and right of the face where $\tilde{\textbf{H}}_f$ is being evaluated.

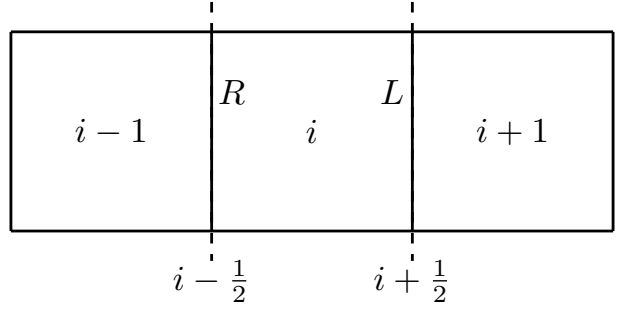
\begin{figure}
\centering
\fontsize{10}{12}\selectfont

\begin{tikzpicture}
 \centering

\draw[thick] (0,0) -- (6,0);   
\draw[thick] (0,-2) -- (6,-2); 

\draw[thick](0,0) -- (0,-2);
\draw[thick](2,0) -- (2,-2);
\draw[thick](4,0) -- (4,-2);
\draw[thick](6,0) -- (6,-2);

\draw (1, -1.0) node[black] {$i-1$};
\draw (3, -1.0) node[black] {$i$};
\draw (5, -1.0) node[black] {$i+1$};

\draw[dashed, thick](2,0.3) -- (2,-2.3);
\draw[dashed, thick](4,0.3) -- (4,-2.3);

\draw (2, -2.5) node[black]{$i-\frac{1}{2}$};
\draw (4, -2.5) node[black]{$i+\frac{1}{2}$};

\draw (2.2, -0.6) node[black] {$R$};
\draw (3.8, -0.6) node[black] {$L$};

%

\end{tikzpicture}
\caption{Stencil of a cell $i$ for the MUSCL reconstruction scheme to obtain the left $(L)$ state at the face $i+1/2$ and the right $(R)$ state at the face $i-1/2$.}
\label{fig:MUSCL_setup}
\end{figure}

We use the higher order Monotone Upwind Schemes for Scalar Conservation Laws (MUSCL) reconstruction strategy \cite{VANLEER1979} to reconstruct the left and the right states at the interface. In the MUSCL reconstruction strategy, the left $(L)$ and the right $(R)$ states at the faces of a cell $i$ can be computed using the cell averages of cells $i-1$, $i$, and $i+1$. These three cells together constitute the stencil $(i-1,\, i,\, i+1)$ for cell $i$ and are shown in Figure \ref{fig:MUSCL_setup}.

The left and right states at the faces of a cell $i$ using the MUSCL reconstruction scheme in terms of primitive variables are determined as follows.
\begin{equation}\label{eq:Muscl_left_wol}
 \textbf{W}^L_{i+\frac{1}{2}} = \textbf{W}_i + \frac{1}{4}\left[(1 - k)\Delta^-\textbf{W}_i + (1 + k)\Delta^+\textbf{W}_i\right]
\end{equation}
\begin{equation}\label{eq:Muscl_right_wol}
 \textbf{W}^R_{i-\frac{1}{2}} = \textbf{W}_{i} - \frac{1}{4}\left[(1 + k)\Delta^-\textbf{W}_{i} + (1 - k)\Delta^+\textbf{W}_{i}\right]
\end{equation}
where, $\textbf{W}$ is the vector of primitive variables, i.e., $\textbf{W} = (\rho, u, v, p)^T$. $\textbf{W}^L_{i+\frac{1}{2}}$ is the desired left state at the face $i+\frac{1}{2}$ and $\textbf{W}^R_{i-\frac{1}{2}}$ is the desired right state at the face $i-\frac{1}{2}$. The forward difference ($\Delta^+$) and the backward difference ($\Delta^-$) operators are defined as $\Delta^+\textbf{W}_i = \textbf{W}_{i+1} - \textbf{W}_i$, $\Delta^-\textbf{W}_i = \textbf{W}_{i} - \textbf{W}_{i-1}$. The parameter $k$ controls the upwind biasing and is taken here as 1/3 to achieve a quadratic reconstruction for smooth solutions.

Traditionally, a limiter is used in the reconstruction to prevent spurious oscillations due to discontinuities in the solution. The left and right states in the MUSCL reconstruction scheme with a slope limiter, $\Phi$, are given by
\begin{equation}\label{eq:Muscl_left_wl}
 \textbf{W}^L_{i+\frac{1}{2}} = \textbf{W}_i + \frac{1}{4}\left[(1 - k) \, \Phi\left(r\right)\, \Delta^-\textbf{W}_i + (1 + k)\, \Phi\left(\frac{1}{r}\right) \,\Delta^+\textbf{W}_i\right]
\end{equation}
\begin{equation}\label{eq:Muscl_right_wl}
 \textbf{W}^R_{i-\frac{1}{2}} = \textbf{W}_{i} - \frac{1}{4}\left[(1 + k) \, \Phi\left(r\right) \, \Delta^-\textbf{W}_{i} + (1 - k) \, \Phi\left(\frac{1}{r}\right) \, \Delta^+\textbf{W}_{i}\right]
\end{equation}
where, $r = \displaystyle \frac{\Delta^+w}{\Delta^-w}$ is the ratio of forward and backward differences of respective primitive variable ($w$). For the slope limiter with symmetry property $\left( \text{i.e.,} \, \Phi(r) = \Phi\left(\frac{1}{r}\right)\right)$, the limited MUSCL reconstruction equations (\ref{eq:Muscl_left_wl}) and (\ref{eq:Muscl_right_wl}) becomes
\begin{equation}\label{eq:Muscl_left_wl_simplified}
 \textbf{W}^L_{i+\frac{1}{2}} = \textbf{W}_i + \frac{1}{2}\, \Psi^L(r)\, \Delta^-\textbf{W}_i
\end{equation}
\begin{equation}\label{eq:Muscl_right_wl_simplified}
 \textbf{W}^R_{i-\frac{1}{2}} = \textbf{W}_{i} - \frac{1}{2}\, \Psi^R(r)\, \Delta^-\textbf{W}_{i}
\end{equation}
with the limiter function $\Psi$ defined as
\begin{equation}
 \Psi^L(r) = \frac{1}{2}\left[(1 - k) + (1 + k) r \right]\Phi(r)
\end{equation}
\begin{equation}
 \Psi^R(r) = \frac{1}{2}\left[(1 + k) + (1 - k) r\right]\Phi(r)
\end{equation}

In the present investigation, we used a slope limiter $\Phi$ given by
\begin{equation}\label{eq:vanAlbada}
 \Phi(r) = \displaystyle \frac{3r}{2r^2 - r + 2}
\end{equation}
with $k = 1/3$. For this slope limiter, the limiter function $\Psi(r)$ corresponds to the limiter of Hemker and Koren \cite{HEMKER1988}.

The semi-discretized form given in equation (\ref{eq:semidiscrete}) can be solved by integrating it in time to obtain the evolution of dependent variables. For this purpose, an implicit matrix-free Lower-Upper Symmetric Gauss-Seidel (LU-SGS) \cite{SHAROV1997} method is utilized with a Courant-Friedrichs-Lewy (CFL) number set to 1 for steady-state problems. For unsteady problems, we use a Total-Variation-Diminishing Runge-Kutta Method (TVD RK3) \cite{SHU1988}, with a CFL number of 0.3.

\subsection*{Steady-state Test Cases:}
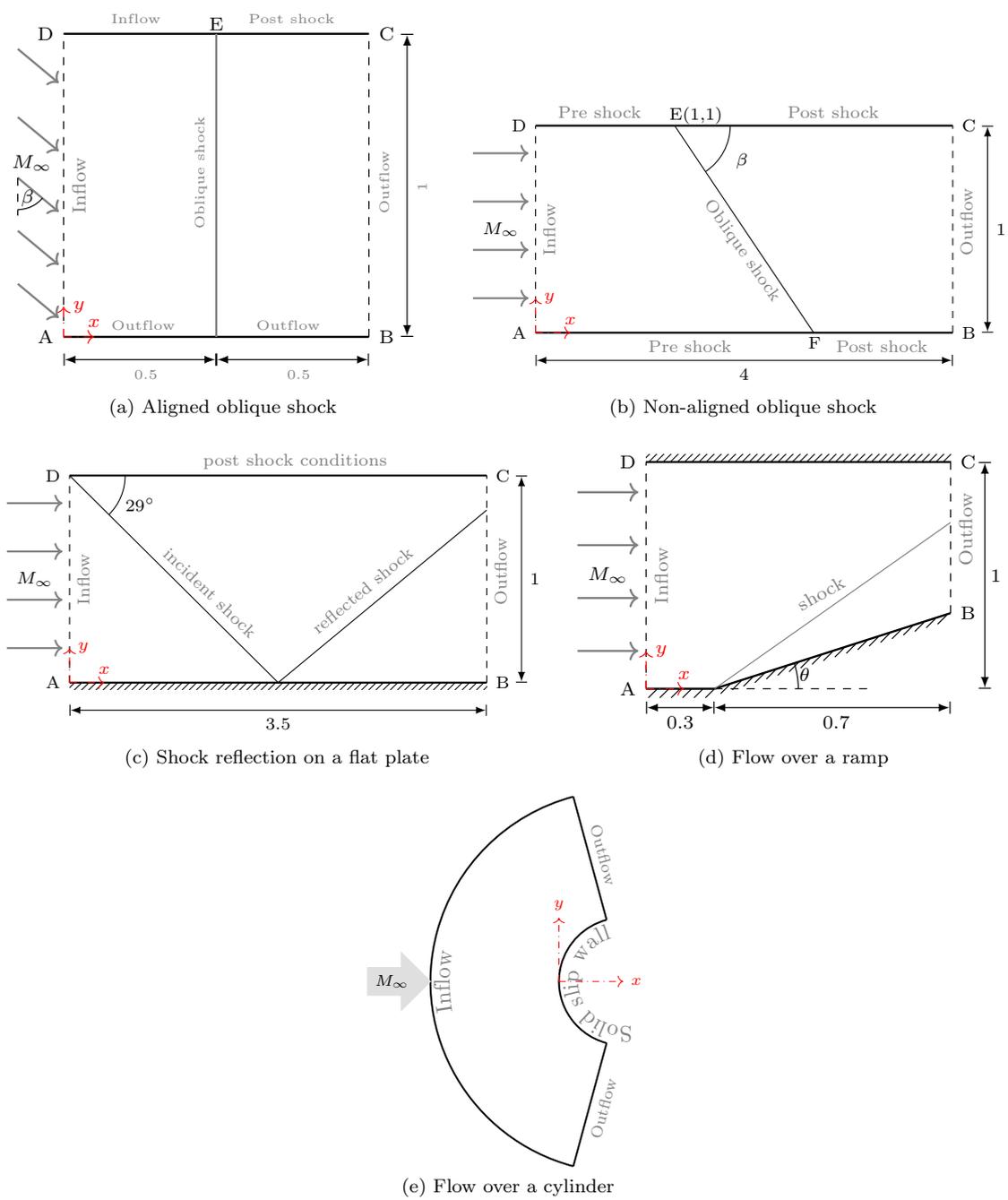
\begin{figure}
\centering
\subfloat[Aligned oblique shock]{

 \begin{tikzpicture}[font=\scriptsize]

 \centering

\draw[thick] (0,0)node[left]{A} -- ++(0:4cm)node[right]{B};       
\draw[thick] (0,4)node[left, rotate=0]{D} -- ++(0:4cm)node[right, rotate=0]{C};      

\draw(2, 3.92) node[above]{E};

\draw[dashed] (0,0) -- ++(90:4cm);         
\draw[dashed] (4,0) -- ++(90:4cm);      
\draw[thick,gray] (2,0) -- ++(90:4cm);    

\draw (0.5, 4.2)node[gray,right, rotate=0]{\tiny Inflow};
\draw (2.3, 4.2)node[gray,right,rotate=0]{\tiny Post shock};

\draw[>=latex, |<->|] (4.5,0.0) -- ++(90:4cm);
\draw (4.7, 1.8)node[gray,right, rotate=90]{\tiny 1};

\draw[>=latex, |<->|] (0.0,-0.3) -- ++(0:2cm);
\draw (0.8, -0.5)node[gray,right, rotate=0]{\tiny 0.5};

\draw[>=latex, |<->|] (2,-0.3) -- ++(0:2cm);
\draw (2.8, -0.5)node[gray,right, rotate=0]{\tiny 0.5};

\draw[thick,gray, ->] (-0.6,3.8) -- ++(-40:0.7cm);

\draw[thick,gray, ->] (-0.6,2.9) -- ++(-40:0.7cm);
\draw[thick,gray, ->] (-0.6,2.1) -- ++(-40:0.7cm);

\draw (-0.8, 2.3)node[black,right]{$M_{\infty}$};


\draw[dashdotted,red,->] (0,0) -- (0.4,0)node[above]{$x$};
\draw[dashdotted,red, ->] (0,0) -- (0,0.4)node[right]{$y$};

\draw[thick,gray, ->] (-0.6,1.4) -- ++(-40:0.7cm);
\draw[thick,gray, ->] (-0.6,0.7) -- ++(-40:0.7cm);

 \coordinate (A) at (-0.6,2.1);
\coordinate  (C) at (-0.6,1.6);
\coordinate (B) at (0.0,1.6);
\tkzMarkAngle[size=0.42cm](C,A,B)
\tkzLabelAngle[pos = 0.3, font=\scriptsize](C,A,B){$\beta$}
\draw[dashed](-0.6,2.1)--(-0.6,1.6);

\draw (1.8,1.3)node[gray,right, rotate=90]{\tiny Oblique shock};

\draw (4.2,1.5)node[gray,right, rotate=90]{\tiny Outflow};
\draw (0.5,0.15)node[gray,right, rotate=0]{\tiny Outflow};
\draw (2.4,0.15)node[gray,right, rotate=0]{\tiny Outflow};

\draw (0.2,1.5)node[gray,right, rotate=90]{Inflow};

\end{tikzpicture} \label{fig:AlignOS_setup}} \hspace{0.3cm}
\subfloat[Non-aligned oblique shock]{

 \begin{tikzpicture}[font=\scriptsize]

 \centering

\draw[thick] (0,0)node[left]{D} -- (6,0)node[right]{C};   
  
\draw[thick] (0,-3) node[left]{A}-- (6,-3)node[right]{B}; 

\draw[dashed](0,-3) -- (0,0); 

\draw[dashed](6,-3) -- (6,0); 

\draw (2,0) -- (4, -3);		

\coordinate (A) at (2cm,0cm);
\coordinate  (C) at (6.0cm,0.0cm);
\coordinate (B) at (4.0cm,-3.0cm);

\tkzMarkAngle[size=0.8cm](B,A,C)
\tkzLabelAngle[pos = 1.1, font=\scriptsize](B,A,C){$\beta$}


\node [] at (-0.5,-1.5) {$M_{\infty}$};

\draw[dashdotted,red,->] (0,-3) -- (0.5,-3)node[above]{$x$};
\draw[dashdotted,red, ->] (0,-3) -- (0,-2.5)node[right]{$y$};

\draw (3.5, 0.2)node[gray,right]{Post shock};
\draw (0.2, 0.2)node[gray,right]{Pre shock};

\draw (4.2, -3.2)node[gray,right]{Post shock};
\draw (1.5, -3.2)node[gray,right]{Pre shock};

\draw (2.4, -1.0)node[gray,right,rotate=-55]{Oblique shock};

\draw[>=latex, |<->|] (6.5,0) -- (6.5,-3);
\draw[>=latex, |<->|] (0,-3.4) -- (6,-3.4);

\draw (6.5, -1.5)node[black, right]{1};
\draw (3.0, -3.4)node[black, below]{4};

\draw (1.8,0.15)node[black, right]{E(1,1)};
\draw (3.8,-3.15)node[black, right]{F};
\draw (6.2,-2.0)node[gray,right, rotate=90]{Outflow};
\draw (0.2,-2.0)node[gray,right, rotate=90]{Inflow};

\draw[thick,gray, ->] (-0.9,-0.4) -- ++(0:0.8cm);
\draw[thick,gray, ->] (-0.9,-1.1) -- ++(0:0.8cm);
\draw[thick,gray, ->] (-0.9,-1.8) -- ++(0:0.8cm);
\draw[thick,gray, ->] (-0.9,-2.5) -- ++(0:0.8cm);

\end{tikzpicture}\label{fig:NonAlignOS_setup}} \\
\subfloat[Shock reflection on a flat plate]{

 \begin{tikzpicture}[font=\scriptsize]

 \centering

\draw[thick] (0,0)node[left]{D} -- (6,0)node[right]{C};   
  
\draw[thick] (0,-3) node[left]{A}-- (6,-3)node[right]{B}; 
\fill[pattern = north east lines] ($ (0,-3) $) rectangle ($ (6,-3.1) $);

\draw[dashed](0,-3) -- (0,0); 

\draw[dashed](6,-3) -- (6,0); 

\draw (0,0) -- (3, -3);		

\coordinate (A) at (0.0cm,0.0cm);
\coordinate  (C) at (6.0cm,0.0cm);
\coordinate (B) at (3.0cm,-3.0cm);

\tkzMarkAngle[size=0.8cm](B,A,C)
\tkzLabelAngle[pos = 1.1, font=\scriptsize](B,A,C){\ang{29}}

\draw (3,-3) -- (6, -0.5);




\draw[dashdotted,red,->] (0,-3) -- (0.5,-3)node[above]{$x$};
\draw[dashdotted,red, ->] (0,-3) -- (0,-2.5)node[right]{$y$};


\draw (1.8, 0.2)node[gray,right]{post shock conditions};

\draw (1.3, -1.0)node[gray,right,rotate=-44]{incident shock};
\draw (5.0, -1.0)node[gray,left,rotate = 40]{reflected shock};


%

\draw[>=latex, |<->|] (6.5,0) -- (6.5,-3);
\draw[>=latex, |<->|] (0,-3.4) -- (6,-3.4);

\draw (6.5, -1.5)node[black, right]{1};
\draw (3.0, -3.6)node[black]{3.5};

\draw[thick,gray, ->] (-0.9,-0.4) -- ++(0:0.8cm);
\draw[thick,gray, ->] (-0.9,-1.1) -- ++(0:0.8cm);
\draw[thick,gray, ->] (-0.9,-1.8) -- ++(0:0.8cm);
\draw[thick,gray, ->] (-0.9,-2.5) -- ++(0:0.8cm);
\node [] at (-0.5,-1.5) {$M_{\infty}$};

\draw (6.2,-2.0)node[gray,right, rotate=90]{Outflow};
\draw (0.2,-2.0)node[gray,right, rotate=90]{Inflow};

\end{tikzpicture}\label{fig:SR_setup}} \hspace{0.3cm}
\subfloat[Flow over a ramp]{

 \begin{tikzpicture}[font=\scriptsize]

\centering


\draw[thick] (0,0)node[left]{D} -- (4,0)node[right]{C};   


\draw[thick] (0,-3) node[left]{A}-- (0.9,-3); 

\draw[thick] (0.9,-3) -- (4,-2)node[right]{B}; 



\draw[dashed](0,-3) -- (0,0); 

\draw[dashed](4,-2) -- (4,0); 

\draw[dashed] (0.9,-3) -- (3.0, -3);		

\coordinate (A) at (0.9cm,-3.0cm);
\coordinate  (C) at (3.0cm,-3.0cm);
\coordinate (B) at (4.0cm,-2.0cm);

\tkzMarkAngle[size=1.1cm](C,A,B)
\tkzLabelAngle[pos = 1.2, font=\scriptsize](C,A,B){$\theta$}

\draw[gray] (0.9,-3) -- (4, -0.8);
\draw (2.7, -1.5)node[gray,left, rotate = 30]{shock};




\draw[dashdotted,red,->] (0,-3) -- (0.5,-3)node[above]{$x$};
\draw[dashdotted,red, ->] (0,-3) -- (0,-2.5)node[right]{$y$};



%



\draw[>=latex, |<->|] (0,-3.25) -- (0.9,-3.25);
\draw (0.45, -3.45)node[black]{0.3};

\draw[>=latex, <->|] (0.9,-3.25) -- (4.0,-3.25);
\draw (2.5, -3.45)node[black]{0.7};


\draw[>=latex, |<->|] (4.45,0) -- (4.45, -3);
\draw (4.4, -1.5)node[black, right]{1};

\draw[thick,gray, ->] (-0.9,-0.4) -- ++(0:0.8cm);
\draw[thick,gray, ->] (-0.9,-1.1) -- ++(0:0.8cm);
\draw[thick,gray, ->] (-0.9,-1.8) -- ++(0:0.8cm);
\draw[thick,gray, ->] (-0.9,-2.5) -- ++(0:0.8cm);
\node [] at (-0.5,-1.5) {$M_{\infty}$};

\draw (4.2,-1.5)node[gray,right, rotate=90]{Outflow};
\draw (0.2,-2.0)node[gray,right, rotate=90]{Inflow};

\fill[pattern = north east lines] ($ (0,0) $) rectangle ($ (4,0.1) $);








\def\slope{(1/3.1} 
    \def\intercept{-3.3} 


    \foreach \x in {1.1,1.25,...,4} {
        \pgfmathsetmacro\y{\slope*\x+\intercept} 
        \draw[] (\x,\y) -- ++(-0.12,-0.12); 
    }

    \foreach \x in {0.15,0.3,...,0.95} {
        \draw[] (\x,-3) -- ++(-0.12,-0.12); 
    }


\end{tikzpicture}\label{fig:Ramp_setup}} \\
\subfloat[Flow over a cylinder]{ \begin{tikzpicture}[font=\scriptsize]
 \centering

\draw [decoration={text along path, text={Solid slip wall},text align={right}, text color={gray}},decorate] (-0.7 * 0.258819, -0.7 * 0.965926) arc (255:105:0.7cm);


\draw[thick] (-0.258819, -0.965926) arc (255:105:1cm){};
\draw[thick] (-3 * 0.258819, -3 * 0.965926) arc (255:105:3cm)node[midway,sloped,above,rotate = 0]{};


\draw [decoration={text along path, text={Inflow},text align={center}, text color={gray}},decorate] (-2.7 * 0.258819, -2.7 * 0.965926) arc (255:105:2.7cm);



\draw[thick] (-0.258819, -0.965926) -- (-3 * 0.258819, -3 * 0.965926) node[gray,midway,sloped,below]{Outflow};
\draw[thick] (-0.258819, 0.965926) -- (-3 * 0.258819, 3 * 0.965926) node[gray,sloped,above,midway]{Outflow};

\draw[dashdotted,red,->] (-1,0) -- (-0.0,0) node[right] {$x$};
\draw[dashdotted,red, ->] (-1,0) -- (-1,1.0) node[above] {$y$};

\node [fill=gray!25, style={single arrow}] at (-3.6,0) {$M_{\infty}$};

\end{tikzpicture}\label{fig:BB_setup}}
\caption{Schematic for various steady-state test cases, illustrating the varying geometric and flow conditions, and boundary conditions. (not drawn to scale).}
\label{fig:Steady_setup}
\end{figure}

Test Case 1: Aligned Oblique Shock. We refer to shocks that coincide with the grid lines as aligned shocks. The computational domain for the test case is $[0, 1] \times [0, 1]$ and is discretized into a grid of $200 \times 200$ cells. Figure (\ref{fig:AlignOS_setup}) illustrates the computational setup of the test case. We compute this problem for a Mach number of 3 with shock angles of \ang{30} and \ang{40}. Boundary conditions are applied as shown in Figure (\ref{fig:AlignOS_setup}). \\ \\

Test Case 2: Non-Aligned Oblique Shock. We refer to shocks that do not coincide with the grid lines as non-aligned shocks. The computational domain for this test case is $[0, 4] \times [0, 1]$ and is divided into a grid of $800 \times 200$ cells. We compute this problem for a Mach number of 3 with shock angles of \ang{30} and \ang{40}. Boundary conditions are applied as shown in Figure (\ref{fig:NonAlignOS_setup}). \\ \\

Test Case 3: Shock Reflection On A Flat Plate. For this test case, the computational domain is $[0, 3.5] \times [0, 1]$ and is divided into $700 \times 200$ cells. The computational setup for this test case is shown in Figure (\ref{fig:SR_setup}). Supersonic inflow and outflow boundary conditions are applied along the edges AD and BC, respectively. Along the edge AB solid slip wall boundary conditions and along DC post-shock conditions characterized by an inlet Mach number of 2.9 and a shock wave angle of \ang{29}. \\ \\

Test Case 4: Flow Over A Ramp. The computational domain for this test case is $[0, 1] \times [0, 1]$ and is divided into $200 \times 200$ cells. A detailed illustration of this test case setup is presented in Figure (\ref{fig:Ramp_setup}). We compute this problem for an inlet Mach number of 3 with ramp angles of \ang{20} and \ang{30}. Supersonic inlet and outlet boundary conditions are applied along the edges AD and BC, respectively and solid slip wall boundary conditions are applied along the edges AB and DC. \\ \\

Test Case 5: Flow Over A Cylinder. The computational grid consists of 160 (radial) $\times$ 320 (circumferential) cells which are equally spaced in each direction. The cylinder radius is 0.5 units and the inflow boundary is two times the radius away from the wall. The computational setup for this problem is shown in Figure (\ref{fig:BB_setup}). We compute this problem for a Mach number of 3. The supersonic inflow boundary condition is used at the inlet. The solid slip wall boundary condition is applied at the cylinder wall boundary, and the other two are supersonic outflow boundaries. \\ \\

\subsection*{Unsteady Test Cases:}
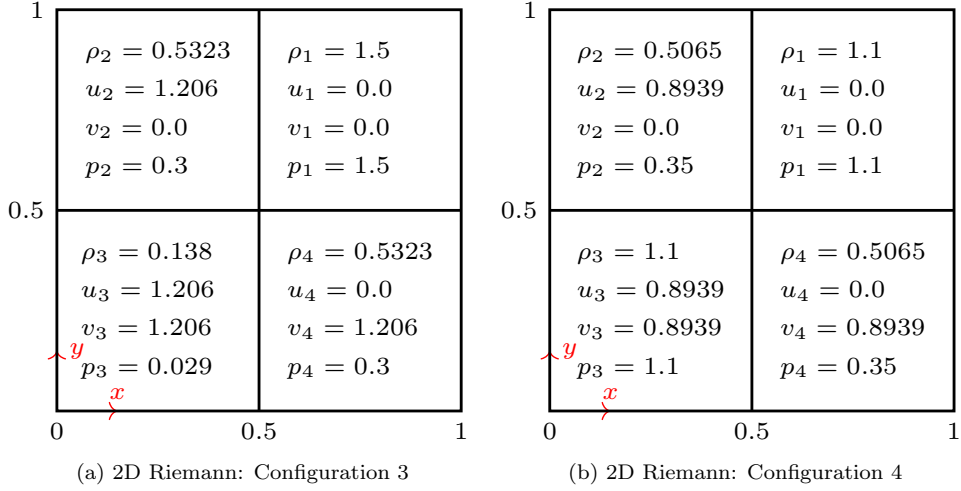
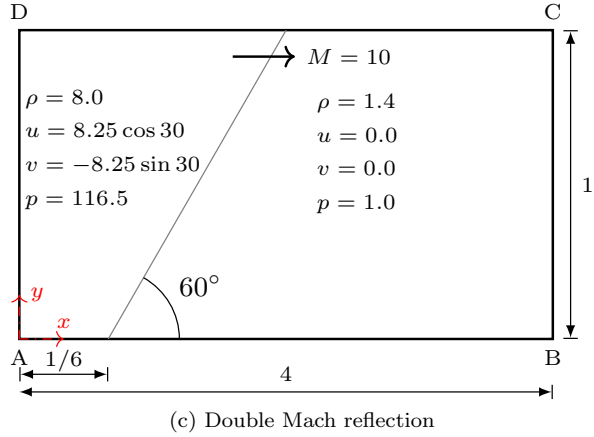
\begin{figure}
\centering
\subfloat[2D Riemann: Configuration 3]{

 \begin{tikzpicture}[font=\scriptsize]

 \centering



    \draw[dashdotted,red,->] (0,0) -- (0.6,0)node[above]{$x$};
    \draw[dashdotted,red, ->] (0,0) -- (0,0.6)node[right]{$y$};

    \draw[thick] (0,0) rectangle (4,4);

    \draw[thick] (0,2) -- (4,2); 
    \draw[thick] (2,0) -- (2,4); 


   \node[align=left] at (2.8,3.0) {$\rho_1=1.5$ \\[0.1cm] $u_1=0.0$ \\[0.1cm] $v_1=0.0$ \\[0.1cm] $p_1=1.5$};


    \node[align=left] at (1,3.0) {$\rho_2=0.5323$ \\[0.1cm] $u_2=1.206$ \\[0.1cm] $v_2=0.0$ \\[0.1cm] $p_2=0.3$};


      \node[align=left] at (0.9,1.0) {$\rho_3=0.138$ \\[0.1cm] $u_3=1.206$ \\[0.1cm] $v_3=1.206$ \\[0.1cm] $p_3=0.029$};


     \node[align=left] at (3,1.0) {$\rho_4=0.5323$ \\[0.1cm] $u_4=0.0$ \\[0.1cm] $v_4=1.206$ \\[0.1cm] $p_4=0.3$};

     \draw (0.0, 0.0) node[below]{0};
     \draw (2.0, 0.0) node[below]{0.5};
     \draw (4.0, 0.0) node[below]{1};
     \draw (0.0, 2.0) node[left]{0.5};
     \draw (0.0, 4.0) node[left]{1};

\end{tikzpicture} \label{fig:config3_setup}}
\subfloat[2D Riemann: Configuration 4]{

 \begin{tikzpicture}[font=\scriptsize]

 \centering



    \draw[dashdotted,red,->] (0,0) -- (0.6,0)node[above]{$x$};
    \draw[dashdotted,red, ->] (0,0) -- (0,0.6)node[right]{$y$};

    \draw[thick] (0,0) rectangle (4,4);

    \draw[thick] (0,2) -- (4,2); 
    \draw[thick] (2,0) -- (2,4); 


   \node[align=left] at (2.8,3.0) {$\rho_1=1.1$ \\[0.1cm] $u_1=0.0$ \\[0.1cm] $v_1=0.0$ \\[0.1cm] $p_1=1.1$};


    \node[align=left] at (1,3.0) {$\rho_2=0.5065$ \\[0.1cm] $u_2=0.8939$ \\[0.1cm] $v_2=0.0$ \\[0.1cm] $p_2=0.35$};


      \node[align=left] at (1.0,1.0) {$\rho_3=1.1$ \\[0.1cm] $u_3=0.8939$ \\[0.1cm] $v_3=0.8939$ \\[0.1cm] $p_3=1.1$};


     \node[align=left] at (3,1.0) {$\rho_4=0.5065$ \\[0.1cm] $u_4=0.0$ \\[0.1cm] $v_4=0.8939$ \\[0.1cm] $p_4=0.35$};

     \draw (0.0, 0.0) node[below]{0};
     \draw (2.0, 0.0) node[below]{0.5};
     \draw (4.0, 0.0) node[below]{1};
     \draw (0.0, 2.0) node[left]{0.5};
     \draw (0.0, 4.0) node[left]{1};

\end{tikzpicture} \label{fig:config4_setup}} \\
\subfloat[Double Mach reflection]{

 \begin{tikzpicture}[font=\scriptsize]

 \centering

\coordinate (A) at (0,0); 
\coordinate (B) at (6,0); 
\coordinate (C) at (6,3.5); 
\coordinate (D) at (0,3.5); 

\coordinate (E) at (1,0);
\coordinate (F) at (3,{2.02*tan(60)}); 

\draw[thick] (A) -- (B) -- (C) -- (D) -- cycle;
\draw[gray] (E) -- (F);

\draw[] pic["\small $60^\circ$", draw=, angle radius=8mm, angle eccentricity=1.5] {angle=B--E--F};


\node[align=left] at (1.05,2.15) {$\rho=8.0$ \\[0.1cm] $u=8.25 \cos30$ \\[0.1cm] $v=-8.25 \sin30$ \\[0.1cm] $p=116.5$};

\node[align=left] at (3.8,2.1) {$\rho=1.4$ \\[0.1cm] $u=0.0$ \\[0.1cm] $v=0.0$ \\[0.1cm] $p=1.0$};

\draw[dashdotted,red,->] (0,0) -- (0.5,0)node[above]{$x$};
\draw[dashdotted,red, ->] (0,0) -- (0,0.5)node[right]{$y$};

\draw[>=latex, |<->|] (0,-0.4) -- (1.0,-0.4);
\draw (0.5, -0.25)node[black]{$1/6$};

\draw[>=latex, <->|] (0.0,-0.6) -- (6.0,-0.6);
\draw (3.0, -0.4)node[black]{4};

\draw[>=latex, |<->|] (6.2,0) -- (6.2, 3.5);
\draw (6.2, 1.75)node[black, right]{1};

\draw[thick, ->](2.4, 3.2) -- (3.1, 3.2)node[right]{$M=10$};

\node[below] at (A) {A};
\node[below] at (B) {B};
\node[above] at (C) {C};
\node[above] at (D) {D};

\end{tikzpicture}\label{fig:DMR_setup}}
\caption{Solution domain and initial conditions for unsteady test cases (not drawn to scale).}
\label{fig:Unsteady_setup}
\end{figure}

Test Case 6: 2D Riemann Problems. We solve the two dimensional Euler equations (\ref{eq:2d_euler}) on the computational domain of $[0, 1] \times [0, 1]$ for configurations 3 and 4 of the 2D Riemann problems \cite{LAX1998}. The initial data for these Riemann problems consist of four constant states of the flow variables over the four sub-domains of the unit square. Figures (\ref{fig:config3_setup}) and (\ref{fig:config4_setup}) present the initial data for the configurations 3 and 4, respectively. The computational domian is discretized into a grid of $400 \times 400$ cells. We compute the solution up to $t = 0.3$ for configuration 3 and upto $t = 0.25 $ for configuration 4. \\ \\

Test Case 7: Double Mach Reflection. This problem is originally suggested by Woodward and Collela \cite{WOODWARD1984}. The computational domain for this problem is $[0, 4] \times [0, 1]$ and is divided into a grid of $960 \times 240$ cells. Initially a right-moving Mach 10 normal shock is positioned at $(x,y) = (1/6,0)$ and makes a \ang{60} angle with positive x-axis. Along the bottom boundary (i.e. $y = 0$), we assign the values of initial post-shock flow for the short region from $x = 0$ to $x = \frac{1}{6}$, and for the rest of the boundary, we use the solid slip wall boundary conditions. For the top boundary (i.e. $y = 1$), we set the values to describe the exact motion of the initial Mach 10 normal shock. The computational setup and initial conditions for this problem are given in Figure (\ref{fig:DMR_setup}). We compute the solution up to $t = 0.2$. \\ \\

\section{Troubled-cell indicator}
\label{sec:TCI}
In this section, we introduce a troubled-cell indicator adapted from DG method. This indicator depends solely on the cell-average of density of the target cell and its immediate neighbours.
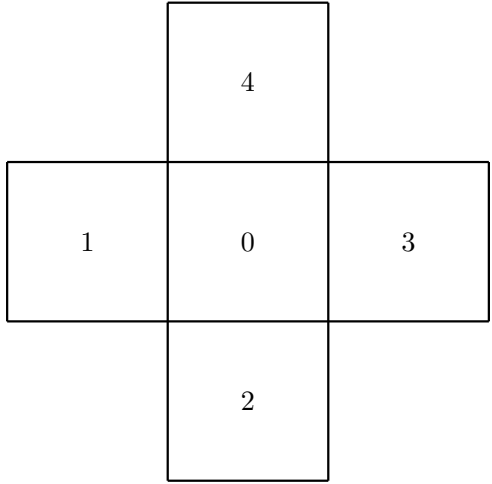
\begin{figure}
\centering

\begin{tikzpicture}
 \centering

\draw[thick] (0,0) -- (6,0);   
\draw[thick] (0,-2) -- (6,-2); 

\draw[thick] (2,-4) -- (2,2);
\draw[thick] (4,-4) -- (4,2);

\draw[thick](0,0) -- (0,-2);
\draw[thick](6,0) -- (6,-2);
\draw[thick](2,-4) -- (4,-4);
\draw[thick](2,2) -- (4,2);

\draw (3, -1) node[black] {$0$};

\draw (1, -1) node[black] {$1$};
\draw (3, -3) node[black] {$2$};
\draw (5, -1) node[black] {$3$};
\draw (3, 1) node[black] {$4$};

\end{tikzpicture}
\caption{Stencil, $S = \{C_0, C_1, C_2, C_3, C_4\}$, used (in the troubled cell indicator) to determine whether the cell $C_0$ is a troubled cell.}
\label{fig:tci_stencil}
\end{figure}
We label the target cell as $C_0$ and denote the stencil to calculate the indicator value as $S = \{C_0, C_1, C_2, C_3, C_4\}$ as shown in Figure (\ref{fig:tci_stencil}). We define the following quantity for the target cell $C_0$,
\begin{equation}\label{eq:TCI}
I_{C_0} = \displaystyle\frac{\sum_{i=1}^4|\bar{\rho}_{C_0} - \bar{\rho}_{C_i}|}{\text{max}_{i\epsilon\{0,1,2,3,4\}}\{\bar{\rho}_{C_i}\}}
\end{equation}
where, $\bar{\rho}_{C_i}$ is the cell average of density of the cell in the stencil. The cell $C_0$ is considered as a troubled cell if $I_{C_0} \geq K$ for a threshold constant $K$.

\subsection*{Performance of the indicator:}
To evaluate the performance of the troubled-cell indicator, we solve the two-dimensional Euler equations (\ref{eq:2d_euler}) for steady-state test cases described in the previous section. We compute the first-order converged solution utilizing the Lax-Friedrichs flux scheme. We use this converged solution to identify the troubled-cells. In this study, various threshold constants within the range of 0.01 to 0.1 are examined, particularly $K = 0.02$, $K = 0.05$, and $K = 0.1$. We omit the results of the constant $K = 0.05$ in the case of constant $K = 0.1$ also identifies the complete shock region.

\begin{figure}
\centering
\subfloat[$\ang{30}, K = 0.02$]{\includegraphics[width=0.5\textwidth, height=0.5\textwidth]{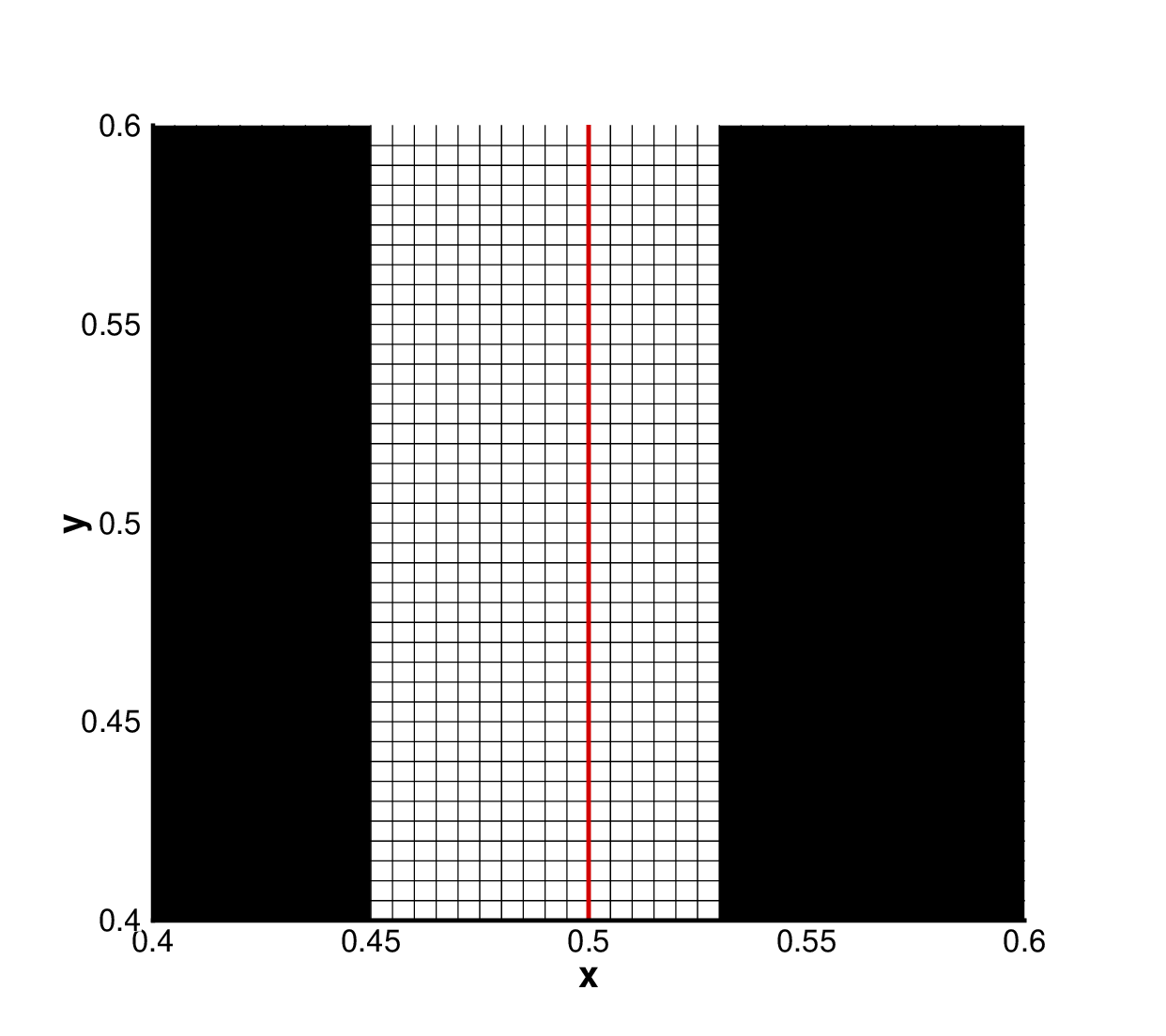}\label{fig:A30_K002_TC}}
\subfloat[$\ang{40}, K = 0.02$]{\includegraphics[width=0.5\textwidth, height=0.5\textwidth]{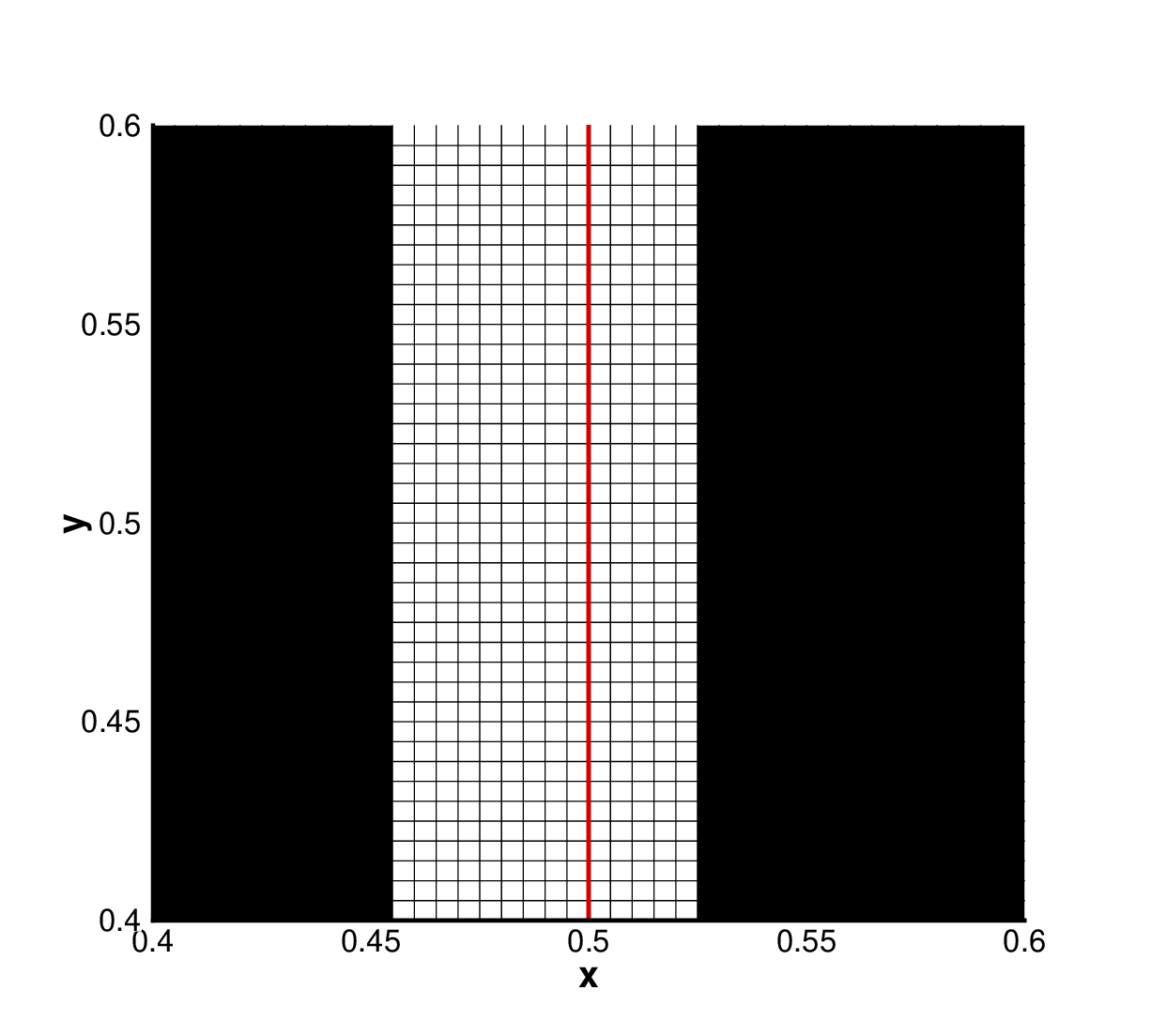}\label{fig:A40_K002_TC}} \\
\subfloat[$\ang{30}, K = 0.05$]{\includegraphics[width=0.5\textwidth, height=0.5\textwidth]{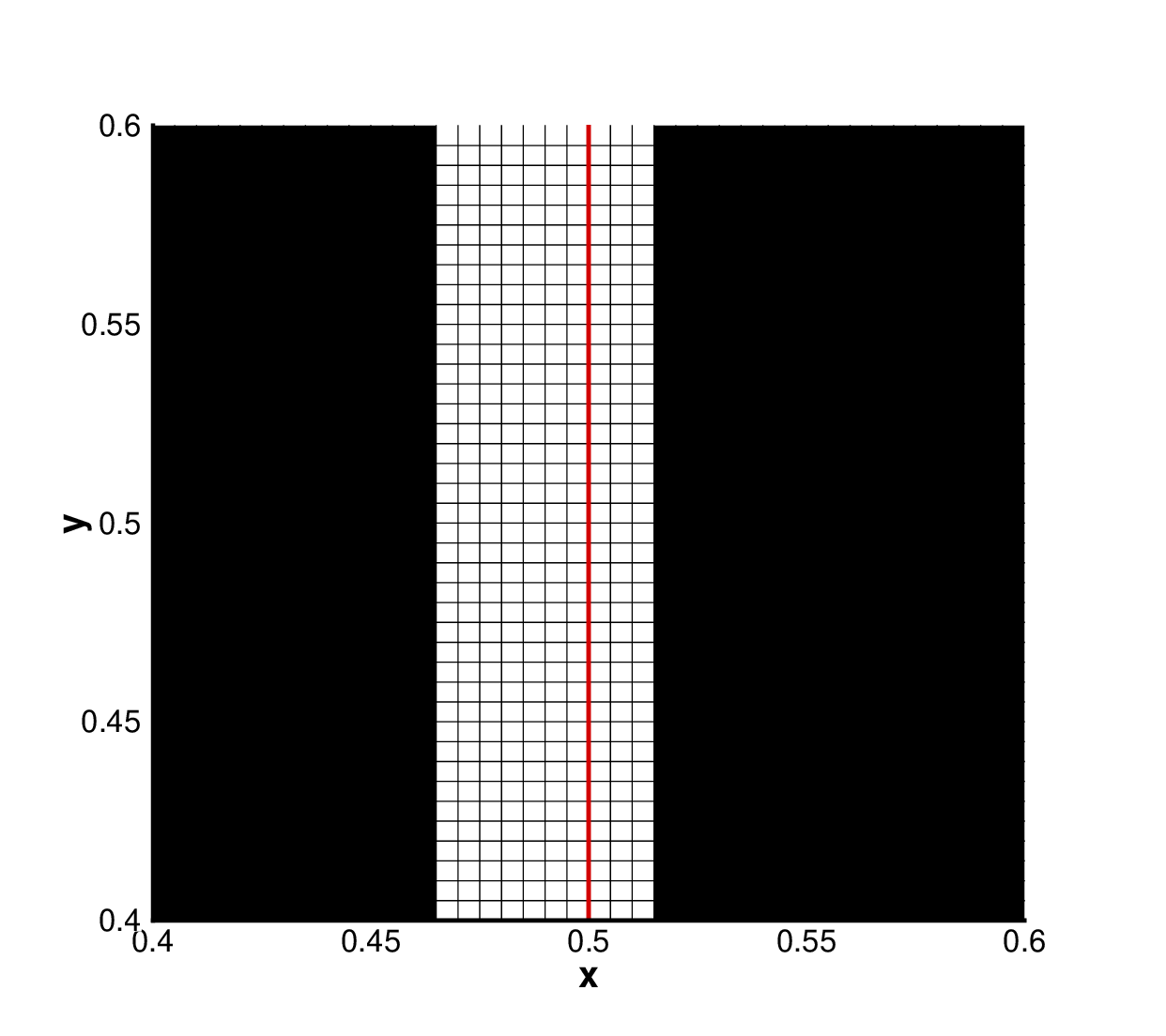}\label{fig:A30_K005_TC}}
\subfloat[$\ang{40}, K = 0.05$]{\includegraphics[width=0.5\textwidth, height=0.5\textwidth]{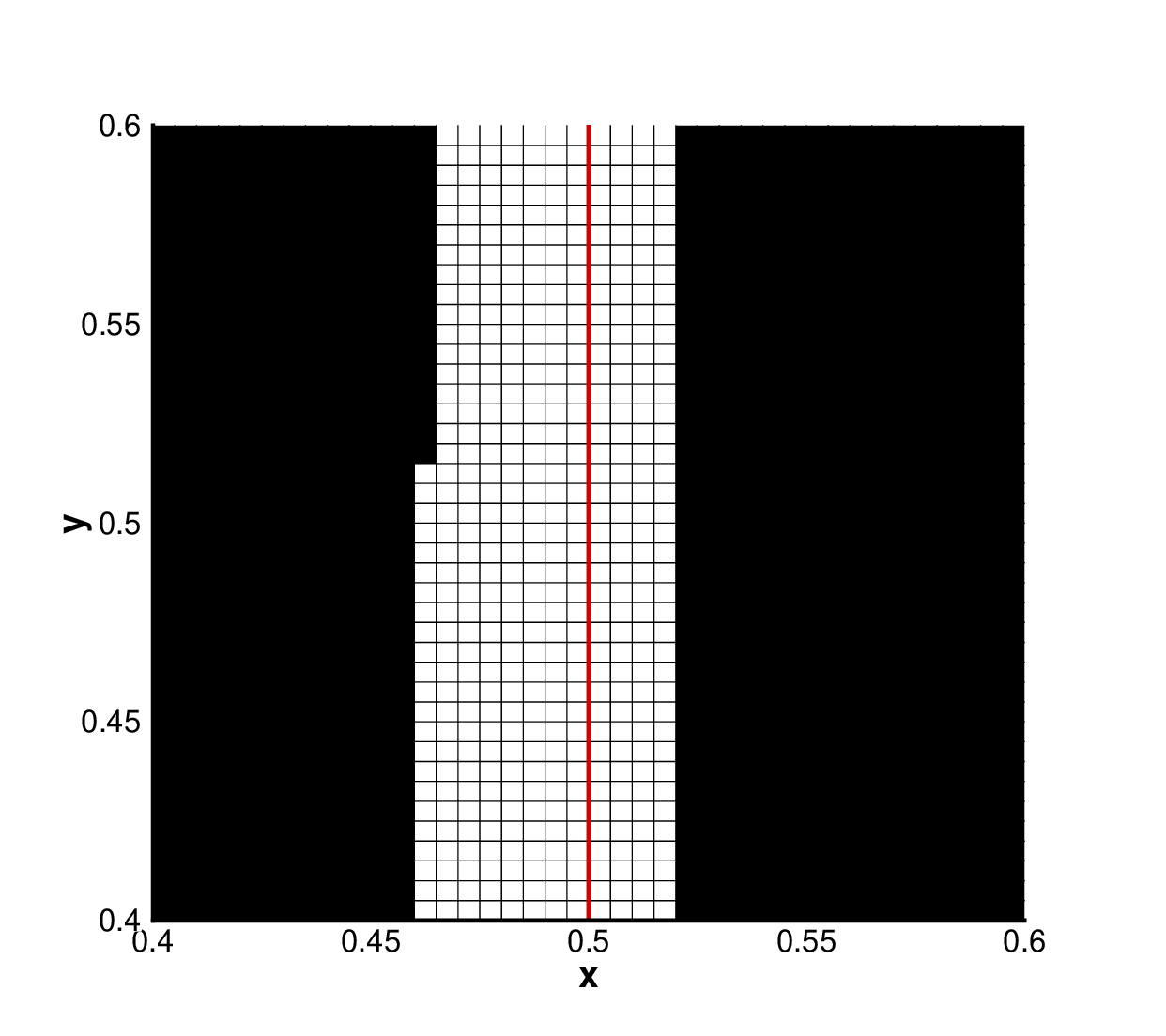}\label{fig:A40_K005_TC}}\\
\caption{Aligned oblique shock. Zoomed-in view of troubled-cells identified by the indicator for two threshold constants. Red line represents the exact shock. Percentage of number of troubled-cells identified in whole computational domain of 200 $\times$ 200 cells: (a) 7.71\% (b) 6.93\% (c) 4.87\% (d) 5.55\%}
\label{fig:AlignOS_TC}
\end{figure}

For the aligned oblique shock test case, the threshold constant $K = 0.1$ failed to identify cells in the post-shock region, particularly for the \ang{30} shock angle. In contrast, the other two threshold constants effectively detected the shock by identifying the troubled cells in the vicinity of the shock. Zoomed-in view of the troubled-cells identified by the indicator for these two threshold constants is shown in Figure ({\ref{fig:AlignOS_TC}}).

\begin{figure}
\centering
\subfloat[$\ang{30}, K = 0.02$]{\includegraphics[width=0.5\textwidth, height=0.5\textwidth]{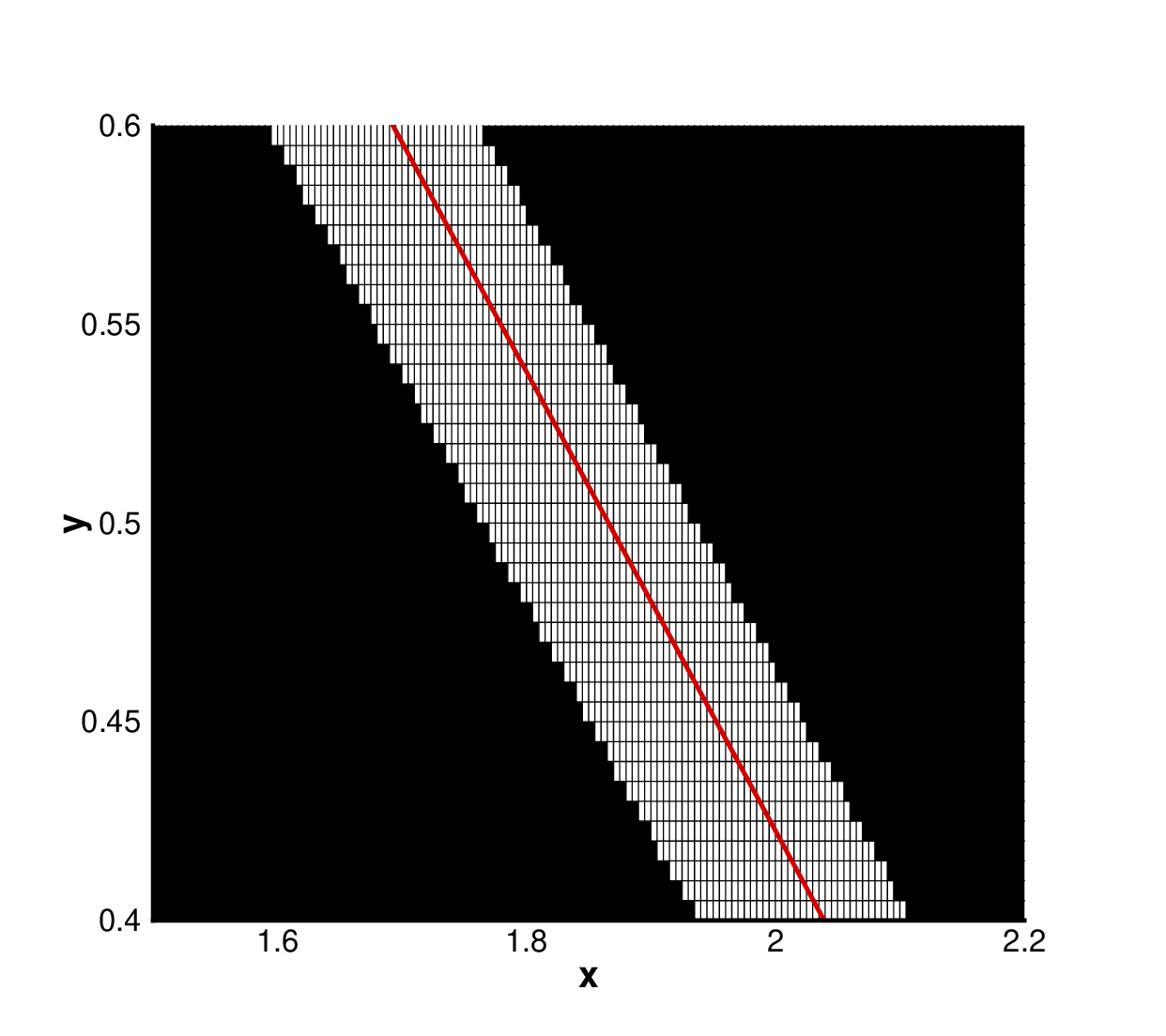}\label{fig:NA30_K002_TC}}
\subfloat[$\ang{40}, K = 0.02$]{\includegraphics[width=0.5\textwidth, height=0.5\textwidth]{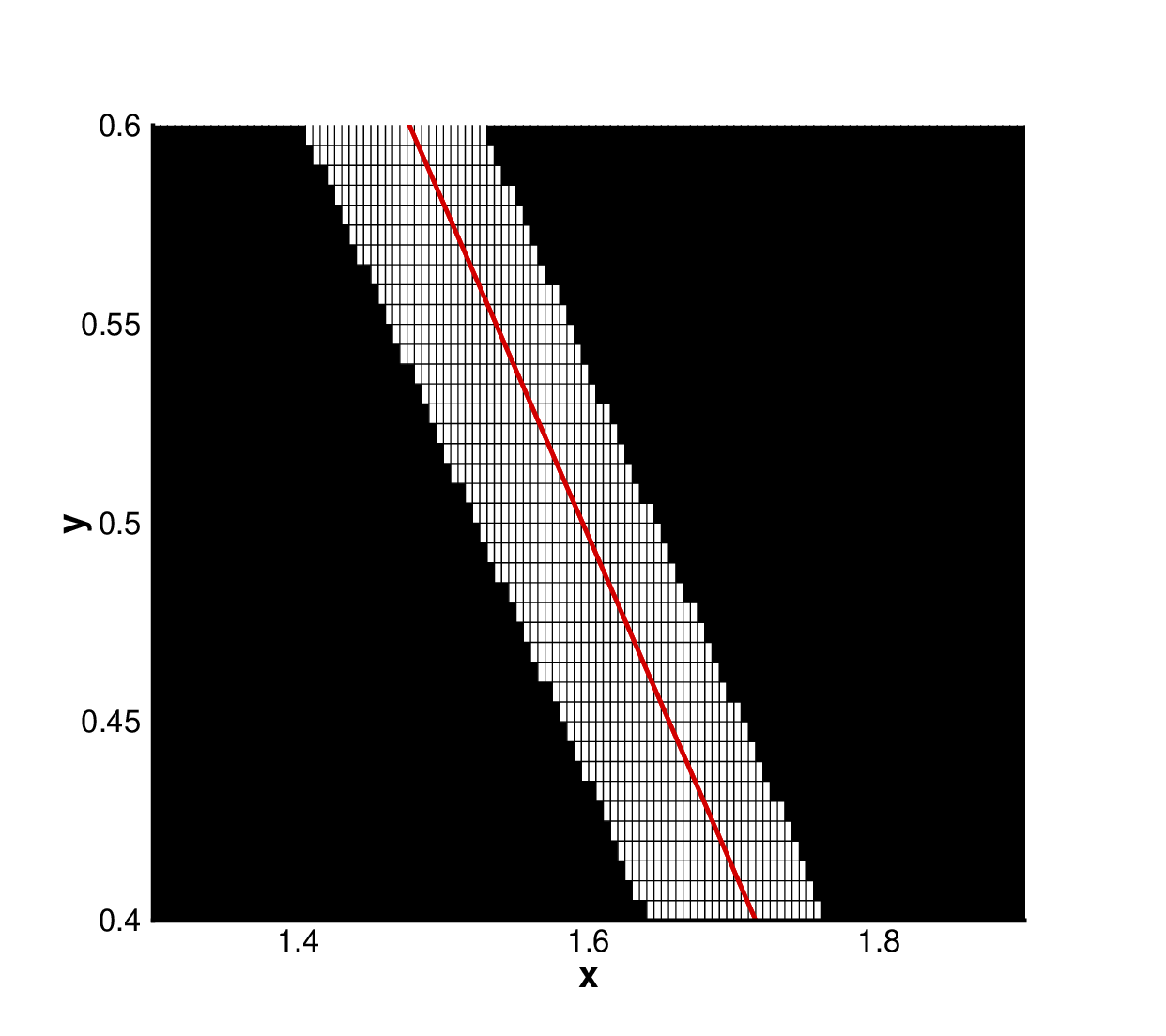}\label{fig:NA40_K002_TC}}\\
\subfloat[$\ang{30}, K = 0.1$]{\includegraphics[width=0.5\textwidth, height=0.5\textwidth]{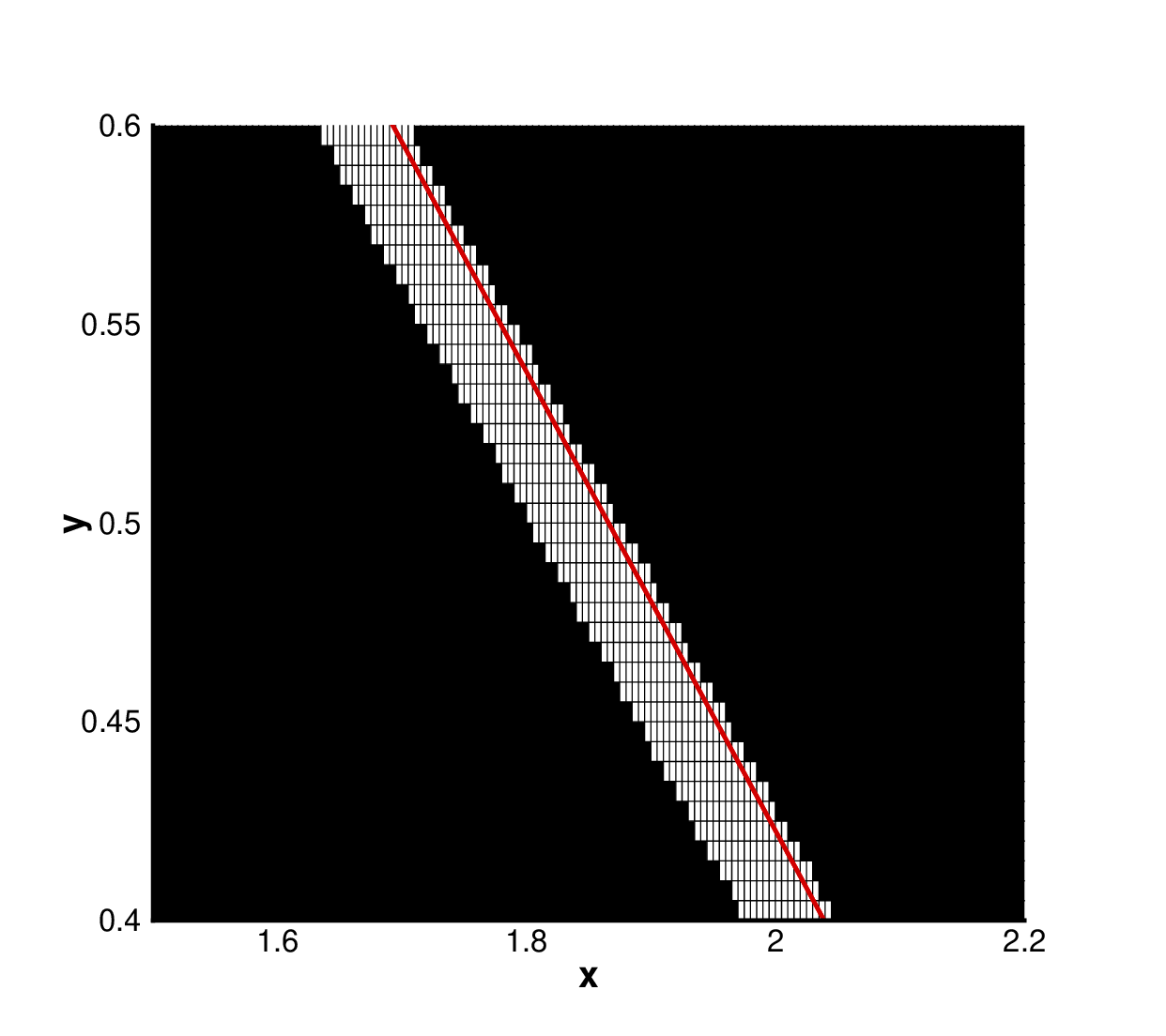}\label{fig:NA30_K01_TC}}
\subfloat[$\ang{40}, K = 0.1$]{\includegraphics[width=0.5\textwidth, height=0.5\textwidth]{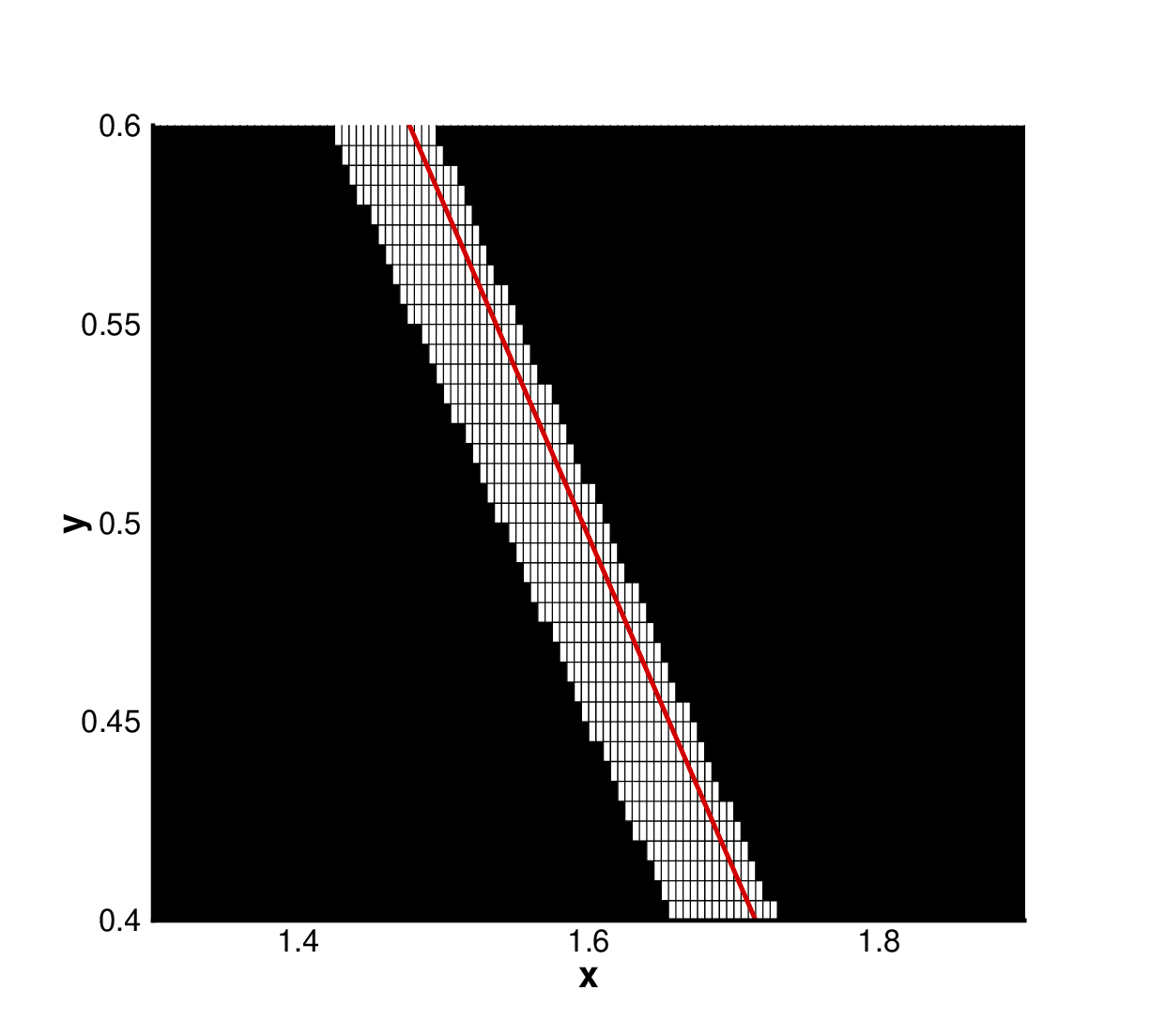}\label{fig:NA40_K01_TC}}\\
\caption{NonAligned oblique shock. Zoomed-in view of troubled-cells identified by the indicator for two threshold constants. Red line represents the exact shock. Percentage of number of troubled-cells identified in whole computational domain of 800 $\times$ 200 cells: (a) 4.11\% (b) 3.01\% (c) 1.81\% (d) 1.76\%}
\label{fig:NonAlignOS_TC}
\end{figure}

For the non-aligned oblique shock test case, the troubled-cell indicator effectively identified cells in the vicinity of the shock for all three threshold constants. Zoomed-in view of the troubled-cells identified for threshold constants $K = 0.02$ and $K = 0.1$ is illustrated in Figure ({\ref{fig:NonAlignOS_TC}}). While the identified troubled-cells completely covers the exact shock location for the threshold constant $K = 0.1$, only a few number of troubled-cells are observed in the post-shock region, particularly for a \ang{30} shock angle.

\begin{figure}
\centering
\subfloat[Incident shock, $K = 0.02$]{\includegraphics[width=0.5\textwidth, height=0.5\textwidth]{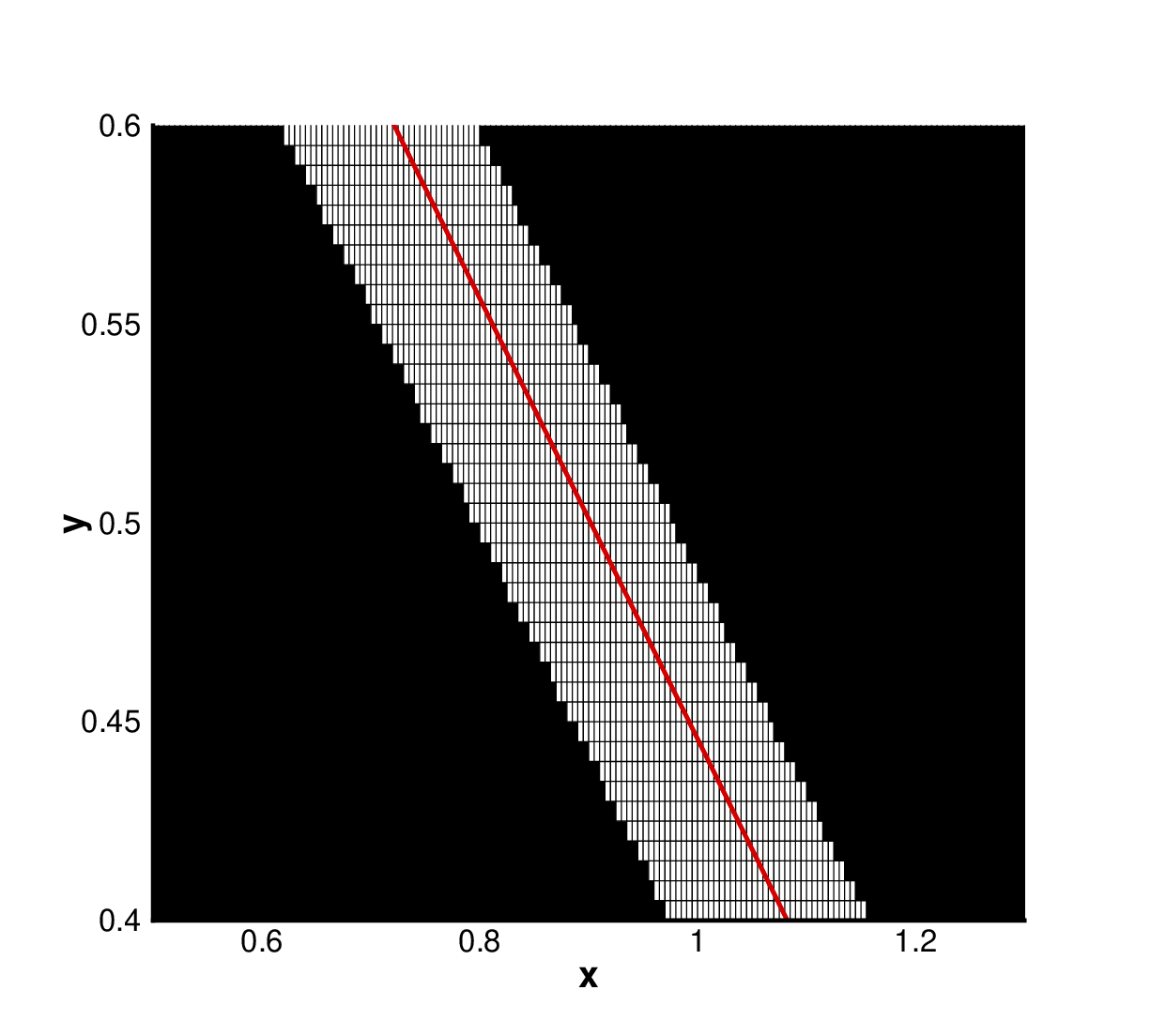}\label{fig:incident_K002_TC}}
\subfloat[Reflected shock, $K = 0.02$]{\includegraphics[width=0.5\textwidth, height=0.5\textwidth]{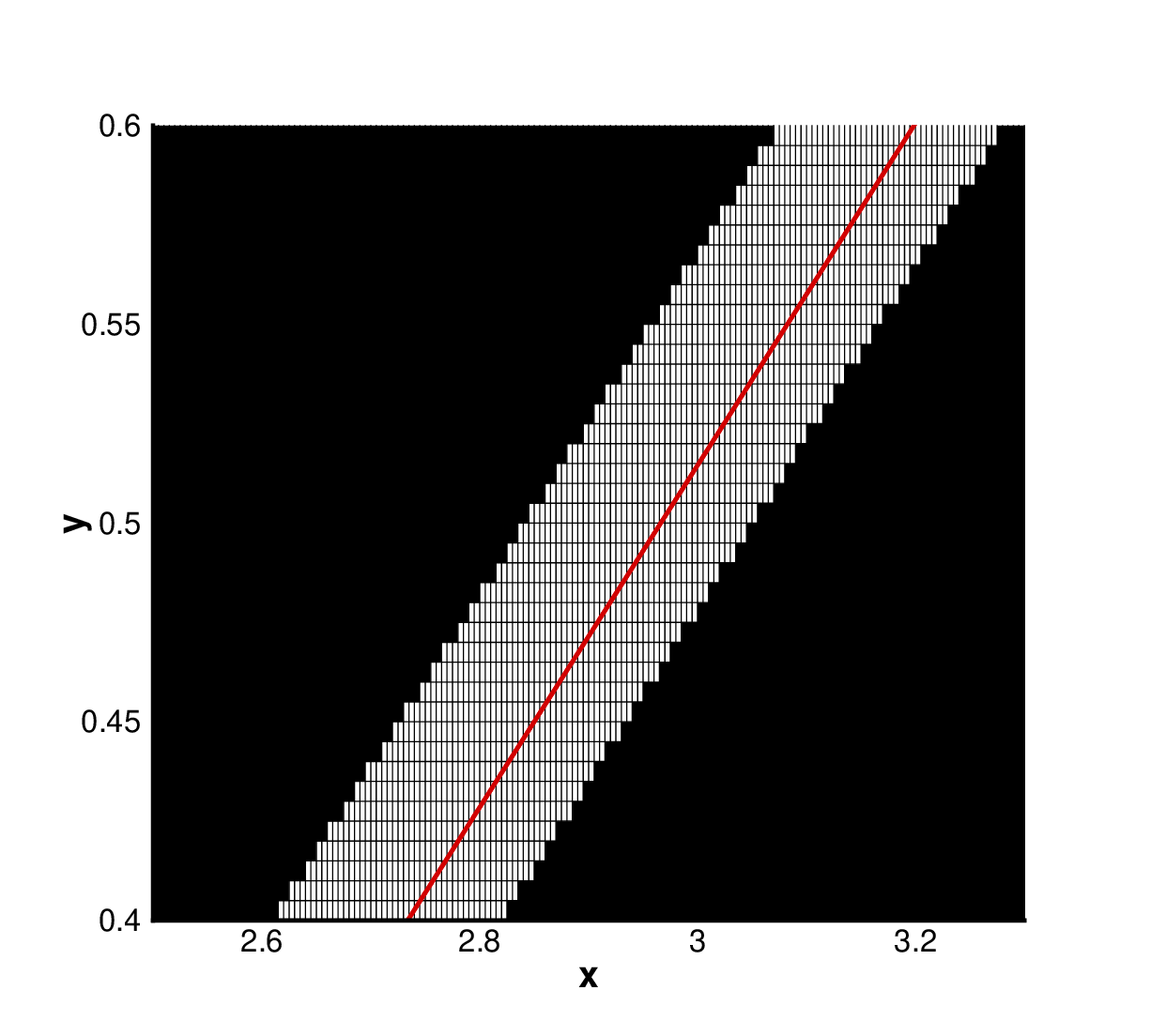}\label{fig:reflected_K002_TC}} \\
\subfloat[Incident shock, $K = 0.05$]{\includegraphics[width=0.5\textwidth, height=0.5\textwidth]{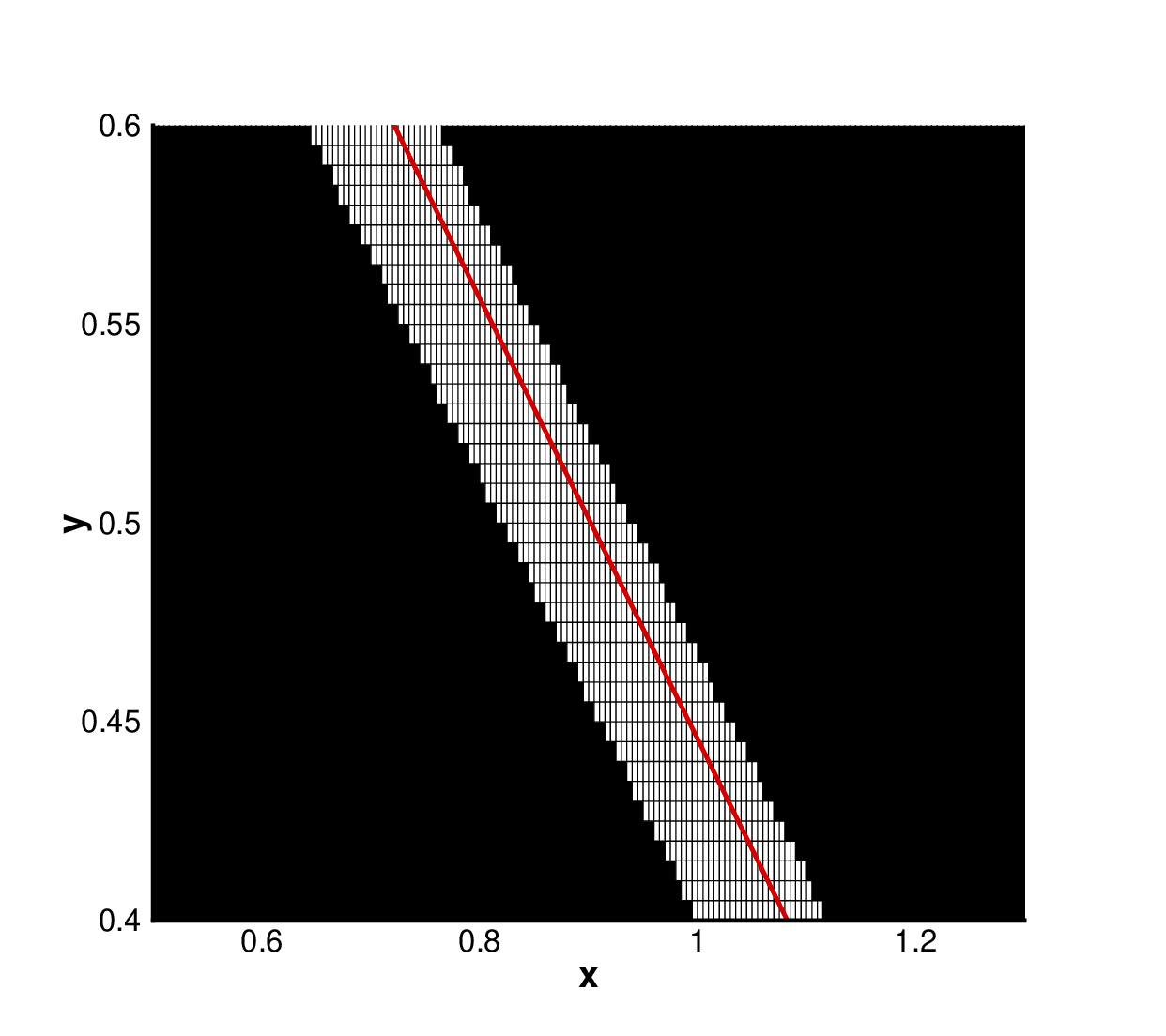}\label{fig:incident_K005_TC}}
\subfloat[Reflected shock, $K = 0.05$]{\includegraphics[width=0.5\textwidth, height=0.5\textwidth]{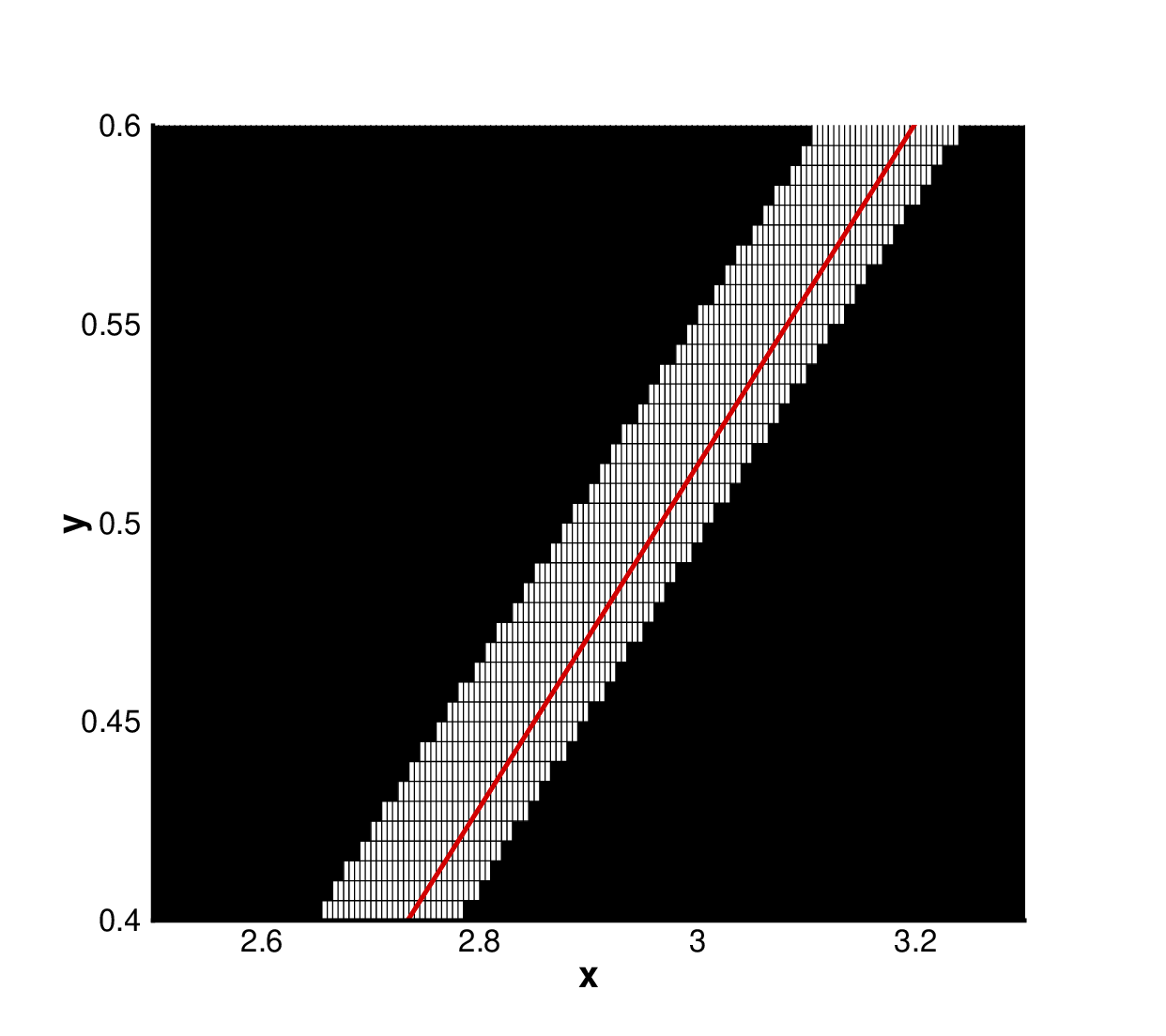}\label{fig:reflected_K005_TC}} \\
\caption{Shock reflection on a flat plate. Zoomed-in view of troubled-cells identified by the indicator for two threshold constants for both incident and reflected shocks. Red line represents the exact shock. Percentage of number of troubled-cells identified in whole computational domain of 700 $\times$ 200 cells: 9.18\% for $K = 0.02$, 6.01\% for $K = 0.05$.}
\label{fig:SR_TC}
\end{figure}

Zoomed-in view of the troubled-cells identified for threshold constants $K = 0.02$ and $K = 0.05$ is illustrated in Figure ({\ref{fig:SR_TC}}), separately highlighting both the incident and reflected shocks for the shock reflection on a flat plate test case. For these threshold constants, a sufficient number of troubled-cells are identified near both the incident and reflected shocks. However, for threshold constant $K = 0.1$, the troubled-cell indicator fails to identify an adequate number of troubled-cells in the vicinity of the reflected shock, particularly close to the plate.

\begin{figure}
\centering
\subfloat[$\ang{20}, K = 0.02$]{\includegraphics[width=0.5\textwidth, height=0.5\textwidth]{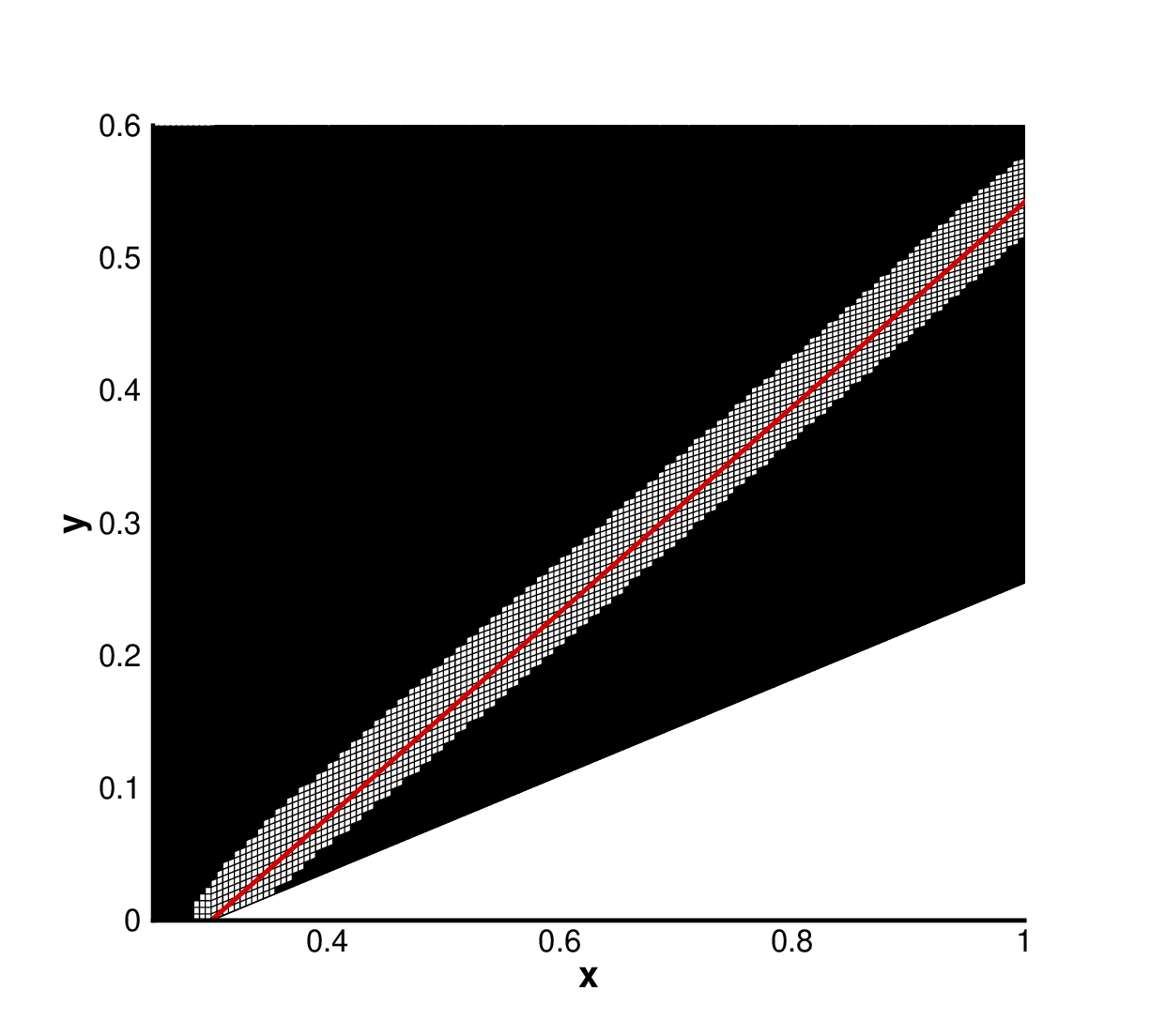}\label{fig:R20_K002_TC}}
\subfloat[$\ang{30}, K = 0.02$]{\includegraphics[width=0.5\textwidth, height=0.5\textwidth]{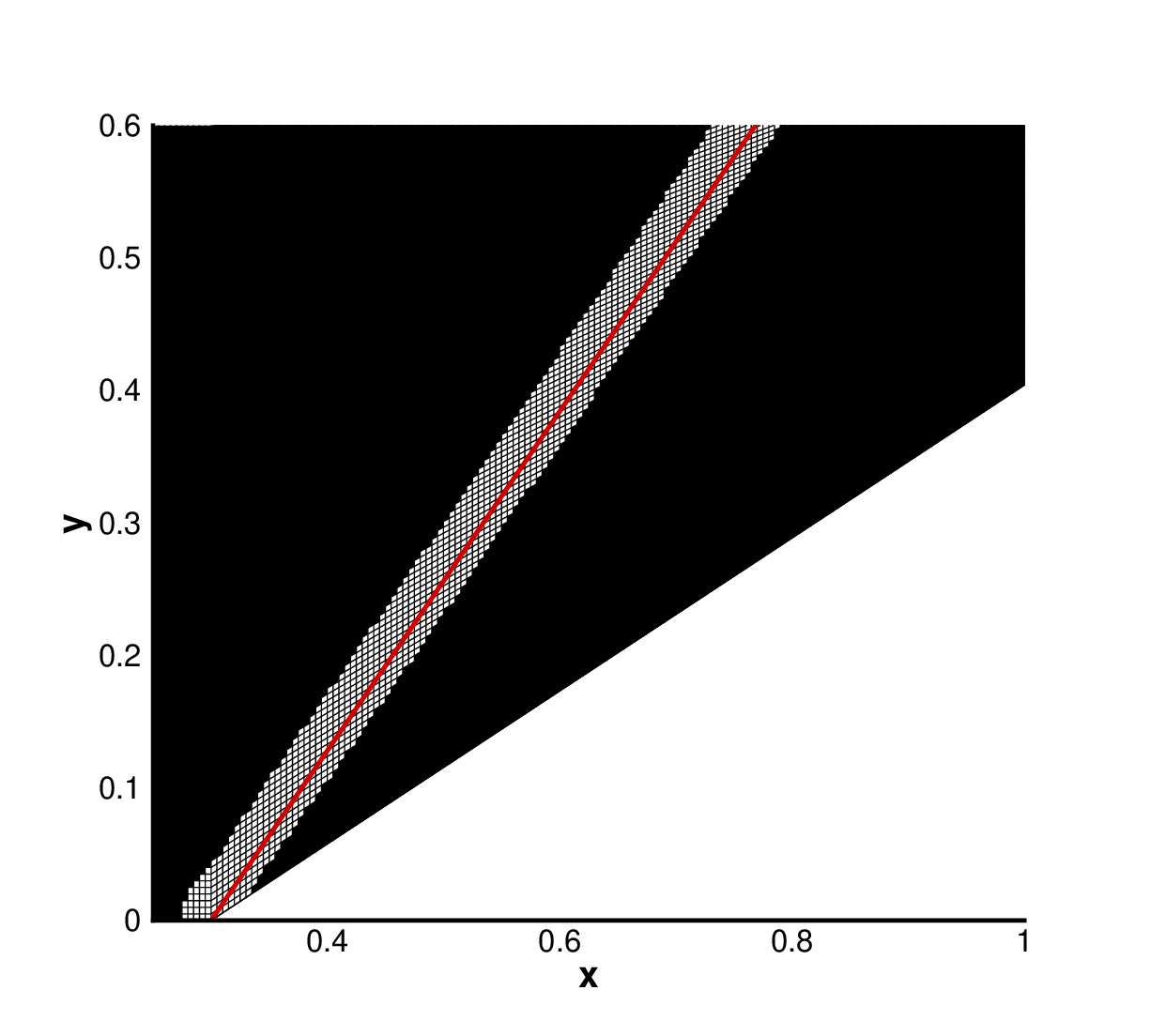}\label{fig:R30_K002_TC}}\\
\subfloat[$\ang{20}, K = 0.1$]{\includegraphics[width=0.5\textwidth, height=0.5\textwidth]{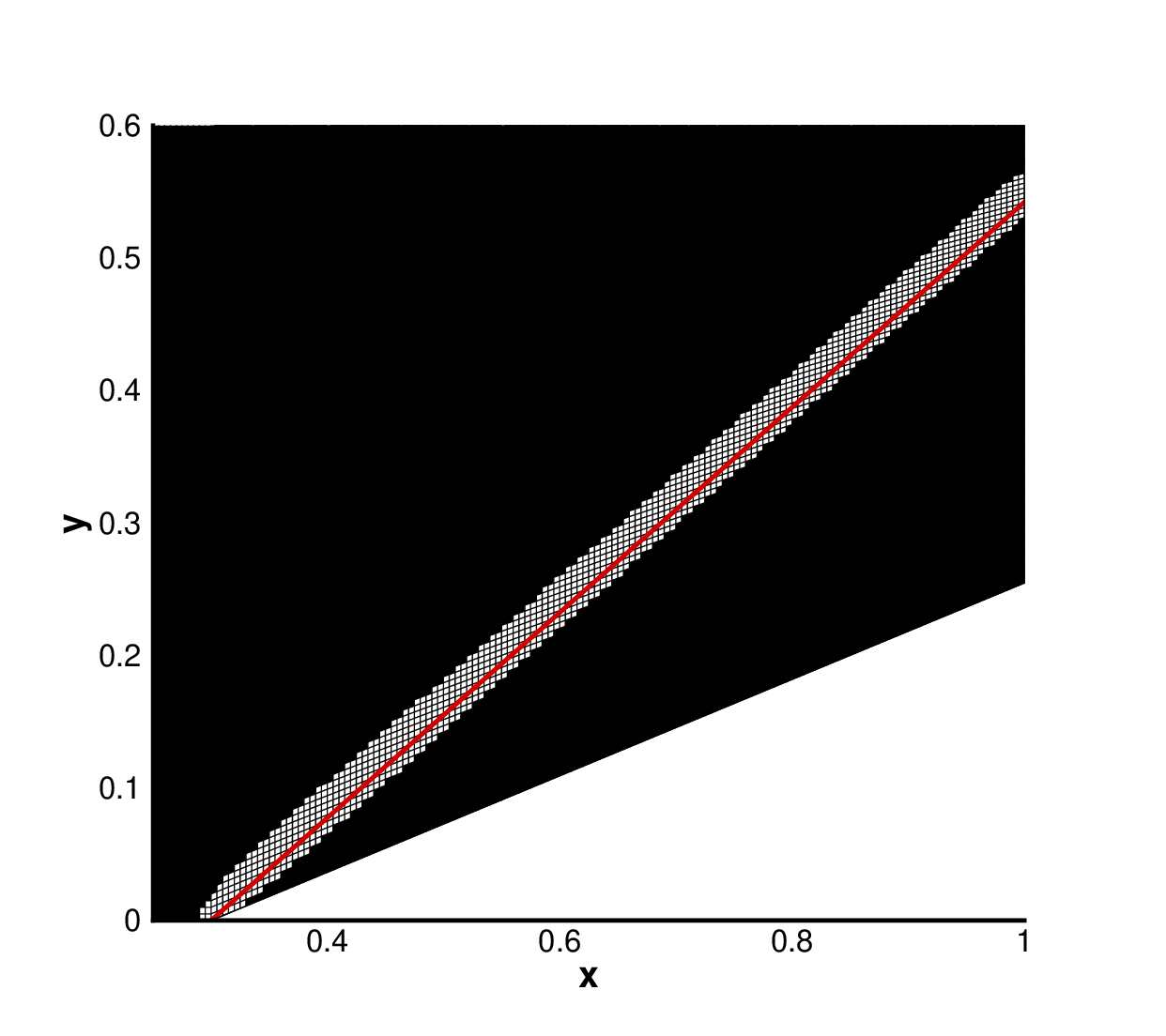}\label{fig:R20_K01_TC}}
\subfloat[$\ang{30}, K = 0.1$]{\includegraphics[width=0.5\textwidth, height=0.5\textwidth]{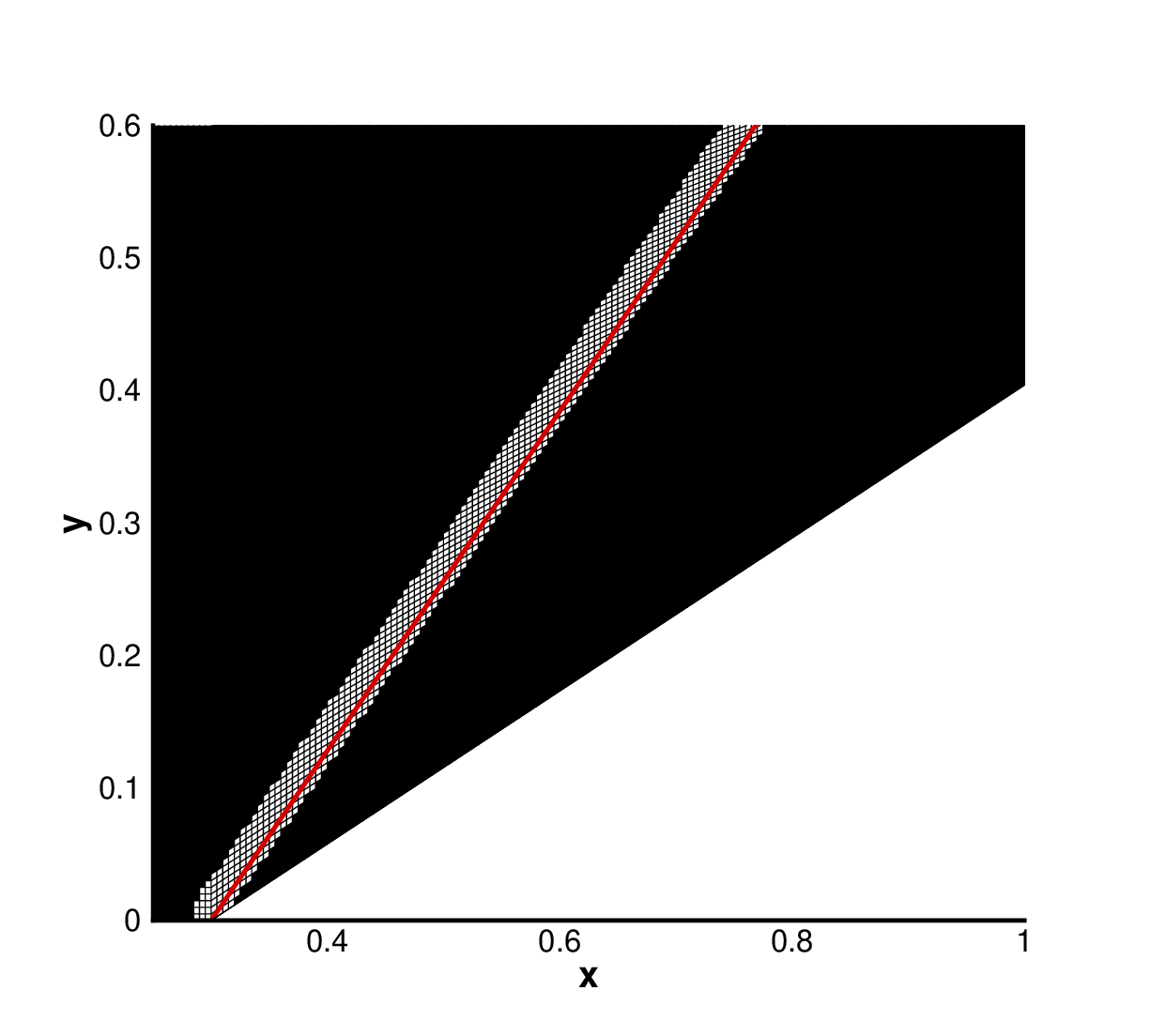}\label{fig:R30_K01_TC}}\\
\caption{Flow over a ramp. Zoomed-in view of troubled-cells identified by the indicator for two threshold constants. Red line represents the exact shock. Percentage of number of troubled-cells identified in whole computational domain of 200 $\times$ 200 cells: (a) 5.05\% (b) 6.74\% (c) 3.12\% (d) 4.45\%}
\label{fig:Ramp_TC}
\end{figure}

Figure ({\ref{fig:Ramp_TC}}) presents the zoomed-in view of the troubled-cells identified for threshold constants $K = 0.02$ and $K = 0.1$ for the test case of flow over a ramp. Similar to the observations from the non-aligned oblique shock test case, a very few number of troubled-cells are observed in the post-shock region for the constant $K = 0.1$.

\begin{figure}
\centering
\subfloat[$K = 0.02$]{\includegraphics[height=0.5\textwidth]{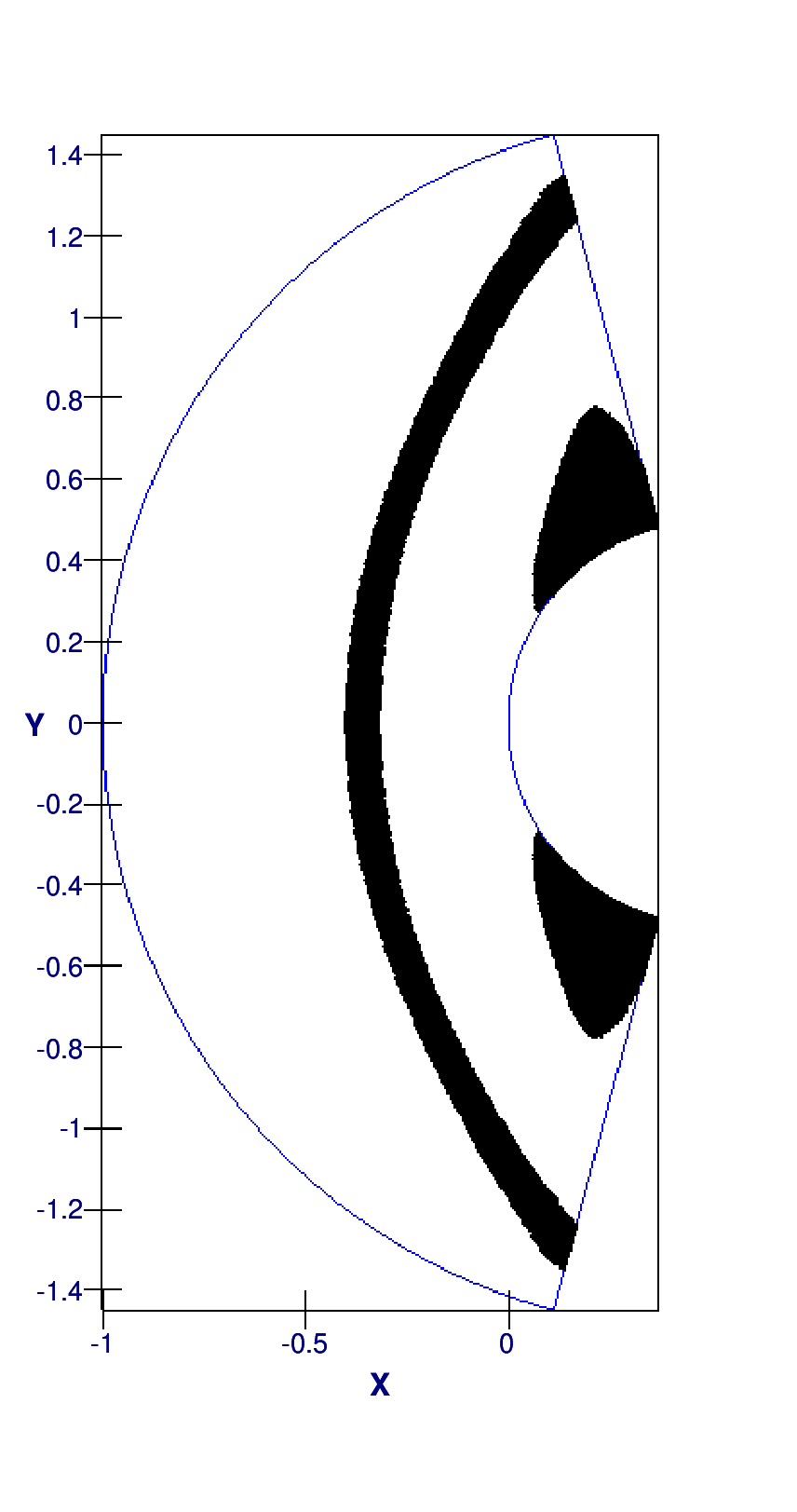}\label{fig:TC_K002}} \hspace{0.3cm}
\subfloat[$K = 0.1$]{\includegraphics[height=0.5\textwidth]{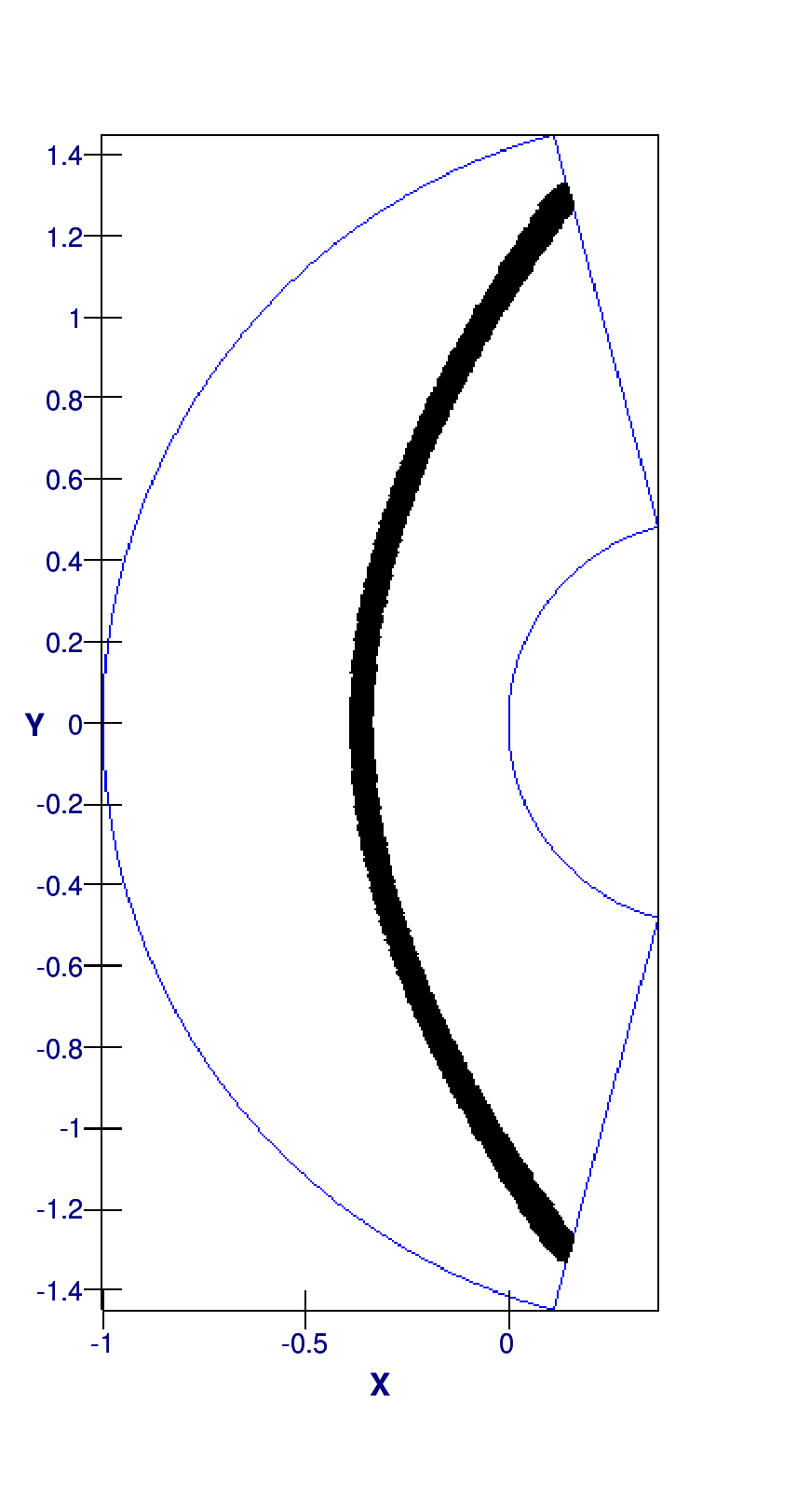}\label{fig:TC_K01}}
\caption{Flow over a cylinder. Troubled-cell region (in black) identified by the indicator for threshold constant $K = 0.02$ and $K = 0.1$, respectively. Percentage of number of troubled-cells identified in whole computational domain of 160 $\times$ 320 cells: (b) 18.26\% (c) 5.86\%}
\label{fig:BB_TC}
\end{figure}

For the flow over a cylinder test case, the troubled-cell indicator effectively identified cells in the vicinity of the shock for all three threshold constants. However, for smaller threshold constants, the region near the cylinder is also identified as troubled. The regions identified by the indicator as troubled for threshold constants $K = 0.02$ and $K = 0.1$ are shown in Figure (\ref{fig:BB_TC}).

Overall, the troubled-cell indicator demonstrated its effectiveness in identifying troubled cells for all steady-state test cases (Test cases 1 to 5). From all these figures of zoomed-in view of the troubled-cells, it can be observed that a greater number of troubled cells are identified in the pre-shock region compared to the post-shock region. This disparity can be attributed to the dissipative nature of the first-order solution obtained using the Lax-Friedrichs flux scheme.

The results indicate a decreasing trend in the number of troubled cells identified as the threshold constant increases. For certain test cases, the threshold constant $K = 0.02$ tend to identify cells located far from the shock region. Conversely, in other cases, the threshold constant $K = 0.1$ fails to adequately capture the shock region. Nevertheless, the threshold constants within this range, particularly $K = 0.05$, perform effectively and provide satisfactory results.

\section{Limiting Strategies}
\label{sec:Limiting}
Once the troubled-cells are identified, we apply the limiter function only in these cells. We refer to this approach of limiting only in troubled-cells as the limiting restricted region approach. We compare the solutions of the limiting restricted region approach with the solutions of the base solver, which we refer to as the limiting everywhere approach.

\subsection{Comparison techniques:}
We compare the solutions of these two limiting approaches using the following: \\ \\
Test of convergence to steady-state: To determine if a steady-state solution has been achieved, a residual norm, denoted by $RN$, is calculated at each iteration. This norm is defined as:
\begin{equation}\label{eq:RN}
 RN = \sqrt{\sum_{i=1}^{N} \sum_{j=1}^{4} \textbf{R}_{ij}^2(\textbf{Q}) \, \Omega_i}
\end{equation}
where, $N$ represents the number of cells in the domain, $\Omega_i$ is the volume of the $i$-th cell and $\textbf{R}_{ij}(\textbf{Q})$ is the $j$-th component of residue for that cell and is given in equation (\ref{eq:residue}). The convergence criterion is set as $10^{-14}$. The simulation is stopped when the convergence criterion is met or at the end of 15,000 iterations.\\ \\
Line plots: Density profiles of the numerical solution along a line crossing the shock are visualized to assess the accuracy of the solver in capturing the shock. \\ \\
Error norms: To evaluate and quantify the accuracy of the solution, we calculate the $L_2$ and $L_{\infty}$ norms of error in density along the line on which density is plotted. These are defined as follows:
\begin{equation}\label{eq:L_2}
 L_2(e) = \sqrt{\frac{1}{N} \sum_{i=1}^{N} |e_i|^2}
\end{equation}
\begin{equation}\label{eq:L_inf}
 L_{\infty}(e) = \max_i |e_i|
\end{equation}
Here, the error, $e$, is defined as
\begin{equation}\label{eq:error}
 e = \rho_{\text{exact}} - \rho_{\text{numerical}}
\end{equation}
where, $\rho_{\text{exact}}$ and $\rho_{\text{numerical}}$ represent the densities of the exact and numerical solutions, respectively.

However, just these norms are not enough to capture the accuracy of the solution, particularly for the case of limiting restricted region approach. We will explain this in detail with example solutions of aligned shock.

We solve the two-dimensional Euler equations (\ref{eq:2d_euler}) for aligned oblique shock configuration shown in Figure (\ref{fig:AlignOS_setup}). The computational domain is initialized with the exact solution corresponding to an inlet Mach number of 3 and a shock wave angle of \ang{40}. Boundary conditions are applied as shown in Figure (\ref{fig:AlignOS_setup}).

\begin{figure}
\centering
\subfloat[Everywhere]{\includegraphics[width=0.5\textwidth]{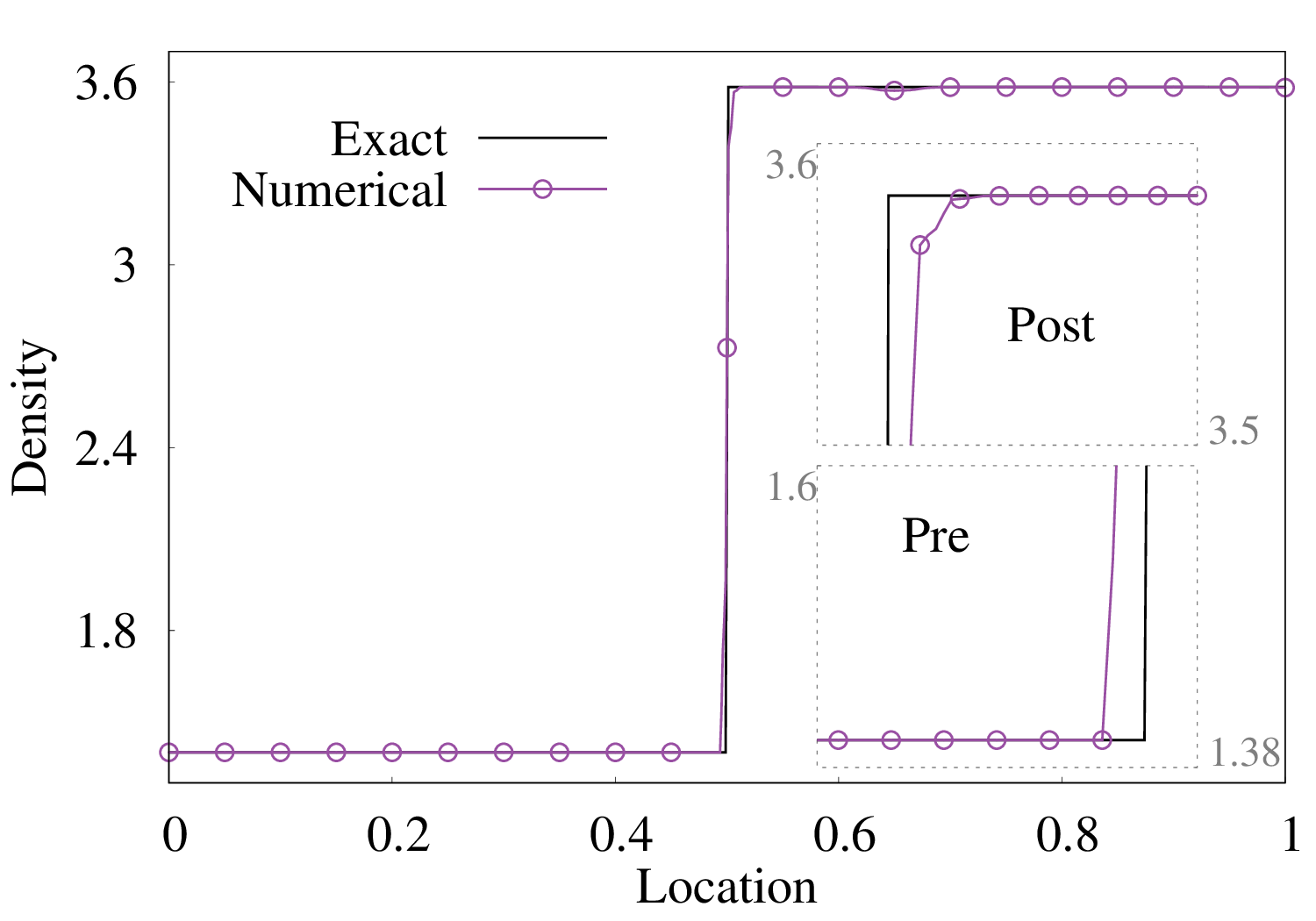}\label{fig:AOS_40_every}}
\subfloat[Restricted]{\includegraphics[width=0.5\textwidth]{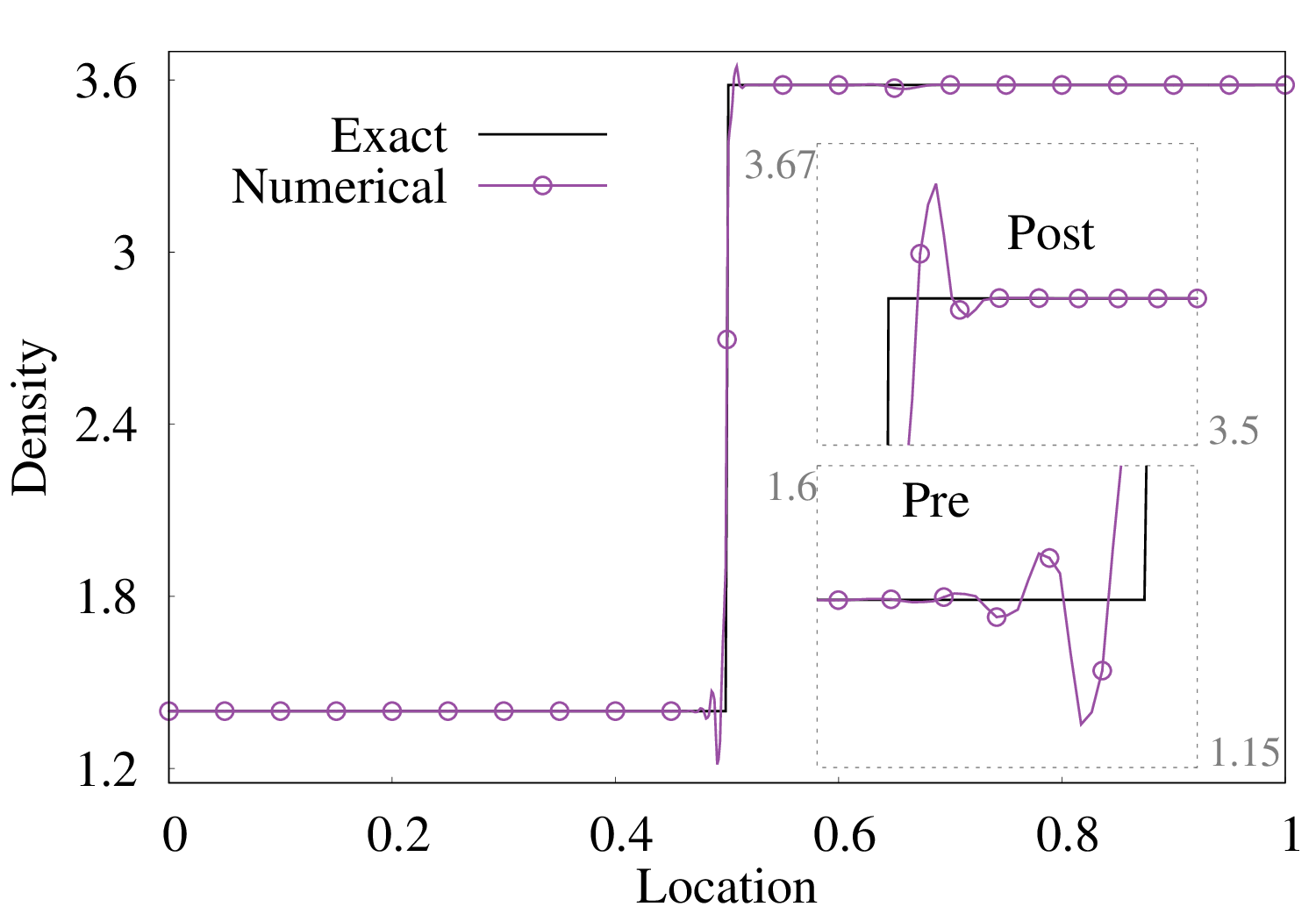}\label{fig:AOS_40_11}}
\caption{Density profiles along the line $y = 0.5$ for both limiting approaches.}
\label{fig:Density_Profiles_40}
\end{figure}

Figure (\ref{fig:Density_Profiles_40}) presents density profiles along the line $y = 0.5$ for both limiting approaches. Figure (\ref{fig:AOS_40_every}) illustrates the solution obtained using the limiting everywhere approach, while Figure (\ref{fig:AOS_40_11}) shows the solution obtained using the limiting restricted region approach. In the latter case, the limiter is applied in just one cell in the pre-shock region and one cell in the post-shock region as an example for this particular problem. A minor trough observed in the post-shock region is attributed to a small artifact generated at the shock emergence boundary due to complicated numerical boundary conditions \cite{ZHANG2013}. In keeping with our decision to use an unaugmented base solver and given that the primary focus is on the near vicinity of the shock, we will not delve deeper into this artifact in this paper.

To focus on the vicinity of the shock, the computations of the $L_2$ and $L_{\infty}$ norms are confined to a restricted domain: 20 cells in pre-shock and 20 cells in post-shock from the shock position. By restricting the domain, we also avoid the artifact mentioned earlier, ensuring a clearer evaluation of the solution near the shock. Table (\ref{tab:Solution_Compare_norms}) presents the $L_2$ and $L_{\infty}$ norms of the error in density for the two solutions shown in Figure (\ref{fig:Density_Profiles_40}). The $L_2$ and $L_{\infty}$ norms for both solutions are quite comparable, in fact slightly better for the solution of limiting restricted region approach. This can lead to the conclusion that the solution to the limiting restricted region is ``better" than limiting everywhere despite the presence of significant oscillations (overshoots and undershoots) in the solution.

\begin{table}
\centering
\caption{$L_2$ and $L_{\infty}$ norms of the density error for the solutions of both limiting approaches.}
\begin{tabular}{c c c}
\toprule
\makecell{Limiting} & $L_2$ & $L_{\infty}$ \\
\midrule
\makecell{everywhere \\ (Figure (\ref{fig:AOS_40_every}))} & 0.069191 & 0.561910\\
\midrule
\makecell{restricted \\ (Figure (\ref{fig:AOS_40_11}))} & 0.063761 & 0.498884\\
\bottomrule
\end{tabular}%
\label{tab:Solution_Compare_norms}
\end{table}

If the limiter function is not applied in a sufficient number of troubled-cells in the vicinity of the shock, then there will be overshoots and undershoots in the solution as shown in Figure (\ref{fig:AOS_40_11}). These oscillations can be captured by the total variation since it is a measure of how oscillatory a solution is. We compute the total variation of the error, $e$ (Equation (\ref{eq:error})), in density. This is defined as follows:
\begin{equation}\label{eq:TV}
TV(e) = \sum_{i=1}^{N-1} |e_{i+1} - e_i|
\end{equation}

\textbf{Monotonic solution property:} If the solution monotonically increases (or decreases) from the exact solution in the pre-shock region, then the total variation of the error and $L_{\infty}$ norm values are equal in that region. Similarly, If the solution monotonically increases (or decreases) to the exact solution in post-shock region, then the total variation of the error and $L_{\infty}$ norm values are equal in that region.

We use this property of monotonic solutions to quantify the quality of the solutions obtained from two limiting approaches and compare them. For that, we compute the total variation of the error in density and $L_{\infty}$ norm of the error in density separately for both pre- and post-shock regions. 20 cells in the pre-shock region from the shock position is the domain for the pre-shock region and 20 cells in the post-shock region from the shock position is the domain for the post-shock region.

\begin{table}
\centering
\caption{$L_{\infty}$ norm and the total variation of the density error for both limiting approaches.}
\begin{tabular}{c c c c}
\toprule
\makecell{Region} & \makecell{Limiting} & $L_{\infty}$ & $TV$  \\
\midrule
\multirow{2}{*}{Pre-shock} & \makecell{everywhere \\ (Figure (\ref{fig:AOS_40_every}))} & 0.561910 & 0.561910 \\ \cline{2-4}
& \makecell{restricted \\ (Figure (\ref{fig:AOS_40_11}))} & 0.498884 & 1.090482 \\
\midrule
\multirow{2}{*}{Post-shock} & \makecell{everywhere \\ (Figure (\ref{fig:AOS_40_every}))} & 0.199483 & 0.199641\\ \cline{2-4}
& \makecell{restricted \\ (Figure (\ref{fig:AOS_40_11}))} & 0.204351 & 0.355214\\
\bottomrule
\end{tabular}%
\label{tab:Solution_Compare_TV}
\end{table}

Table (\ref{tab:Solution_Compare_TV}) presents the $L_{\infty}$ norm and the total variation of the density error for two solutions shown in Figure (\ref{fig:Density_Profiles_40}). For the solution of limiting everywhere approach (Figure (\ref{fig:AOS_40_every})), the $L_{\infty}$ norm is exactly equal to the total variation of the error in the pre-shock region, whereas the total variation is slightly higher in the post-shock region. This indicates the minimal presence of oscillations in the solution, and the solution is  monotonically increasing. Conversely, for the solution of limiting restricted region approach (Figure (\ref{fig:AOS_40_11})), the total variation of the error is significantly higher than the $L_{\infty}$ norm in both the pre- and post-shock regions, indicating the presence of substantial oscillations in the solution.

In the next section, we quantify the quality of the solution in the neighbourhood of the shock for both limiting approaches for various test cases, particularly for the test cases to which the exact solution is known, through the following computations:
\begin{itemize}
 \item Calculate the $L_{\infty}$ norm of the error in density separately for the pre-shock and post-shock regions, each confined to 20 cells from the shock location.
 \item Similarly, evaluate the total variation of the error in density for the pre-shock and post-shock regions.
\item Add the pre- and post-shock values to obtain the corresponding overall value for the entire domain of 40 cells.
\item We define a parameter, $\mu$, called the monotonicity parameter, as
\begin{equation}
 \mu = TV - L_{\infty}
\end{equation}
\item Compute the monotonicity parameter, $\mu$, using the overall total variation and the overall $L_{\infty}$ norm of the error in density.
\end{itemize}
The better solution will be the one with the value of monotonicity parameter $\mu$ is close to zero and lower in comparison to those of other solutions.

\subsection{Results of two limiting approaches}
For all steady-state test cases (Test cases 1 to 5), troubled-cells are identified using the troubled-cell indicator applied to a first-order solution and shown in the previous section. To obtain the high-order solution, the computational domain is initialized using the first-order solution, incorporating information about the identified troubled cells. Throughout the simulation, we keep the troubled cells fixed i.e., limiting in the same troubled-cells at every iteration. The limiter is applied in all computational cells for the case of the limiting-everywhere approach, whereas it is applied only in the troubled cells for the case of the limiting restricted region approach.

\begin{figure}
\centering
\subfloat[Aligned Shock - $\ang{30}$]{\includegraphics[width=0.45\linewidth]{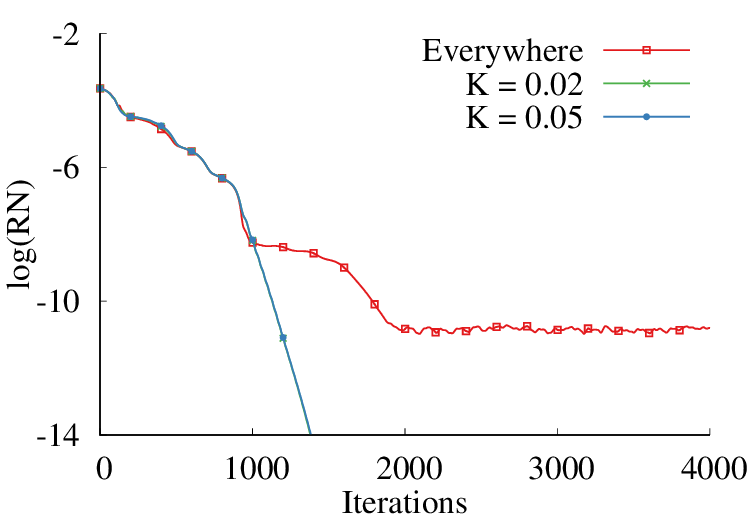}\label{fig:AOS_30_RN}}\hspace{0.2cm}
\subfloat[Aligned Shock - $\ang{40}$]{\includegraphics[width=0.45\linewidth]{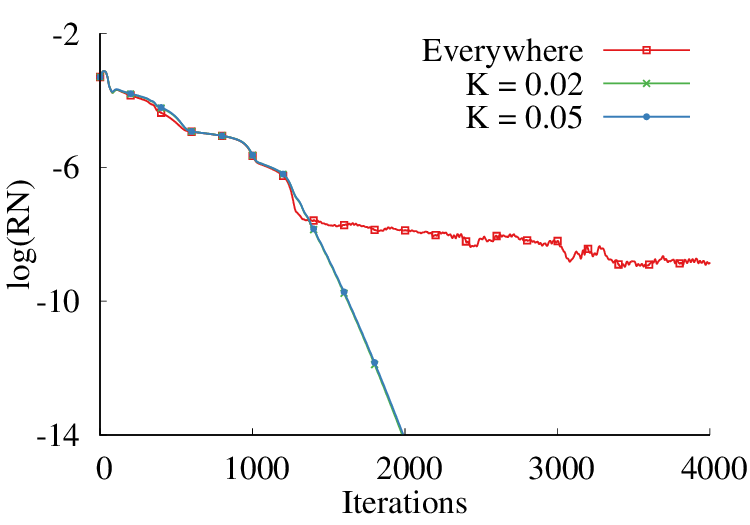}\label{fig:AOS_40_RN}}\\
\subfloat[Non-aligned Shock - $\ang{30}$]{\includegraphics[width=0.45\linewidth]{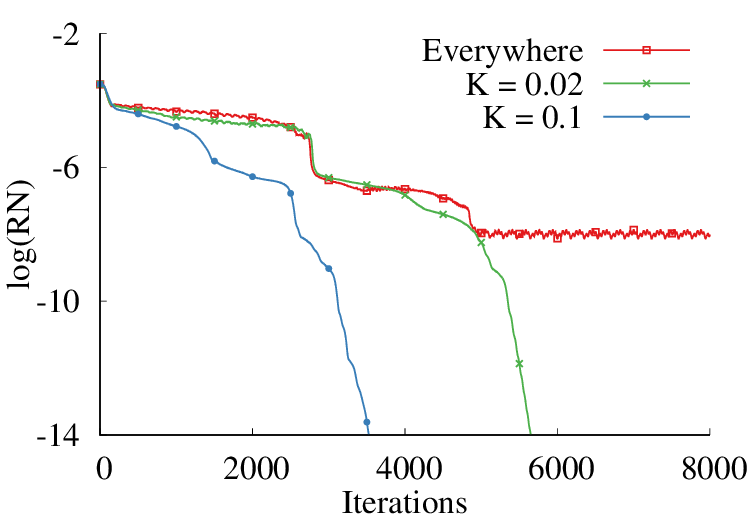}\label{fig:NOS_30_RN}}\hspace{0.2cm}
\subfloat[Non-aligned Shock - $\ang{40}$]{\includegraphics[width=0.45\linewidth]{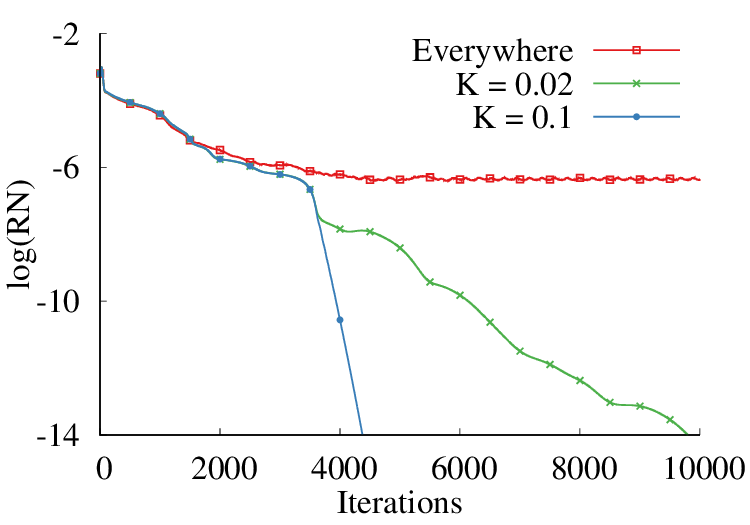}\label{fig:NOS_40_RN}}\\
\subfloat[Ramp - $\ang{20}$]{\includegraphics[width=0.45\linewidth]{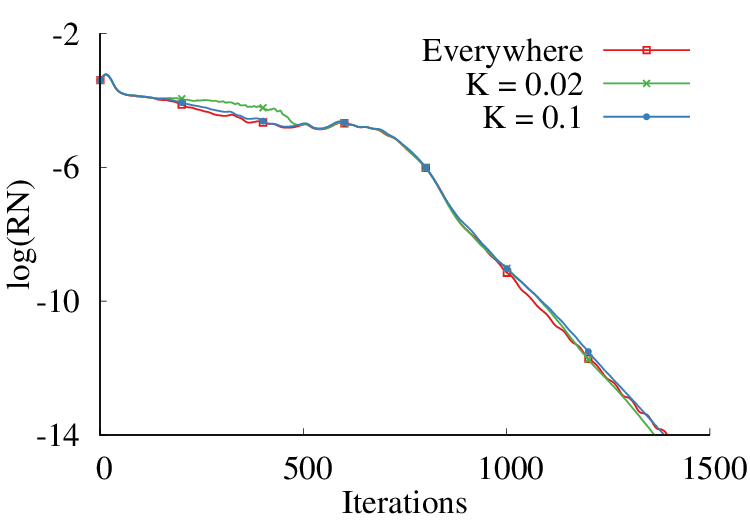}\label{fig:Ramp_20_RN}}\hspace{0.2cm}
\subfloat[Ramp - $\ang{30}$]{\includegraphics[width=0.45\linewidth]{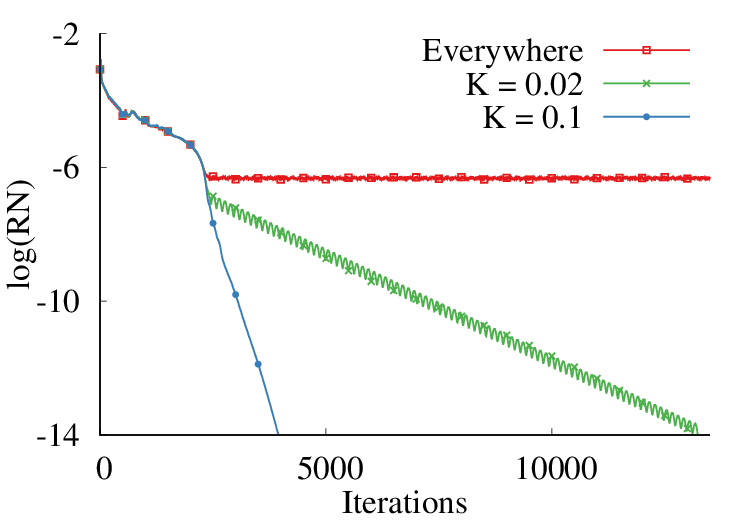}\label{fig:Ramp_30_RN}}\\
\subfloat[Shock Reflection]{\includegraphics[width=0.45\linewidth]{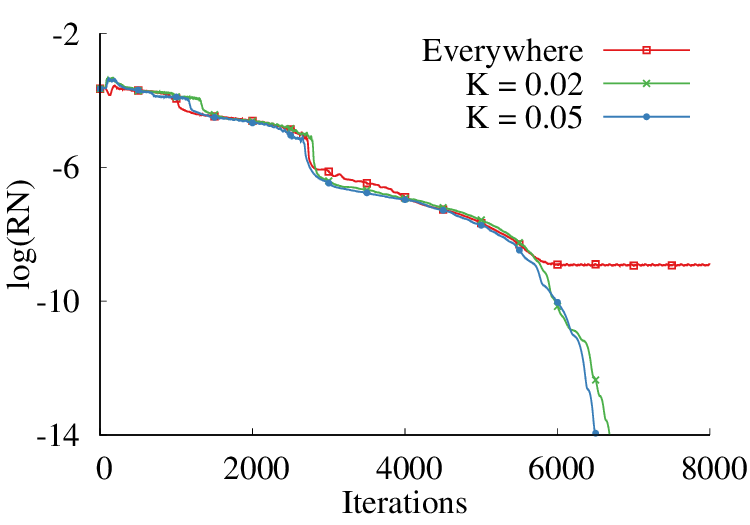}\label{fig:SR_RN}}\hspace{0.2cm}
\subfloat[Blunt Body]{\includegraphics[width=0.45\linewidth]{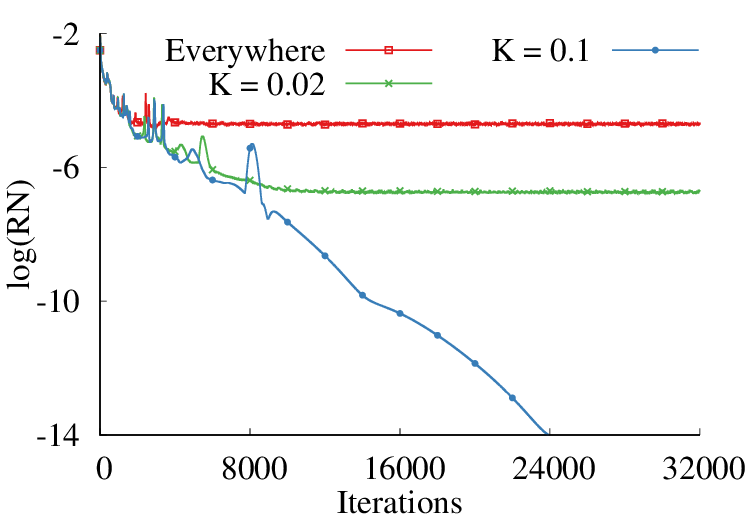}\label{fig:BB_RN}}
\caption{The convergence history of the residual norm as a function of number of iterations for both limiting approaches for all steady-state test cases (Test cases 1 to 5).}
\label{fig:RN_Forall}
\end{figure}

Figure (\ref{fig:RN_Forall}) presents the convergence history of the residual norm, $RN$ (Eq. \ref{eq:RN}), as a function of the number of iterations for both limiting approaches for all steady-state test cases (Test cases 1 to 5). In this figure, the labels `K = 0.02', `K = 0.05' and `K = 0.1' indicate that the limiter is applied only in the troubled-cells identified by the indicator for thresholds constants 0.02, 0.05, and 0.1 respectively. And, the label `Everywhere' indicates the application of limiter in all computational cells. These label definitions hold consistently for all other figures as well.

The results indicate that the convergence of the limiting restricted region approach is significantly better than that of the limiting everywhere approach across all steady-state test cases. Furthermore, within the limiting restricted region approach, as the threshold constant increases, resulting in fewer identified troubled cells, the convergence of the residual norm reaches the convergence criterion ($10^{-14}$) in fewer iterations compared to smaller threshold constants.

\begin{figure}
\centering
\subfloat[Aligned Shock - $\ang{30}$]{\includegraphics[width=0.45\linewidth]{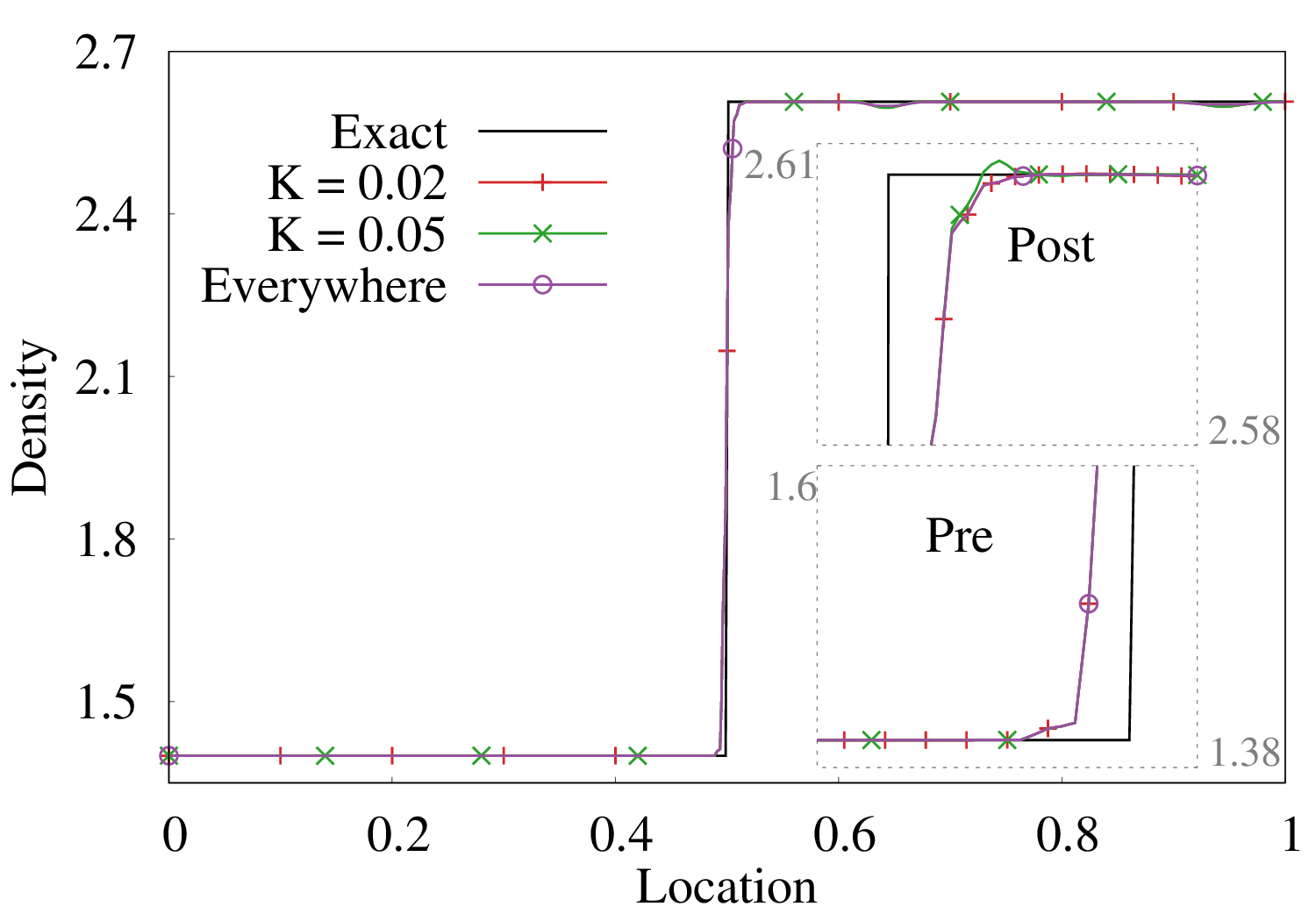}\label{fig:AOS_30_Density}}\hspace{0.2cm}
\subfloat[Aligned Shock - $\ang{40}$]{\includegraphics[width=0.45\linewidth]{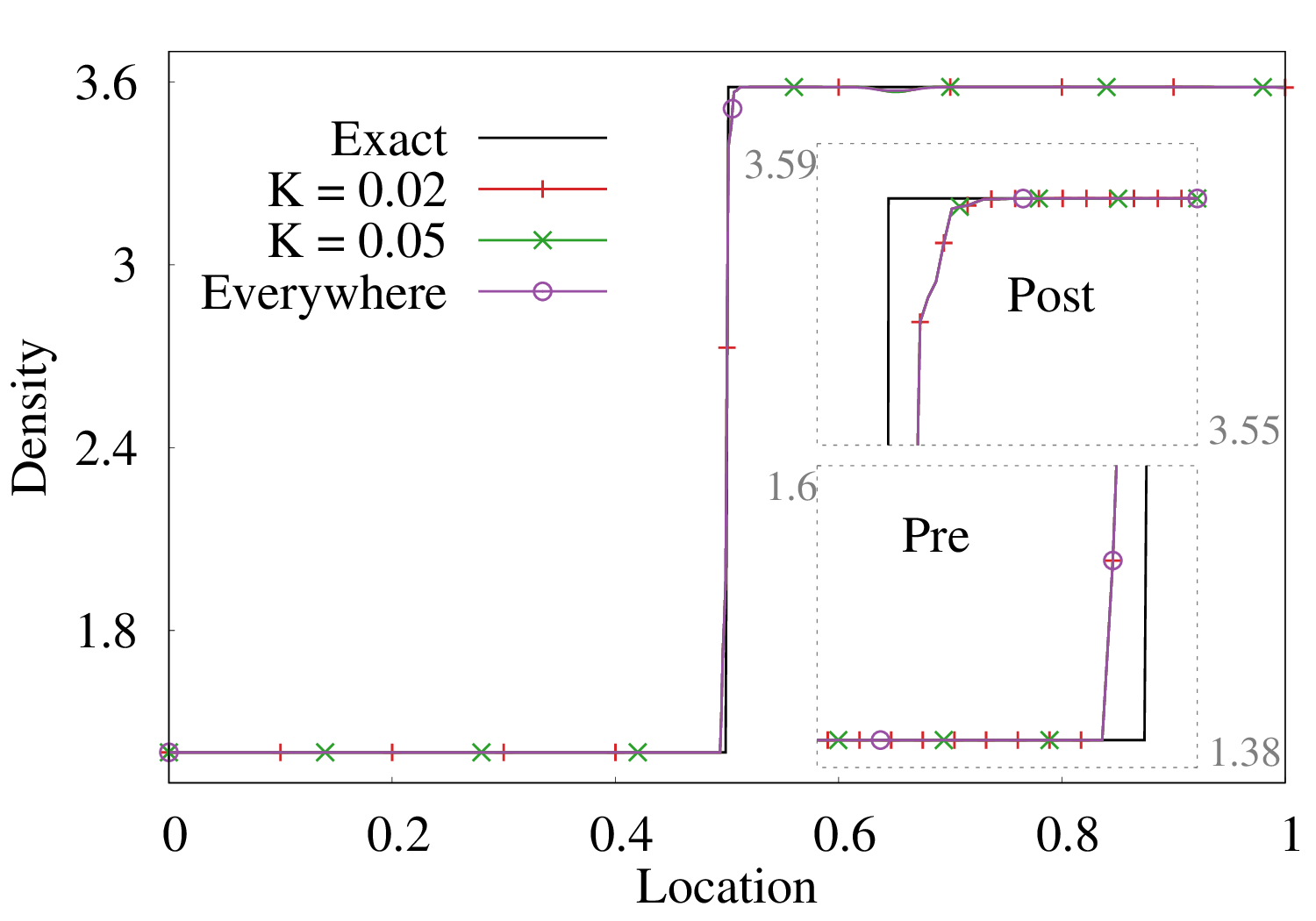}\label{fig:AOS_40_Density}}\\
\subfloat[Non-aligned Shock - $\ang{30}$]{\includegraphics[width=0.45\linewidth]{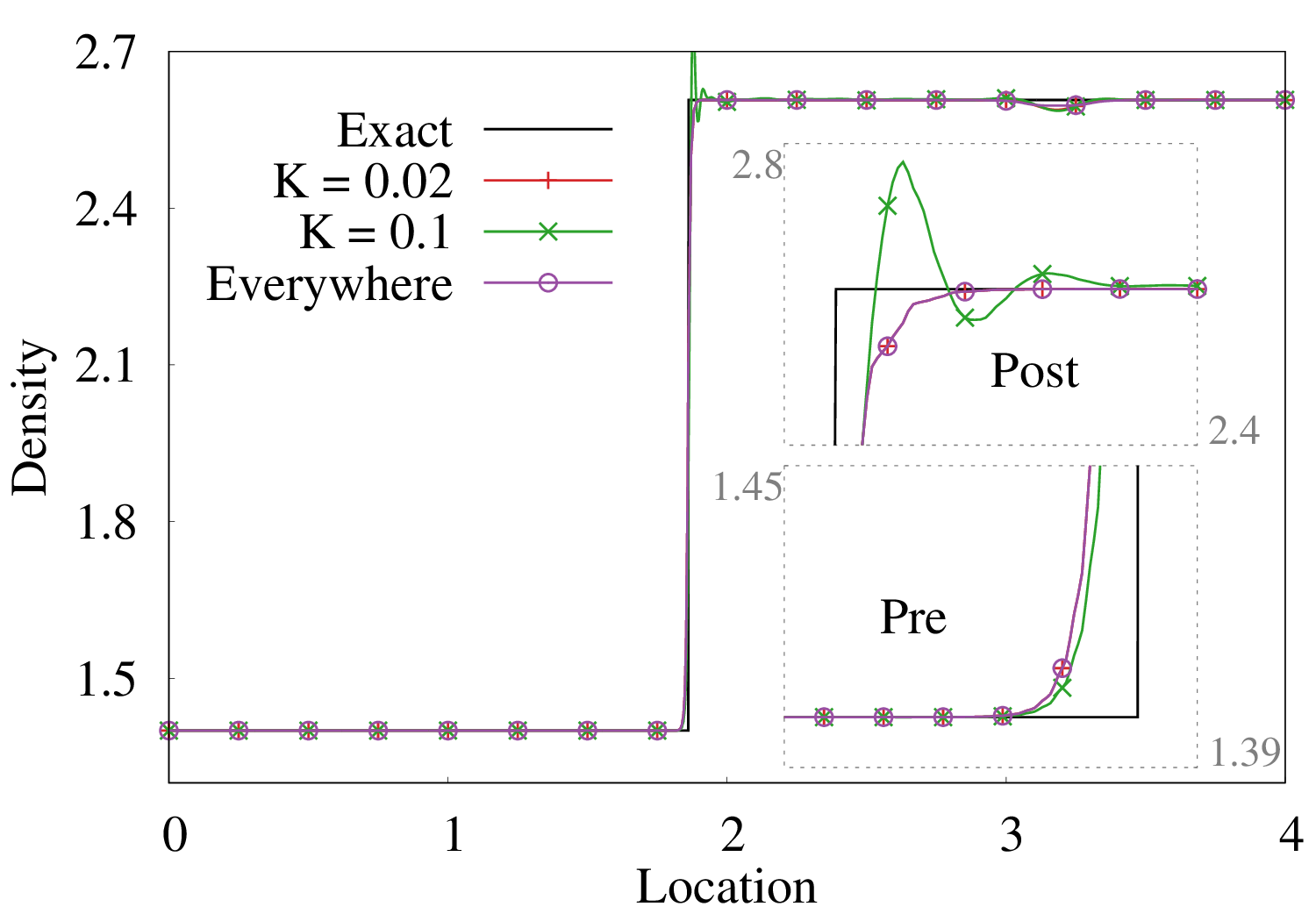}\label{fig:NOS_30_Density}}\hspace{0.2cm}
\subfloat[Non-aligned Shock - $\ang{40}$]{\includegraphics[width=0.45\linewidth]{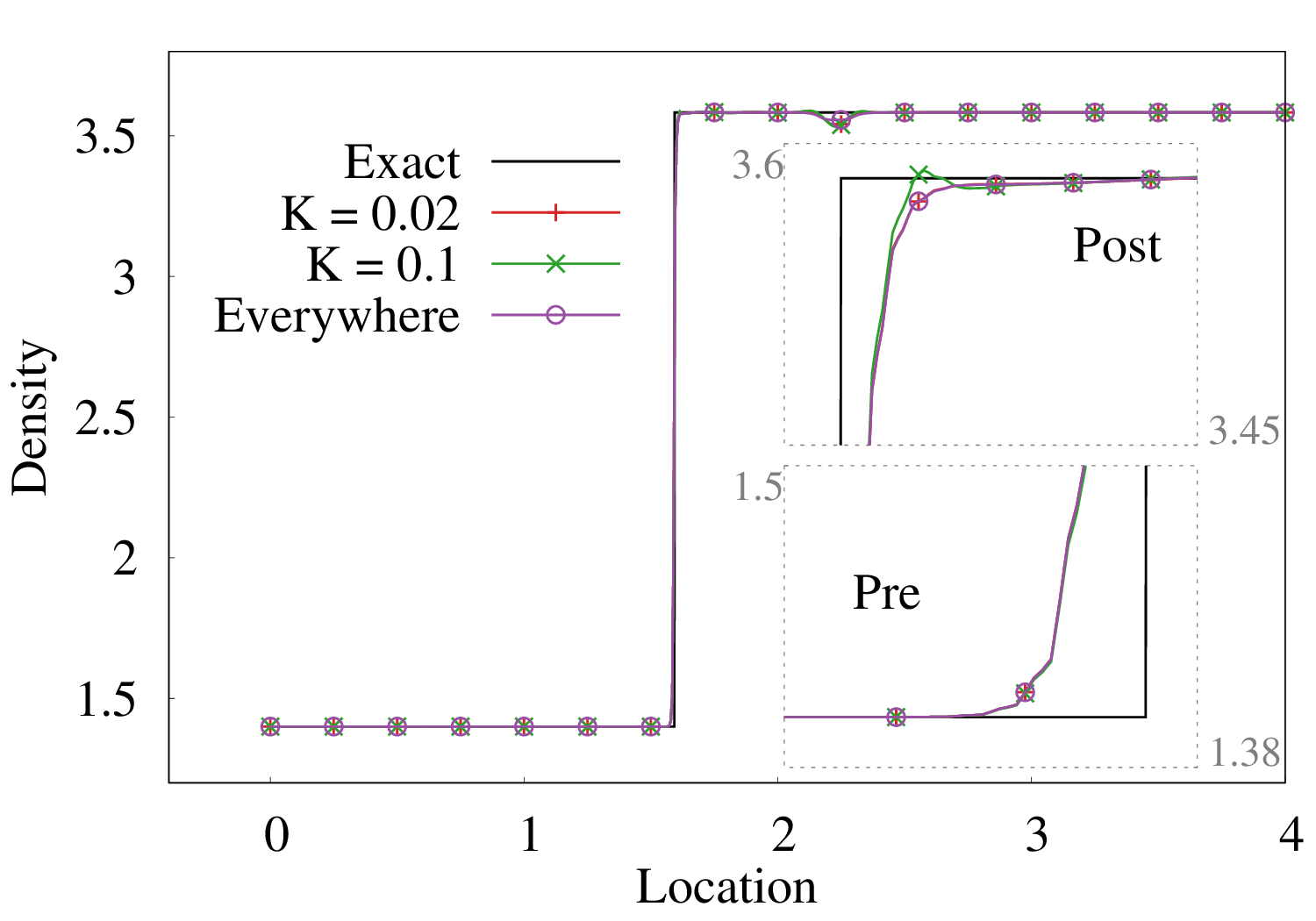}\label{fig:NOS_40_Density}}\\
\subfloat[Ramp - $\ang{20}$]{\includegraphics[width=0.45\linewidth]{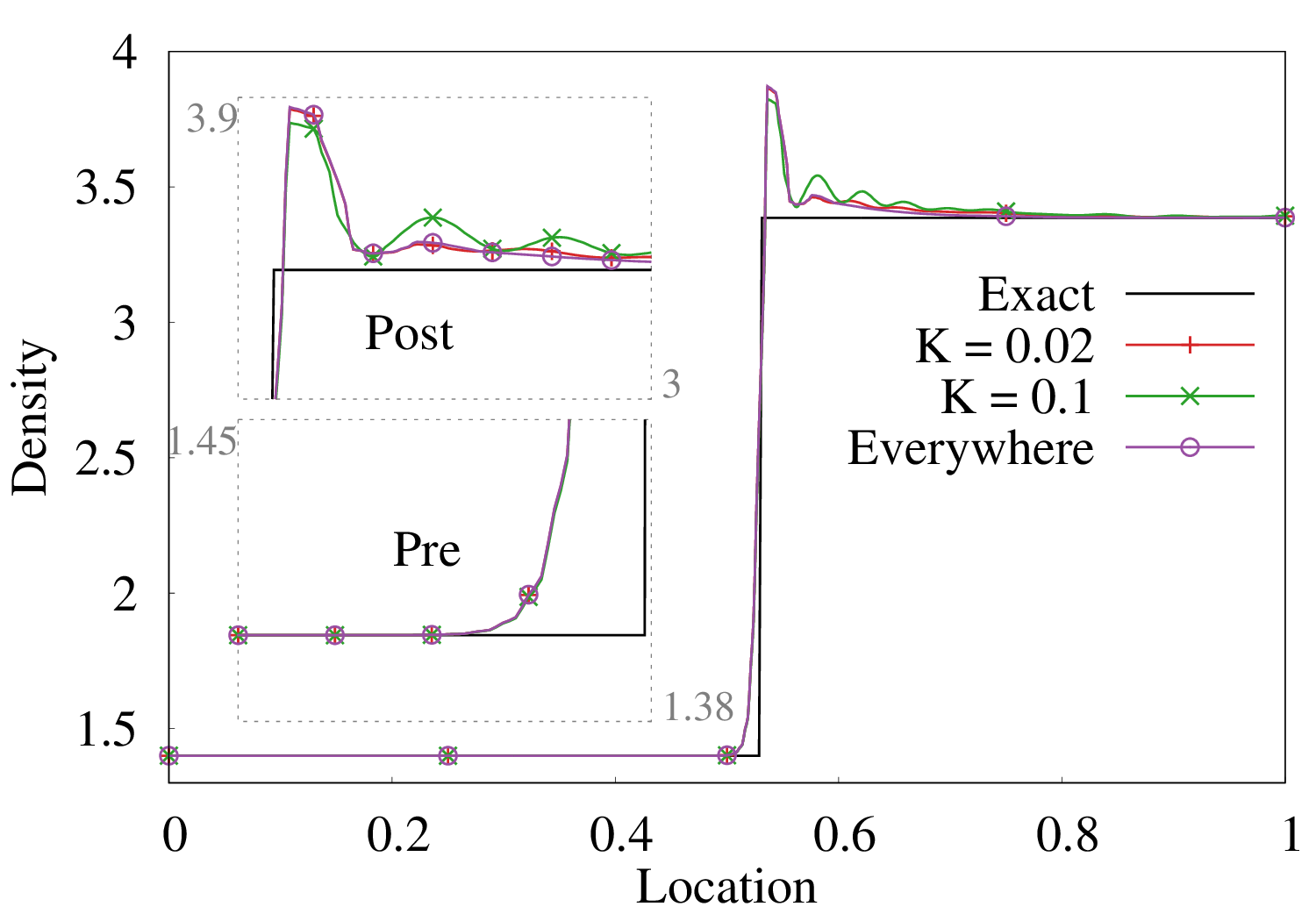}\label{fig:Ramp_20_Density}}\hspace{0.2cm}
\subfloat[Ramp - $\ang{30}$]{\includegraphics[width=0.45\linewidth]{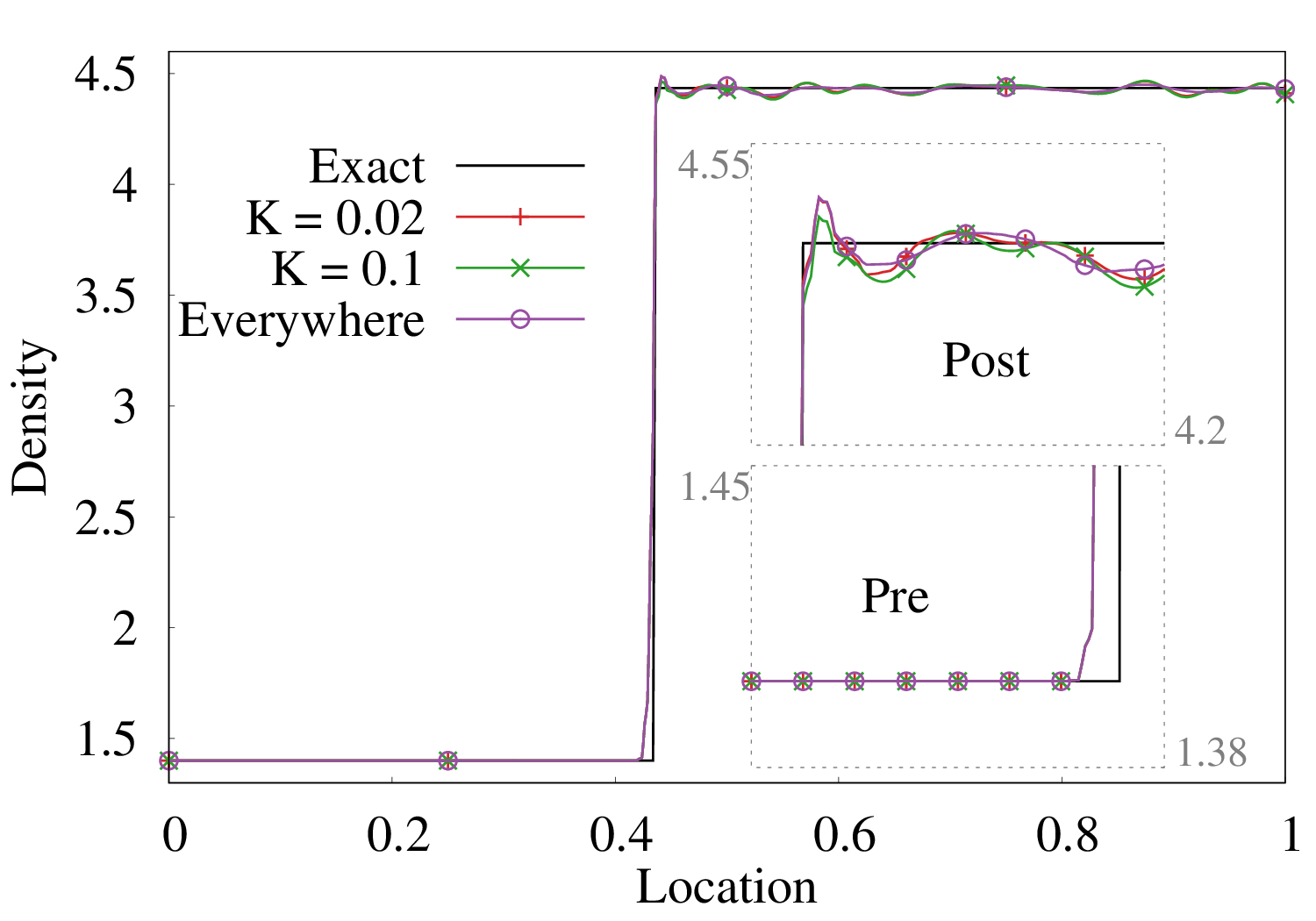}\label{fig:Ramp_30_Density}}\\
\subfloat[Shock Reflection]{\includegraphics[width=0.45\linewidth]{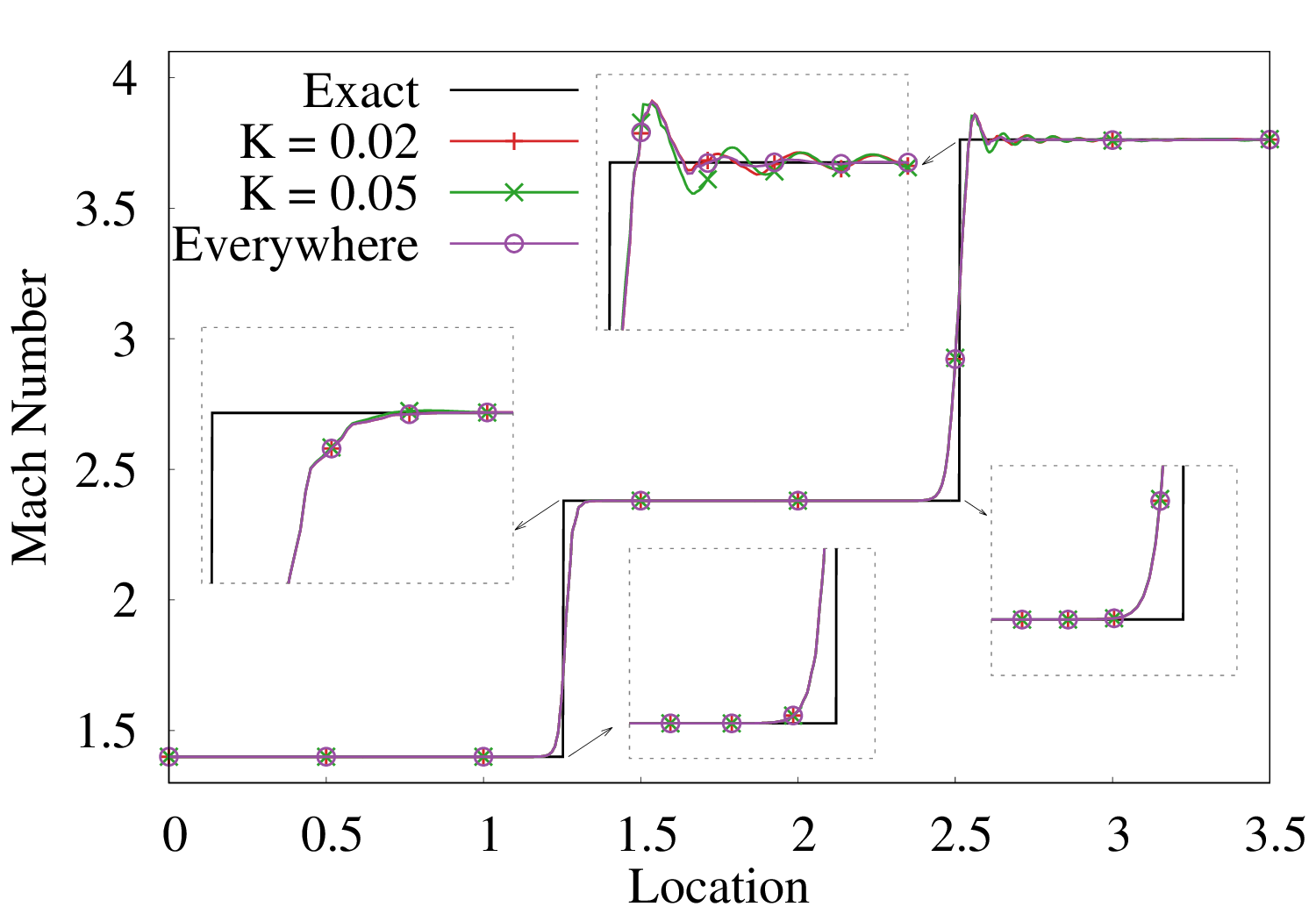}\label{fig:SR_Density}}
\caption{Density profiles along a line for both limiting approaches for test cases 1 to 4. (a), (b), (c), and (d): along the line $y = 0.5$. (e) and (f): along a streamline originating at $y = 0.1$ at the inlet. (g) along the line $y = 0.3$.}
\label{fig:LinePlots_Forall}
\end{figure}

Figure (\ref{fig:LinePlots_Forall}) depicts the density profiles along a line for both limiting approaches for steady-state test cases 1 to 4. For both aligned and non-aligned oblique shocks, the density profiles are plotted along the line $y = 0.5$. In the aligned oblique shock test case, the density profiles of both limiting approaches are almost comparable. However, a slight overshoot is observed for the case of \ang{30} with a threshold constant of $K = 0.05$, as only three troubled cells are identified in the post-shock region (Figure (\ref{fig:A30_K005_TC})) compared to six troubled cells for $K = 0.02$ (Figure (\ref{fig:A30_K002_TC})).

Similarly, in the non-aligned oblique shock test case, the density profiles of the limiting restricted region approach, particularly for $K = 0.02$, closely match those of the limiting everywhere approach. However, overshoots are evident in the solutions for $K = 0.1$. For the \ang{30} non-aligned oblique shock, these overshoots are notably significant due to the limited number of troubled cells identified in the post-shock region (Figure (\ref{fig:NA30_K01_TC})).

For the shock reflection over a flat plate test case, the density profiles are plotted along the line $y = 0.3$. For the flow over a ramp test case, the density profiles are plotted along a streamline originating at $y = 0.1$ at the inlet. Consistent with the results from other test cases, the solutions obtained using the limiting-restricted region approach with a higher $K$ exhibit more oscillations in the post-shock region compared to those with a smaller $K$.

These results can be verified from the Table (\ref{tab:mu_Forall}), where the monotonicity parameter, $\mu$, presented for both limiting approaches for test cases 1 to 4. The monotonicity parameter values are higher for the solutions obtained using the limiting-restricted region approach with a higher $K$ compared to those with a smaller $K$. For the \ang{30} non-aligned oblique shock, the $\mu$ value is significantly higher for $K = 0.1$ compared to $K = 0.02$ due to the presence of strong overshoots in the solution.

\begin{table}
\centering
\caption{Overall $L_{\infty}$ norm of the density error, the total variation of the density error and the monotonicity parameter for both limiting approaches for test cases 1 to 4.}
\begin{tabular}{c c c c c c}
\toprule
\makecell{Test case} & \makecell{Angle} & \makecell{$K$} & \makecell{$TV$}& \makecell{$L_{\infty}$} & $\mu$\\
\midrule
\multirow{6}{*}{Aligned shock} & \multirow{3}{*}{\ang{30}} & $0.02$ & 0.690730 	& 0.688894 	& 1.836e-03\\
& & $0.05$ & 0.692678 & 0.687912 	& 4.766e-03\\
 & & \cellcolor{lightgray}Everywhere & \cellcolor{lightgray} 0.689646 	& \cellcolor{lightgray} 0.688909 & \cellcolor{lightgray} 7.370e-04\\
\cmidrule(){2-6}
& \multirow{3}{*}{\ang{40}} & $0.02$ & 0.761684 	& 0.761424 	& 2.600e-04\\
& & $0.05$ & 0.761714 &	 0.761382 	& 3.320e-04\\
& & \cellcolor{lightgray} Everywhere & \cellcolor{lightgray} 0.761551 &	 \cellcolor{lightgray} 0.761393 & \cellcolor{lightgray}1.580e-04\\

\midrule

\multirow{6}{*}{Non-aligned shock} & \multirow{3}{*}{\ang{30}} & $0.02$ & 1.147028 	& 1.146915 	& 1.130e-04\\
& & $0.1$ & 1.587309 	& 1.123343 &	 4.640e-01\\
& & \cellcolor{lightgray} Everywhere & \cellcolor{lightgray} 1.147028 	& \cellcolor{lightgray} 1.146865 & \cellcolor{lightgray}1.630e-04 \\
\cmidrule(){2-6}
& \multirow{3}{*}{\ang{40}} & $0.02$ & 2.017843 	& 2.016480 	& 1.363e-03\\
& & $0.1$ & 2.033468 & 2.013740 & 1.973e-02\\
&  & \cellcolor{lightgray} Everywhere & \cellcolor{lightgray} 2.017959 & \cellcolor{lightgray} 2.016751 & \cellcolor{lightgray} 1.208e-03\\

\midrule

\multirow{6}{*}{Ramp} & \multirow{3}{*}{\ang{20}} & $0.02$ & 2.667568 	& 1.704789 & 9.628e-01\\
& & $0.1$ & 2.792442 	& 1.675887 	& 1.117e+00 \\
& & \cellcolor{lightgray} Everywhere & \cellcolor{lightgray} 2.687059 & \cellcolor{lightgray} 1.714014 & \cellcolor{lightgray} 9.730e-01 \\
\cmidrule(){2-6}
& \multirow{3}{*}{\ang{30}} & $0.02$ & 1.354126 & 0.924365 & 4.298e-01\\
& & $0.1$ & 1.420589 	& 0.914211 	& 5.064e-01\\
& & \cellcolor{lightgray} Everywhere & \cellcolor{lightgray} 1.240228 	& \cellcolor{lightgray} 0.907188 &	\cellcolor{lightgray} 3.330e-01\\

\midrule

\multirow{6}{*}{Shock Reflection} & \multirow{3}{*}{Incident} & $0.02$ & 0.938981 	& 0.938610 	& 3.710e-04\\
& & $0.05$ & 0.941201 & 0.938667 & 2.534e-03\\
&  & \cellcolor{lightgray} Everywhere & \cellcolor{lightgray} 0.938979 & \cellcolor{lightgray} 0.938619 & \cellcolor{lightgray}3.600e-04\\
\cmidrule(){2-6}
& \multirow{3}{*}{Reflected} & $0.02$ & 1.623307 & 1.332814 & 2.905e-01\\
& & $0.05$ & 1.701964 & 1.332589 & 3.694e-01\\
& & \cellcolor{lightgray} Everywhere & \cellcolor{lightgray} 1.592389 & \cellcolor{lightgray} 1.332786 & \cellcolor{lightgray}2.596e-01\\

\bottomrule
\end{tabular}%
\label{tab:mu_Forall}
\end{table}

\begin{figure}
\centering
\subfloat[$K = 0.02$]{\includegraphics[width=0.3\textwidth, height=0.5\textwidth]{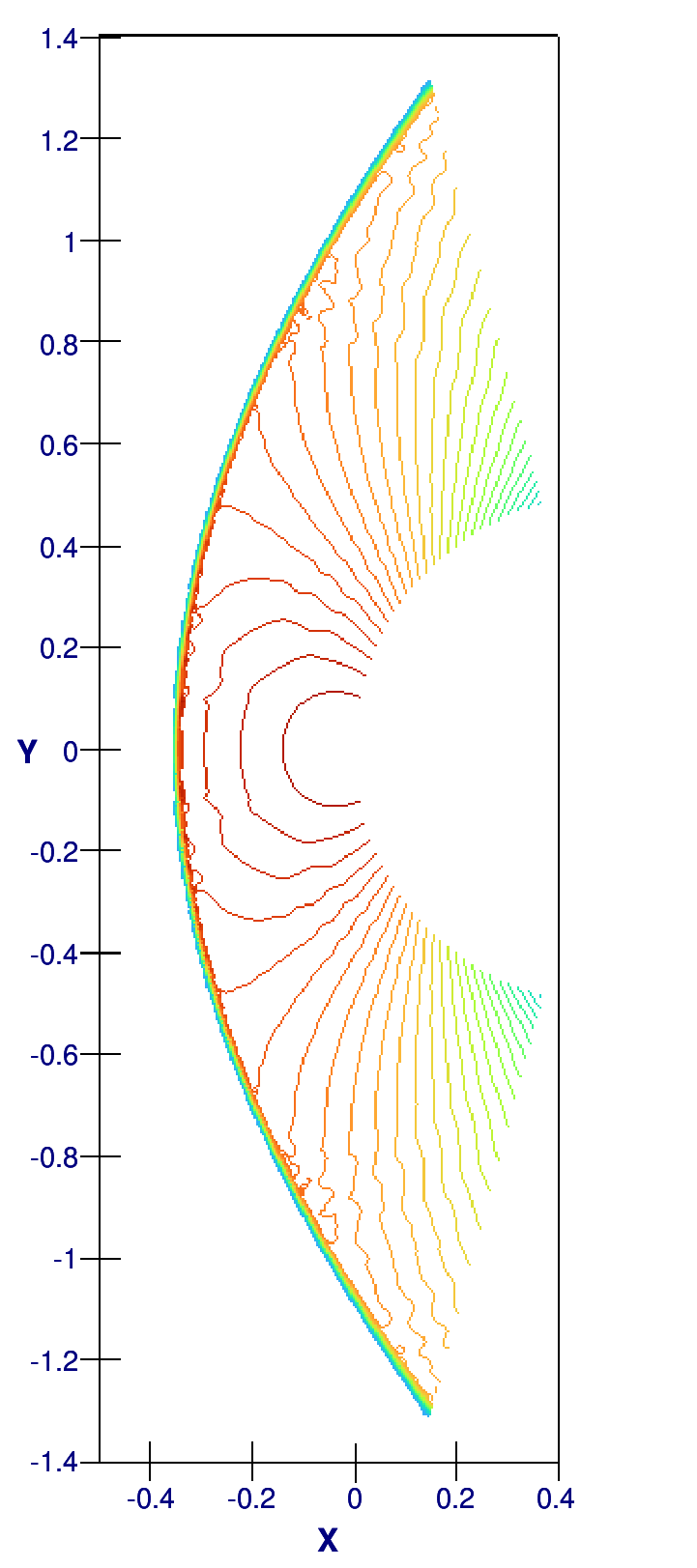}\label{fig:Con_K002}} \hspace{0.3cm}
\subfloat[$K = 0.1$]{\includegraphics[width=0.3\textwidth, height=0.5\textwidth]{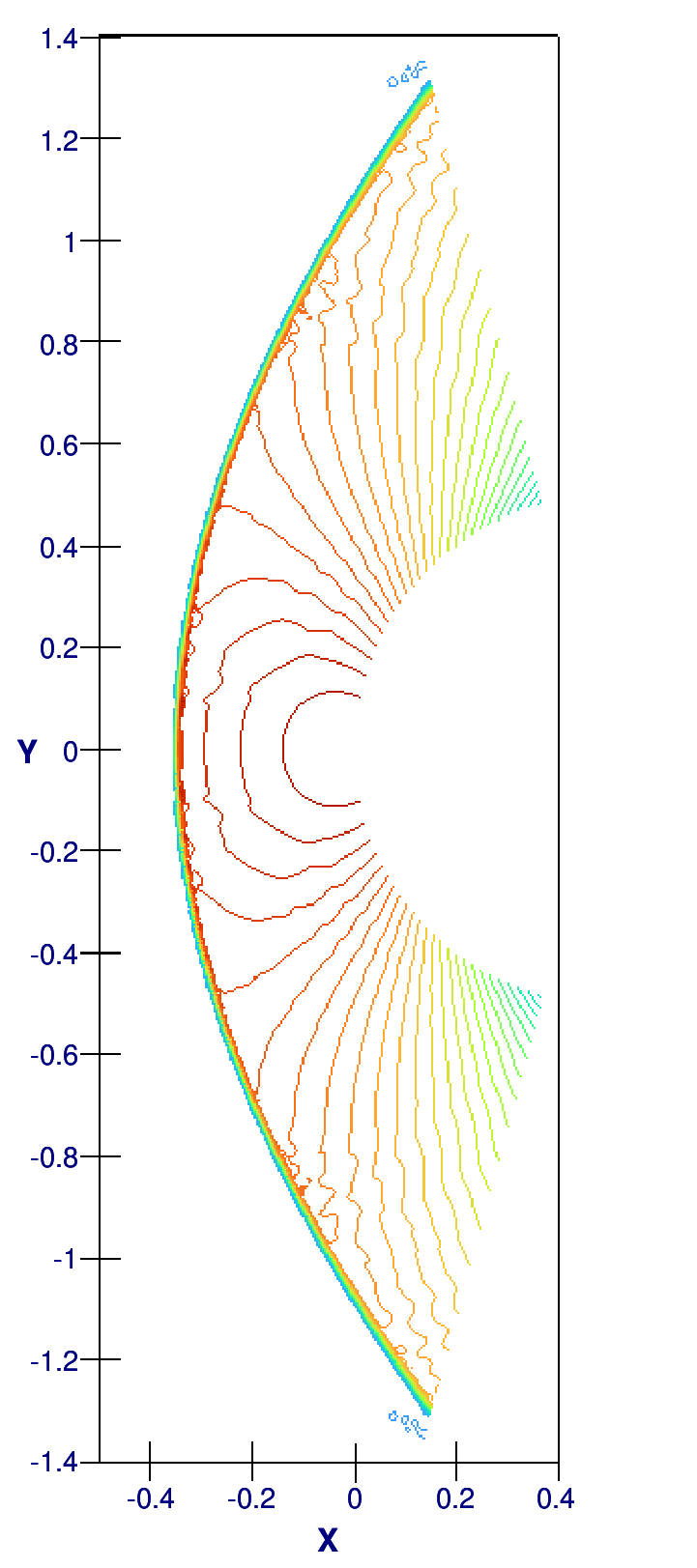}\label{fig:Con_K01}} \hspace{0.3cm}
\subfloat[Everywhere]{\includegraphics[width=0.3\textwidth, height=0.5\textwidth]{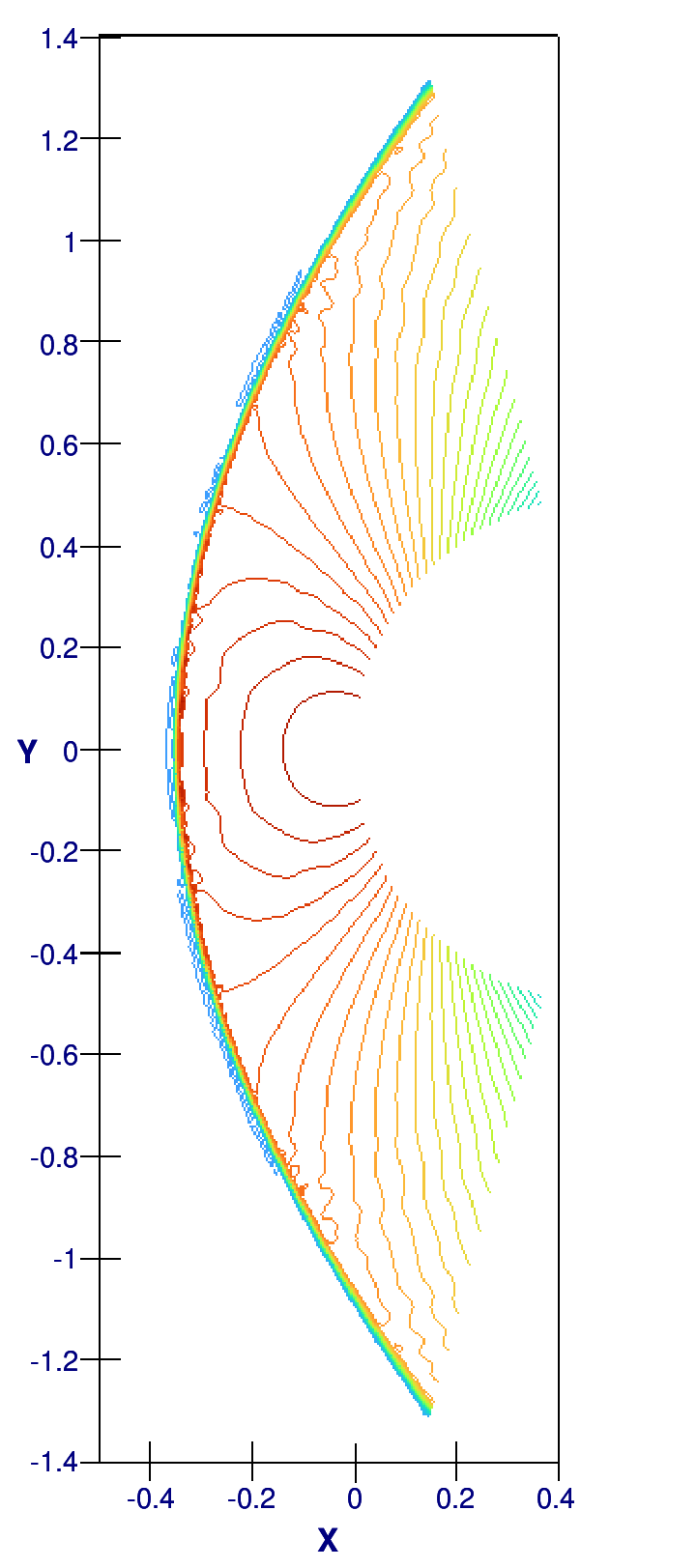}\label{fig:Con_every}}
\caption{Flow over a cylinder. Density contours of the numerical solutions obtained from both limiting approaches. 30 equally spaced contour lines from 1.4 to 6.03.}
\label{fig:BB_Contour}
\end{figure}

For the flow over a cylinder test case, the density contours of the numerical solutions obtained using both limiting approaches are presented in Figure (\ref{fig:BB_Contour}). The results indicate that the solutions from the two approaches closely match, particularly for the smaller threshold constant $K = 0.02$. For the threshold constant $K = 0.1$, the solution exhibits slight oscillations near the ends of bow shock.

Overall, the results of all steady-state test cases 1 to 5 demonstrate that limiting only in the troubled cells is sufficient to get solutions comparable to those obtained with the limiting-everywhere approach. This approach not only significantly reduces computational cost but also enhances convergence. However, as the number of troubled cells decreases, convergence improves while the solution exhibits increased oscillations.

For the unsteady test cases (Test cases 6 and 7), the computational domain is initialized using the given initial conditions. We apply the troubled-cell indicator at every iteration during the simulation to dynamically identify the troubled-cells.

\begin{figure}
\centering
\subfloat[Config3, $K = 0.02$]{\includegraphics[width=0.48\textwidth, height=0.5\textwidth]{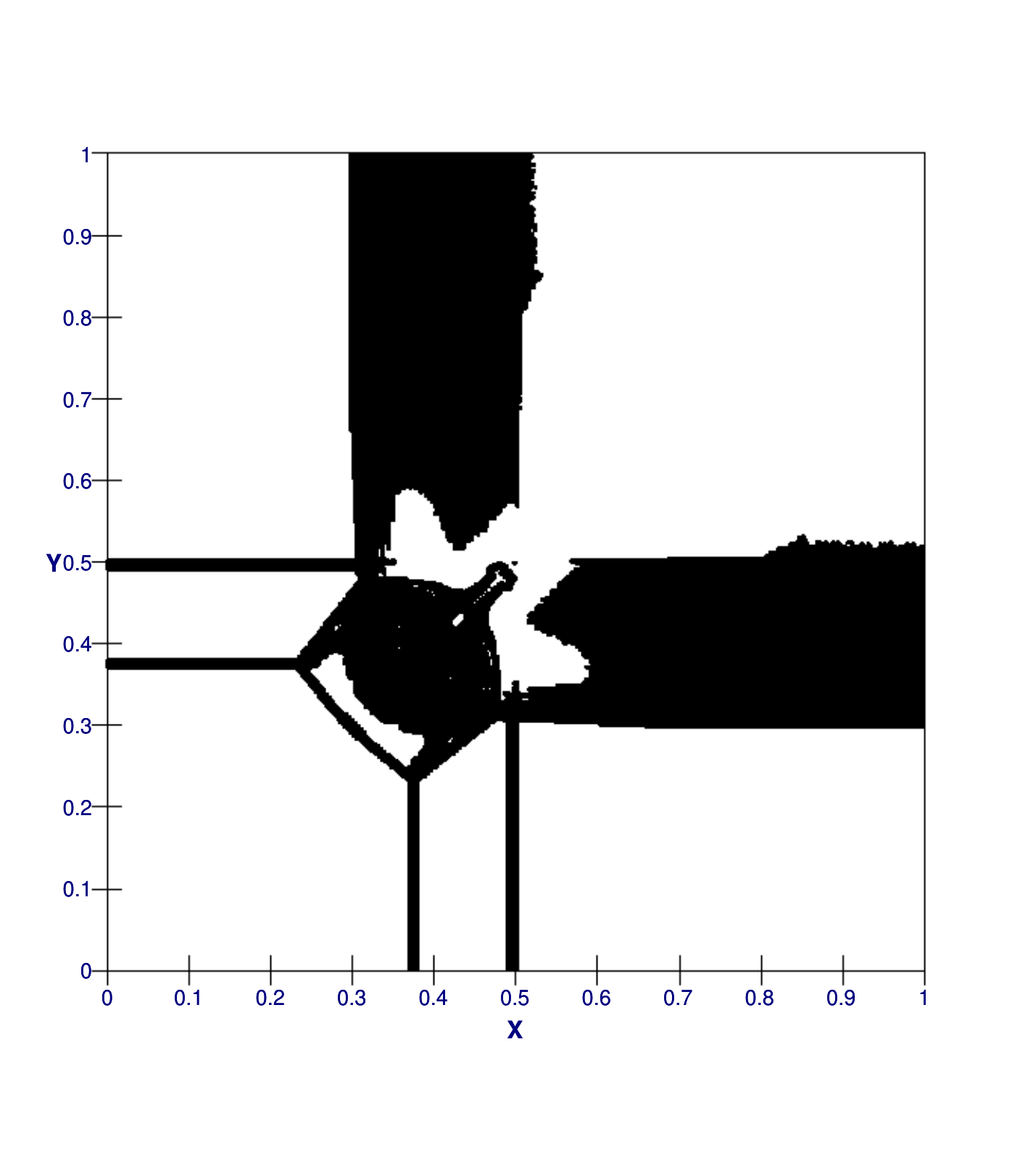}\label{fig:C3_K002_TC}} \hspace{0.3cm}
\subfloat[Config4, $K = 0.02$]{\includegraphics[width=0.48\textwidth, height=0.5\textwidth]{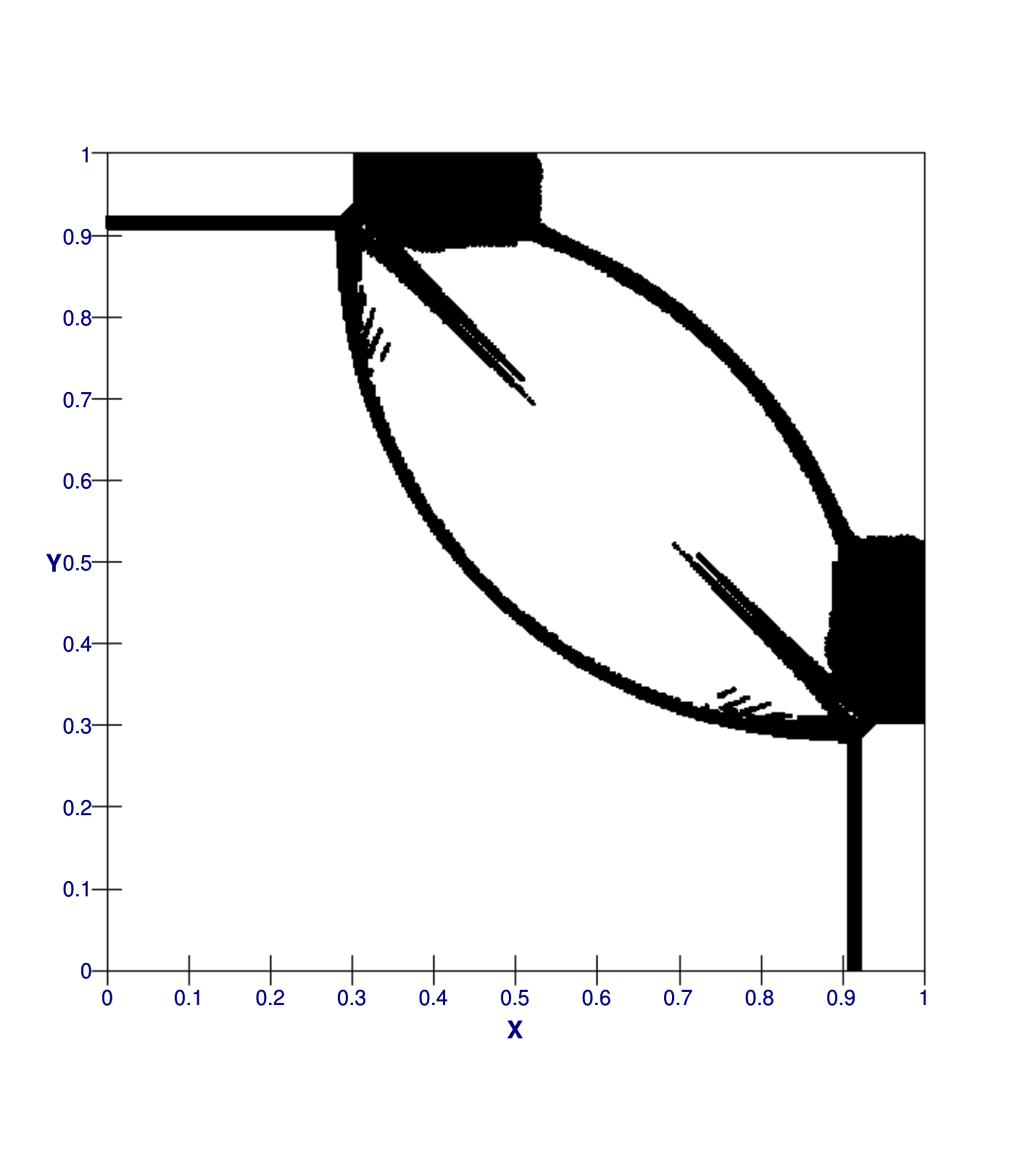}\label{fig:C4_K002_TC}} \hspace{0.3cm}
\subfloat[Config3, $K = 0.1$]{\includegraphics[width=0.48\textwidth, height=0.5\textwidth]{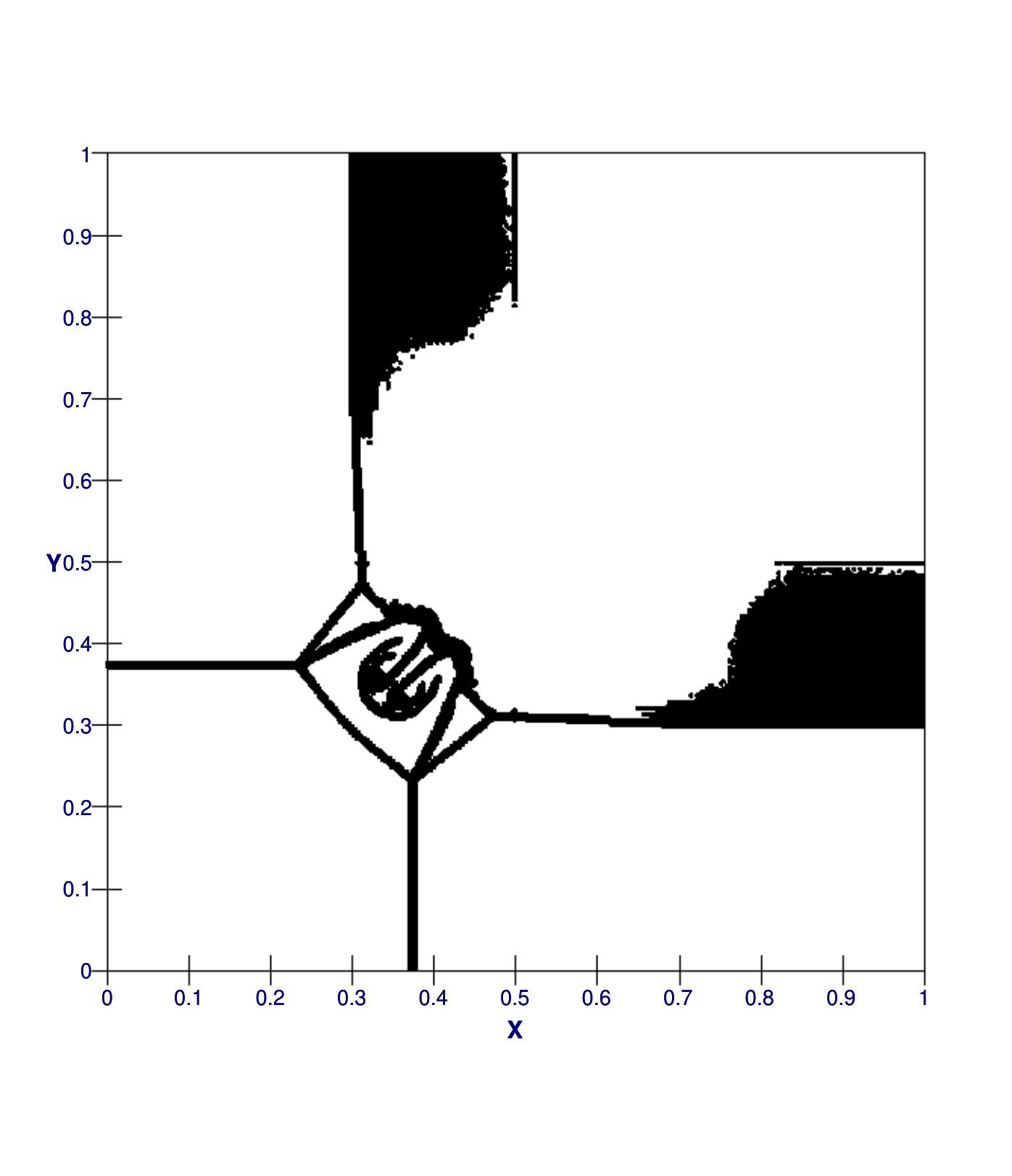}\label{fig:C3_K01_TC}} \hspace{0.3cm}
\subfloat[Config4, $K = 0.1$]{\includegraphics[width=0.48\textwidth, height=0.5\textwidth]{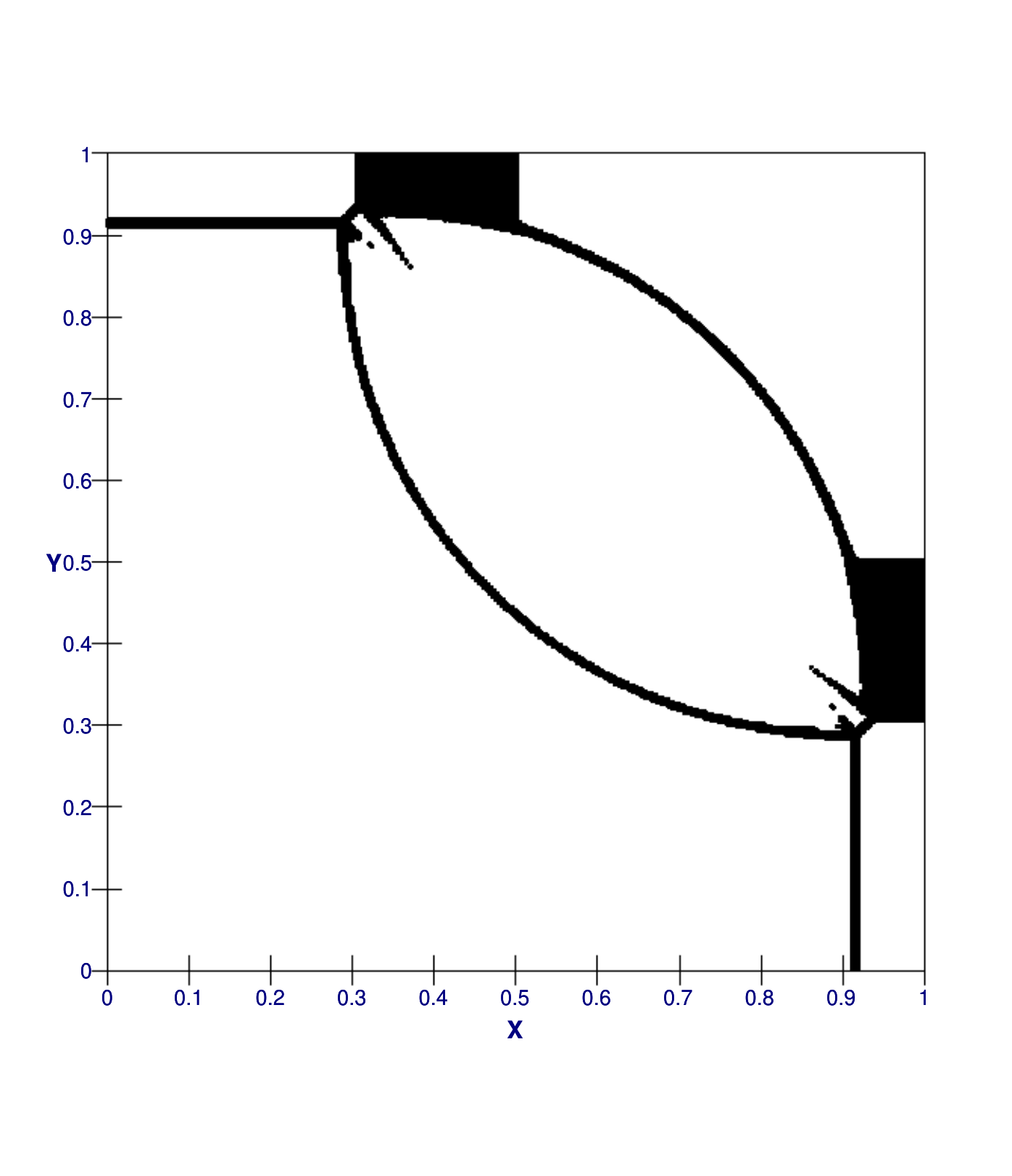}\label{fig:C4_K01_TC}}\\
\caption{2D Riemann problems. Troubled-cell region (in black) identified by the indicator at the last time step for threshold constants $K = 0.02$ and $K = 0.1$ for configurations 3 and 4. Percentage of number of troubled-cells identified in whole computational domain of 400 $\times$ 400 cells: (a) 24.1\% (b) 9.23\% (c) 11.12\% (d) 5.29\%}
\label{fig:Riemann_TC}
\end{figure}

\begin{figure}
\centering
\subfloat[Config3, $K = 0.02$]{\includegraphics[width=0.45\textwidth, height=0.4\textwidth]{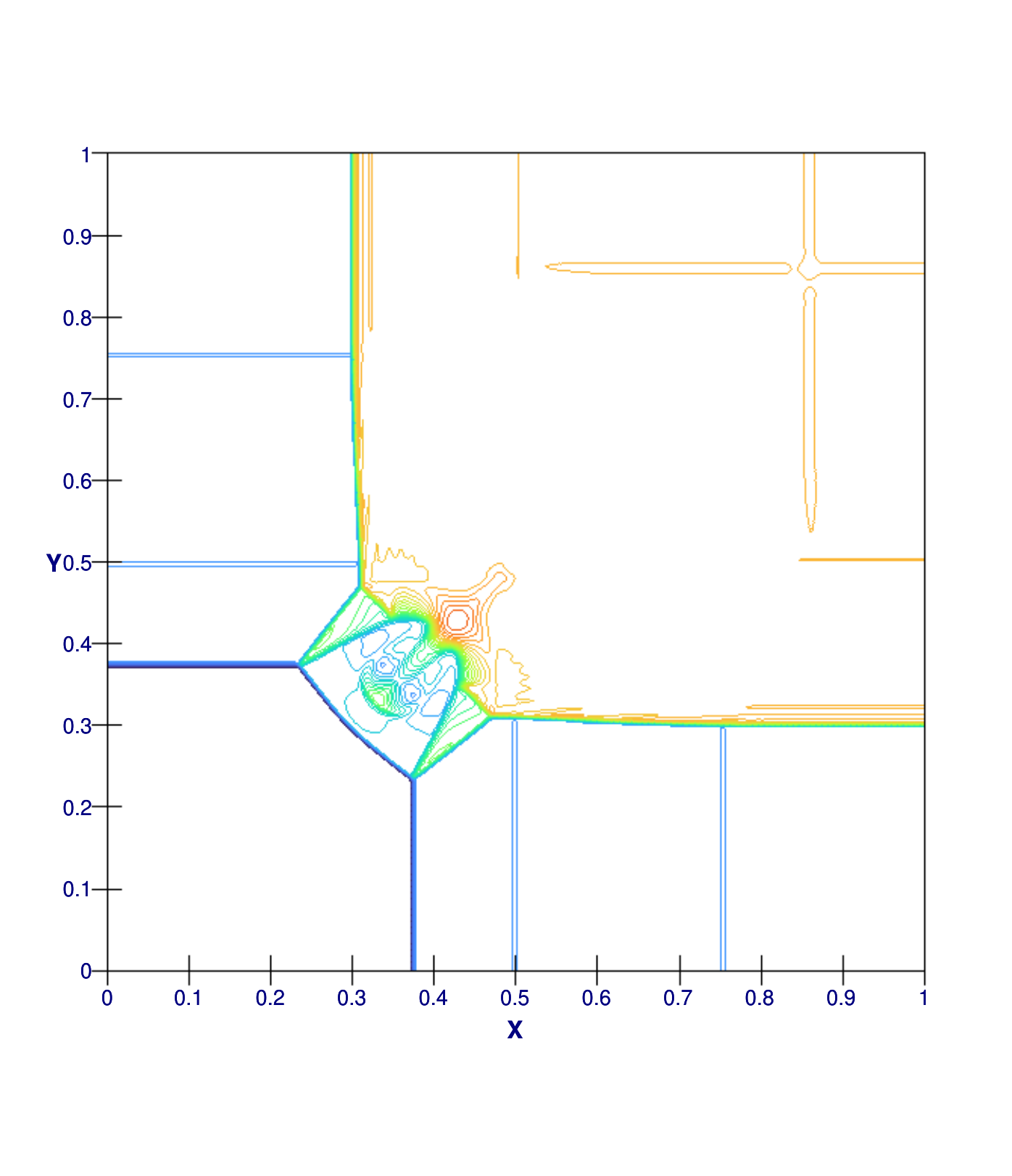}\label{fig:C3_K002_Con}} \hspace{0.2cm}
\subfloat[Config4, $K = 0.02$]{\includegraphics[width=0.45\textwidth, height=0.4\textwidth]{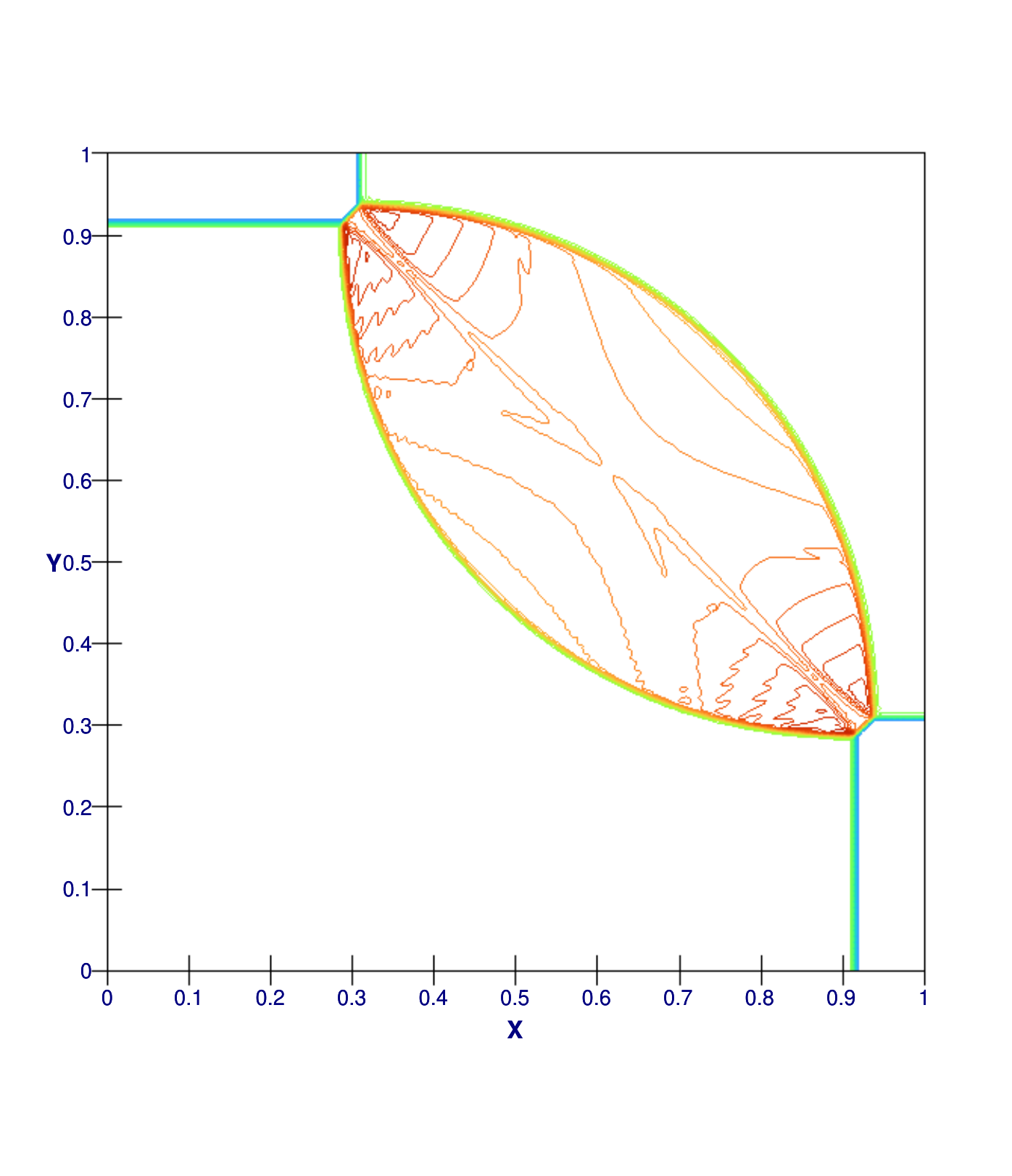}\label{fig:C4_K002_Con}} \\ \vspace{0.05cm}
\subfloat[Config3, $K = 0.1$]{\includegraphics[width=0.45\textwidth, height=0.4\textwidth]{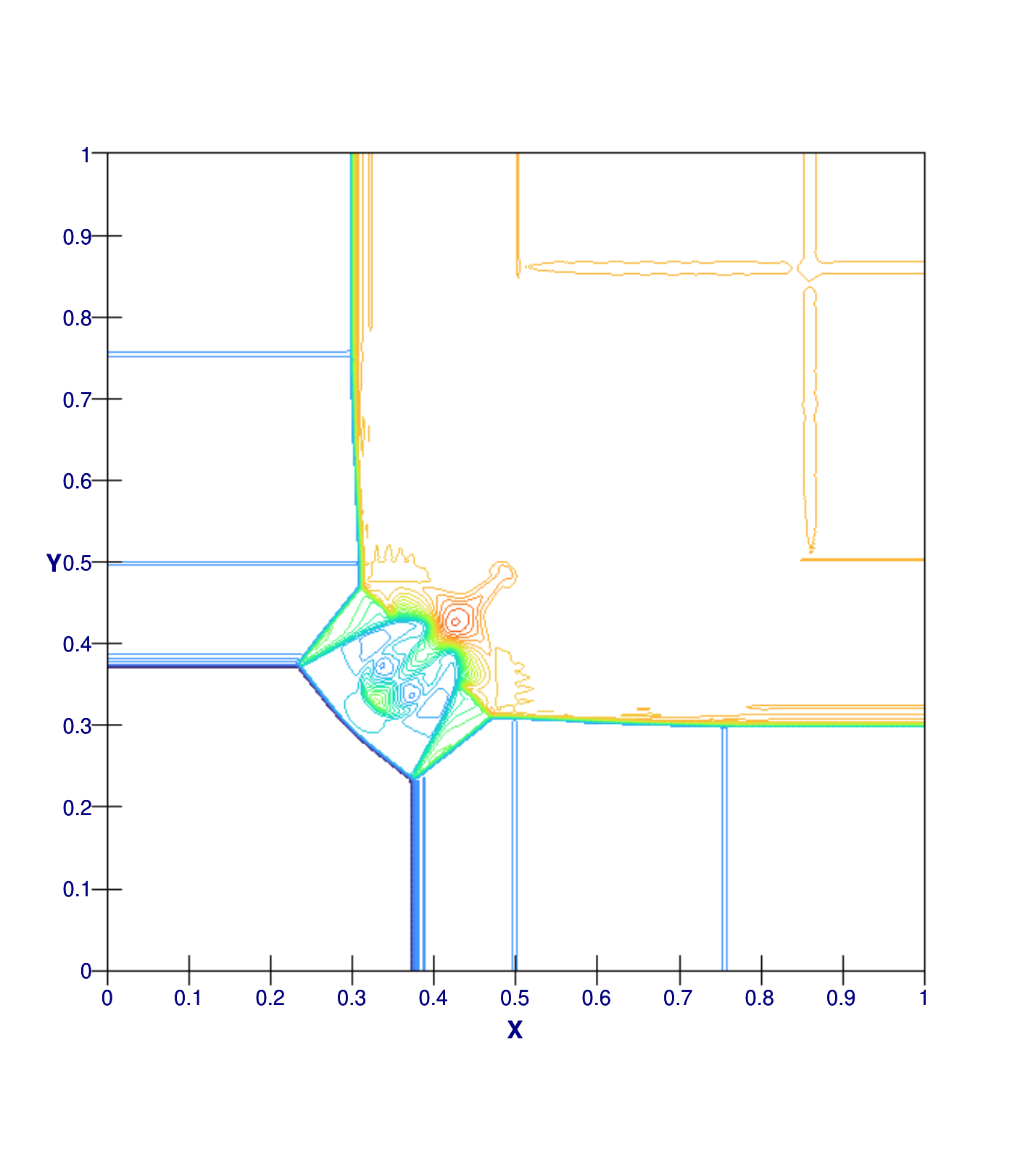}\label{fig:C3_K01_Con}} \hspace{0.2cm}
\subfloat[Config4, $K = 0.1$]{\includegraphics[width=0.45\textwidth, height=0.4\textwidth]{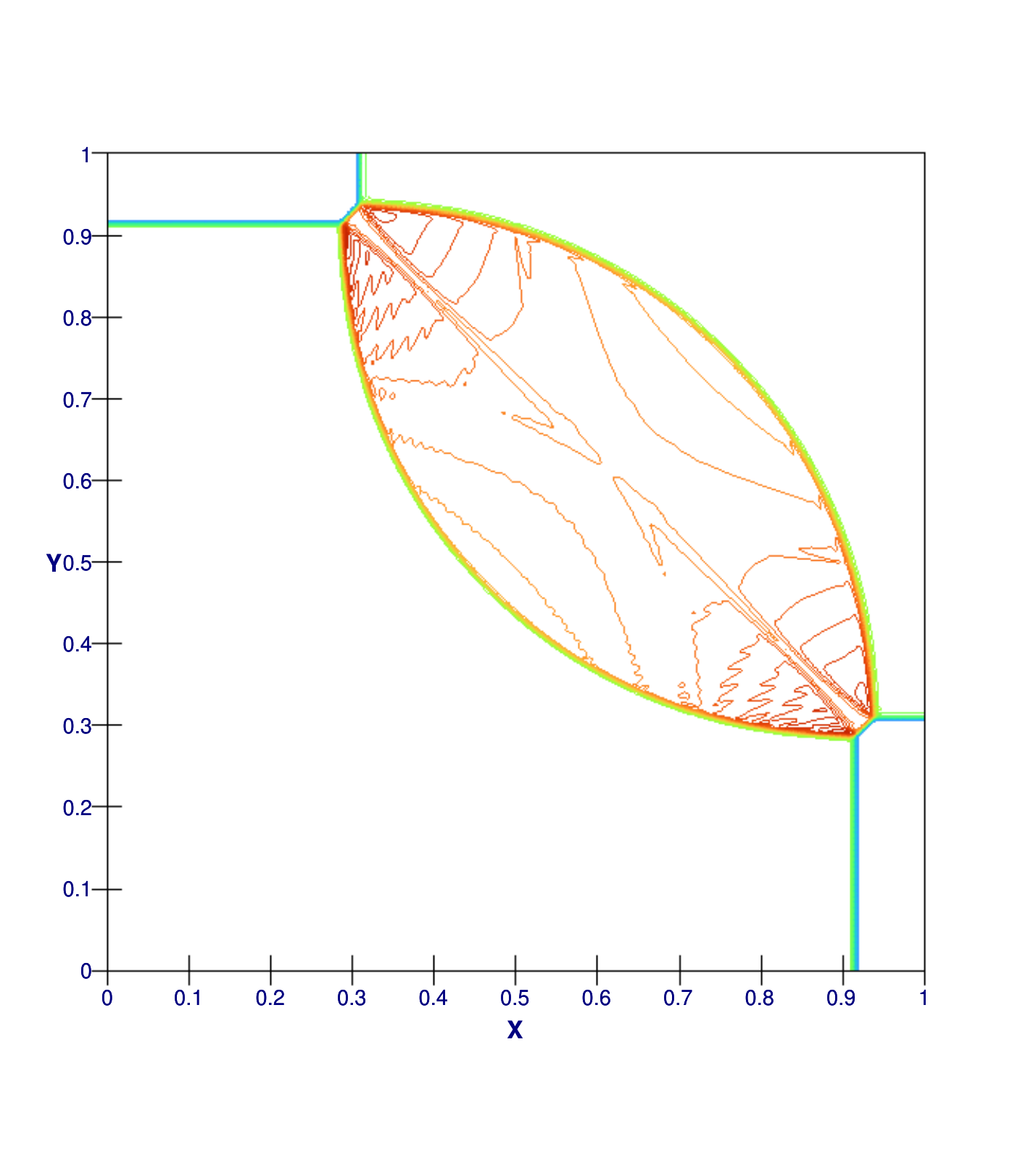}\label{fig:C4_K01_Con}}\\ \vspace{0.05cm}
\subfloat[Config3, Everywhere]{\includegraphics[width=0.45\textwidth, height=0.4\textwidth]{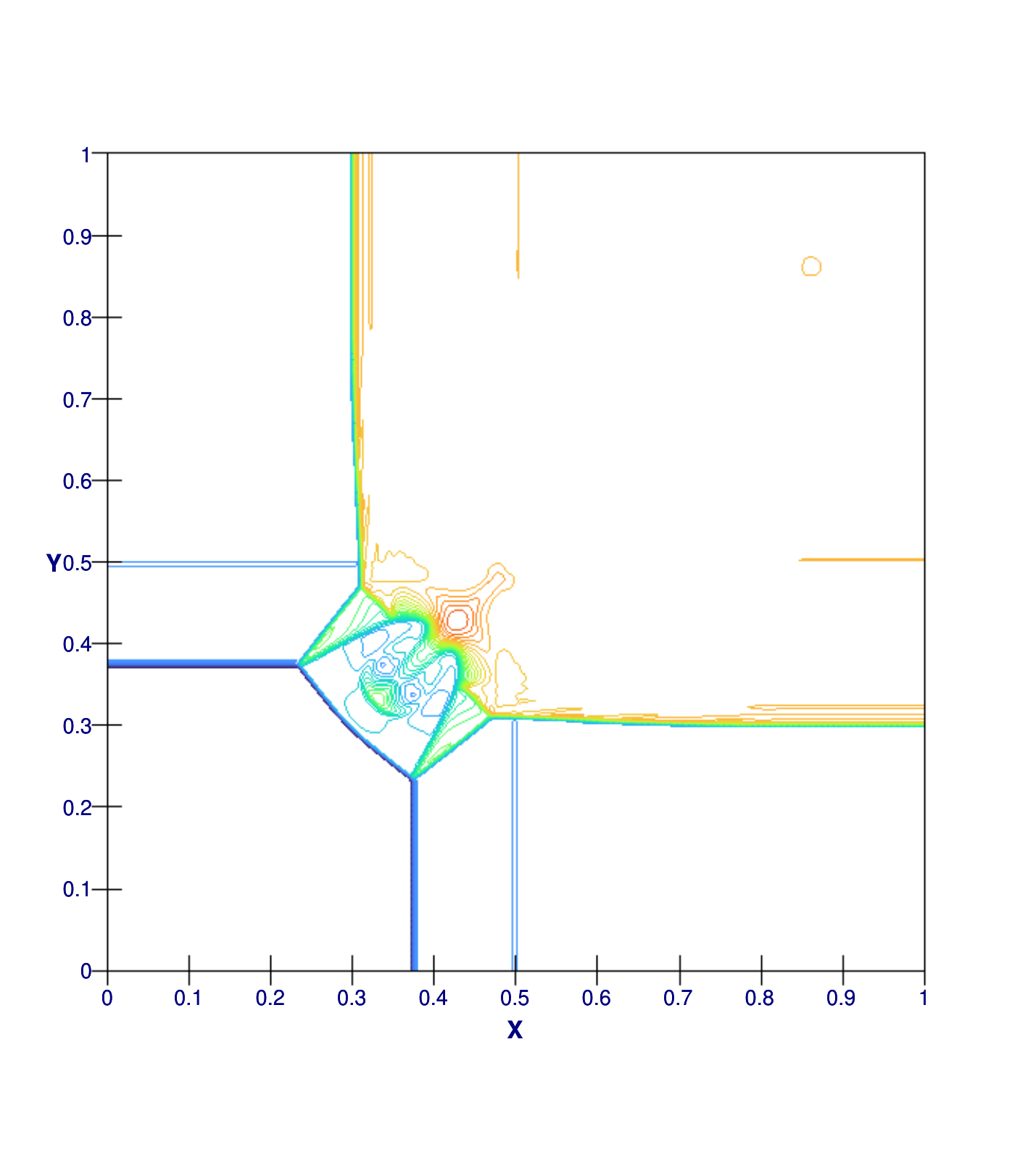}\label{fig:C3_every_Con}} \hspace{0.2cm}
\subfloat[Config4, Everywhere]{\includegraphics[width=0.45\textwidth, height=0.4\textwidth]{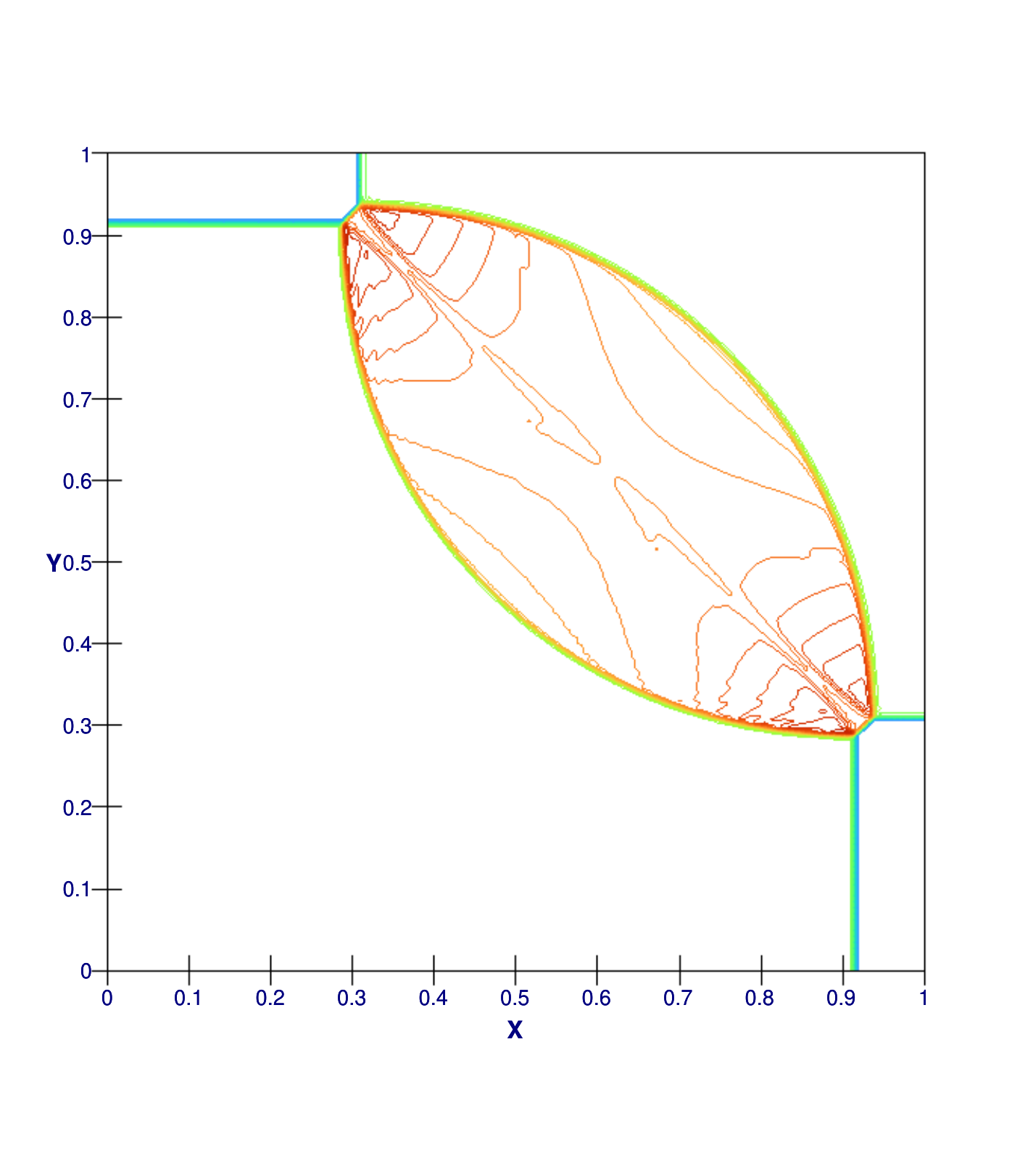}\label{fig:C4_every_Con}}\\
\caption{2D Riemann problems. Density contours of the numerical solutions obtained from both limiting approaches. 30 equally spaced contour lines from 0.13 to 1.79 for configuration 3. 30 equally spaced contour lines from 0.5 to 2.0 for configuration 4.}
\label{fig:Riemann_density}
\end{figure}

Figure (\ref{fig:Riemann_TC}) presents the troubled-cell region identified by the indicator at the final time step for the respective configurations 3 and 4 of 2D Riemann test case for threshold constants $K = 0.02$ and $K = 0.1$. Figure (\ref{fig:Riemann_density}) presents the density contours of the numerical solutions obtained using both limiting approaches for the same configurations.

\begin{figure}
\centering
\subfloat[]{\includegraphics[width=0.5\textwidth, height=0.5\textwidth]{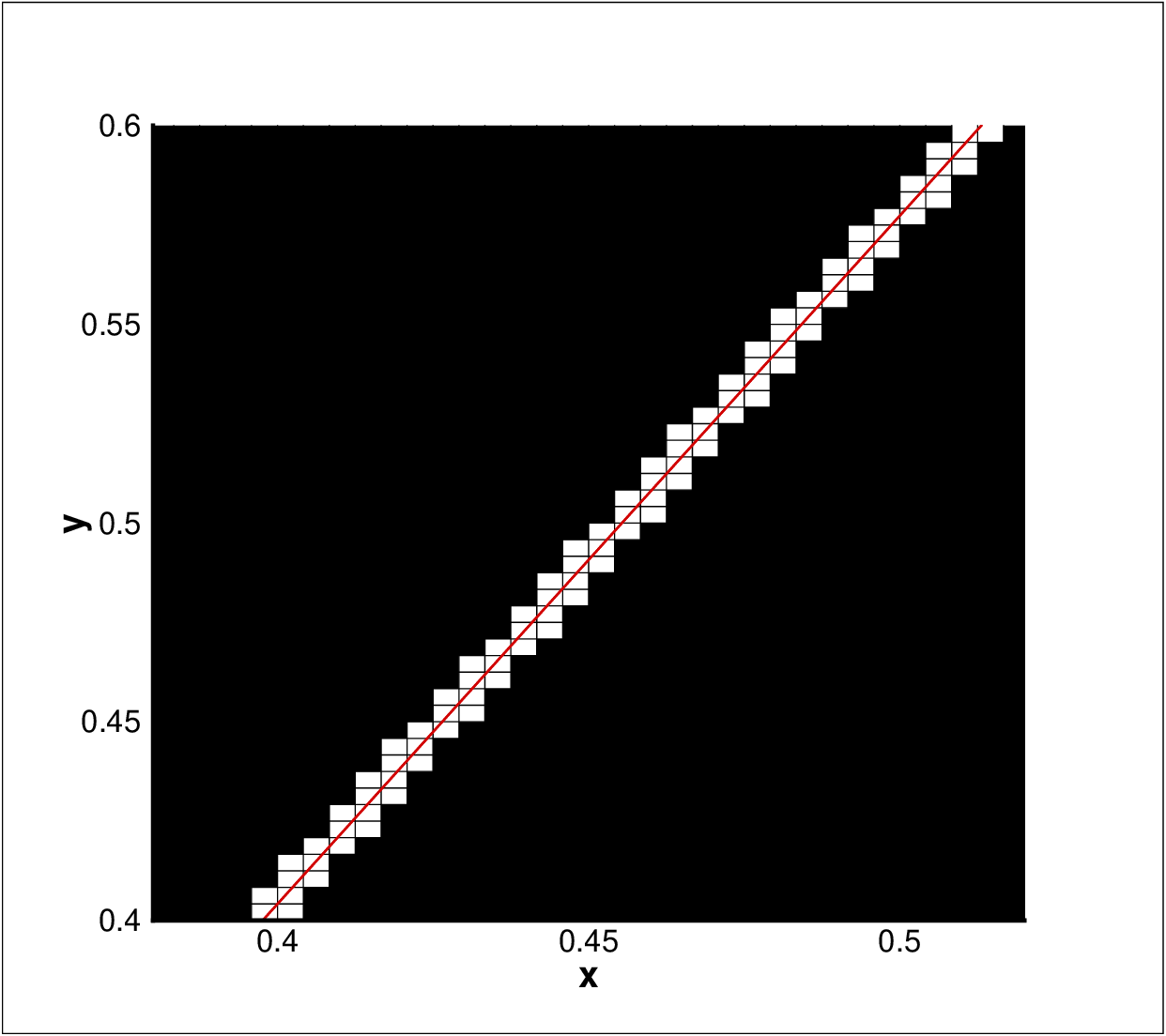}\label{fig:DMR1}}
\subfloat[]{\includegraphics[width=0.5\textwidth, height=0.5\textwidth]{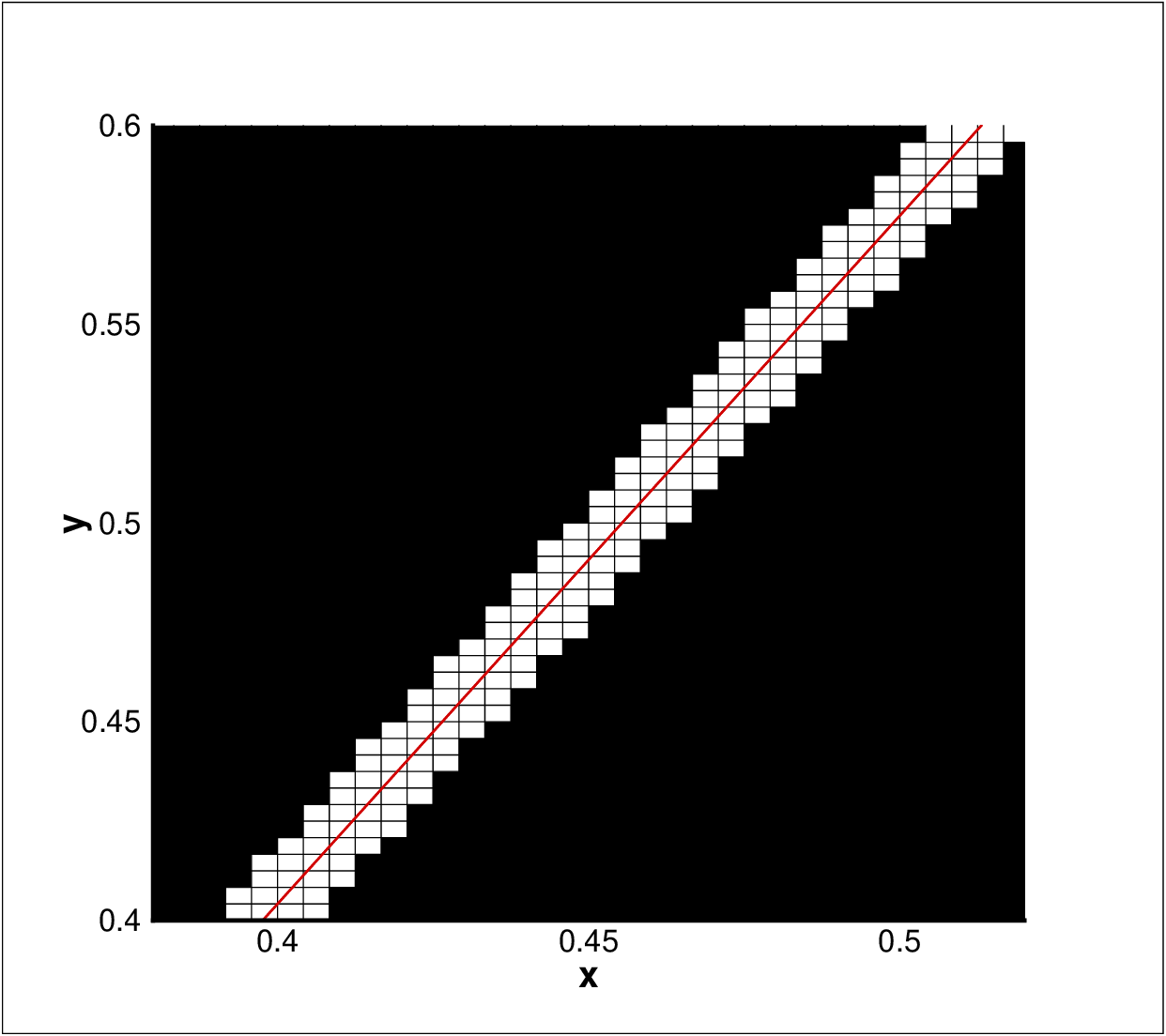}\label{fig:DMR2}}
\caption{Double Mach Reflection. (a) Zoomed-in view of troubled-cells identified by the indicator at the first time step. (b) Zoomed-in view of troubled-cells added on either side of the shock to the troubled-cells shown in (a).}
\label{fig:DMR_TC_first}
\end{figure}

Figure (\ref{fig:DMR1}) shows the troubled-cells identified during the first time step for the double Mach reflection test case. Since, we initialised the computational domain with the exact solution at $t = 0$, only a few troubled-cells are identified near the shock. Applying the limiter only in these identified cells resulted in the solution diverging at the end of first time step. To address this issue, we added a additional layer of troubled-cells on either side of the shock to the troubled-cells identified by the indicator as shown in Figure (\ref{fig:DMR2}). By applying the limiter in these additional troubled cells, along with those identified by the indicator, a stable solution is successfully obtained.

\begin{figure}
\centering
\subfloat[$K = 0.02$]{\includegraphics[width=\textwidth, height=0.5\textwidth]{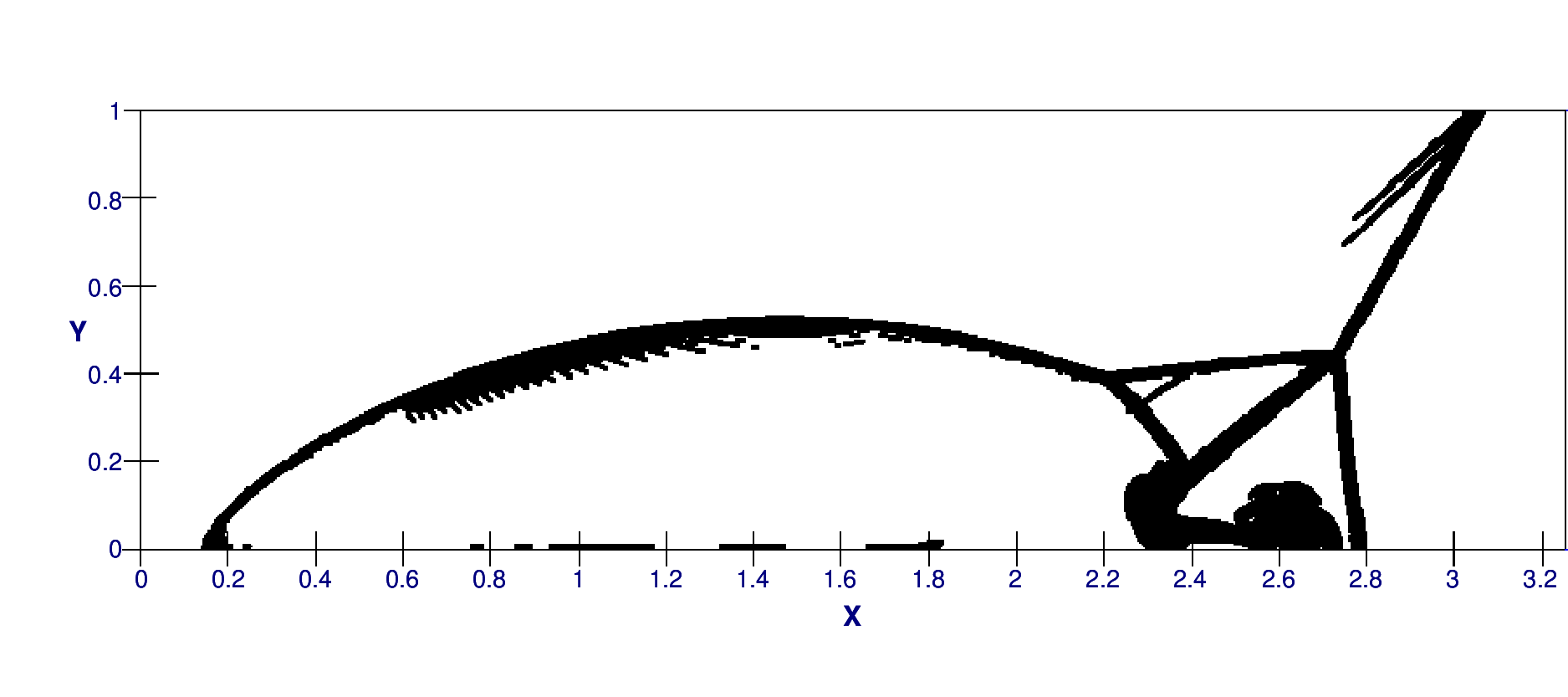}\label{fig:DMR_K002_TC}} \\
\subfloat[$K = 0.1$]{\includegraphics[width=\textwidth, height=0.5\textwidth]{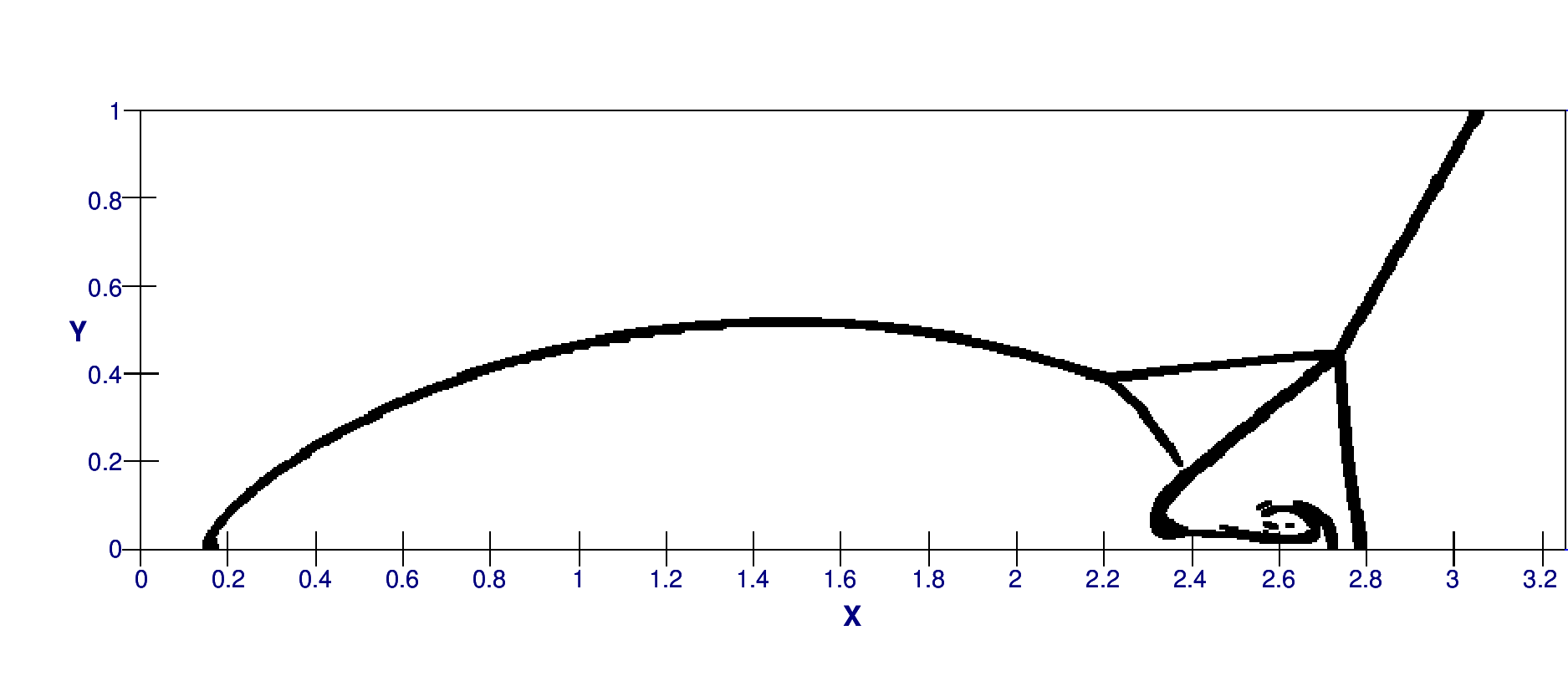}\label{fig:DMR_K01_TC}}
\caption{Double Mach Reflection. Troubled-cell region (in black) identified by the indicator at the last time step for threshold constants $K = 0.02$ and $K = 0.1$. Percentage of number of troubled-cells identified in whole computational domain of 960 $\times$ 240 cells: (a) 4.81\% (b) 2.07\%}
\label{fig:DMR_TC}
\end{figure}

\begin{figure}
\centering
\subfloat[$K = 0.02$]{\includegraphics[width=\textwidth, height=0.45\textwidth]{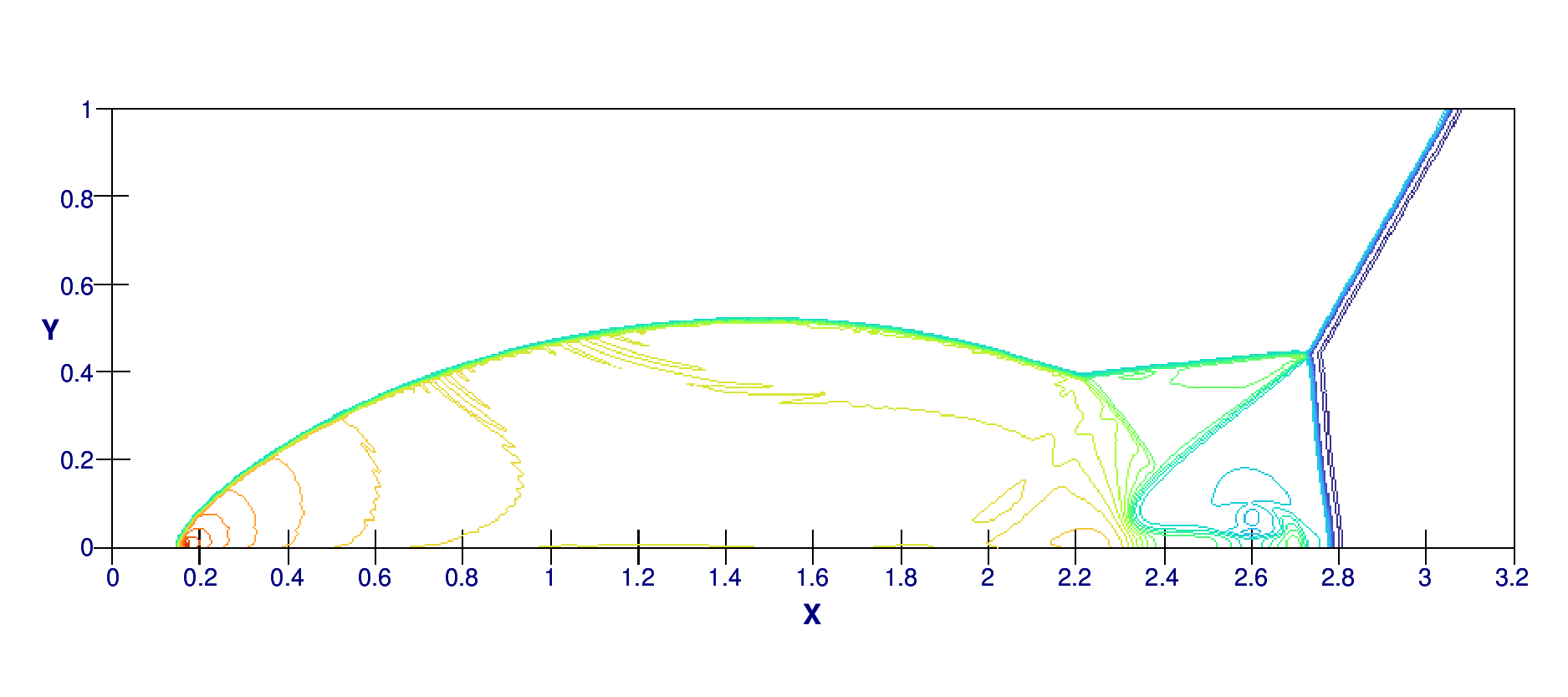}\label{fig:DMR_K002_Con}} \\
\subfloat[$K = 0.1$]{\includegraphics[width=\textwidth, height=0.45\textwidth]{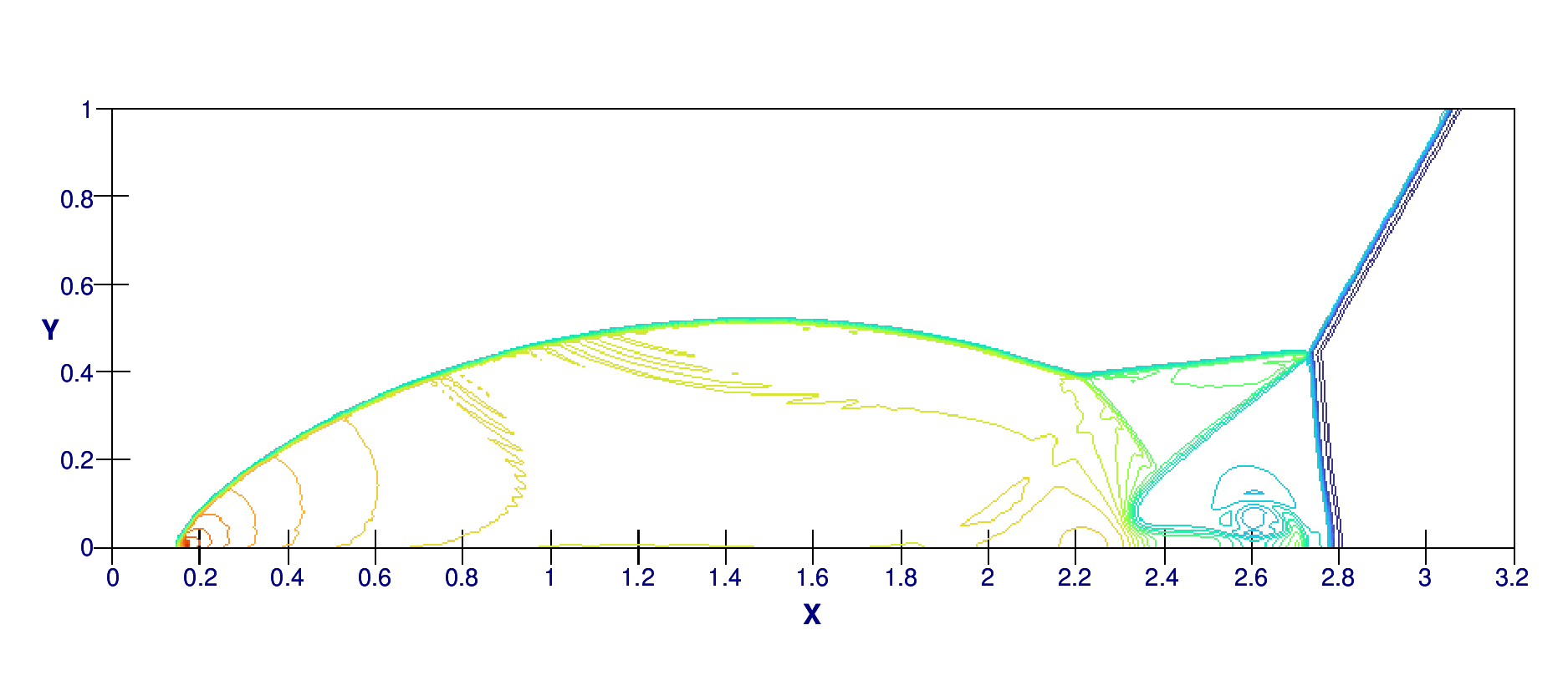}\label{fig:DMR_K01_Con}} \\
\subfloat[Everywhere]{\includegraphics[width=\textwidth, height=0.45\textwidth]{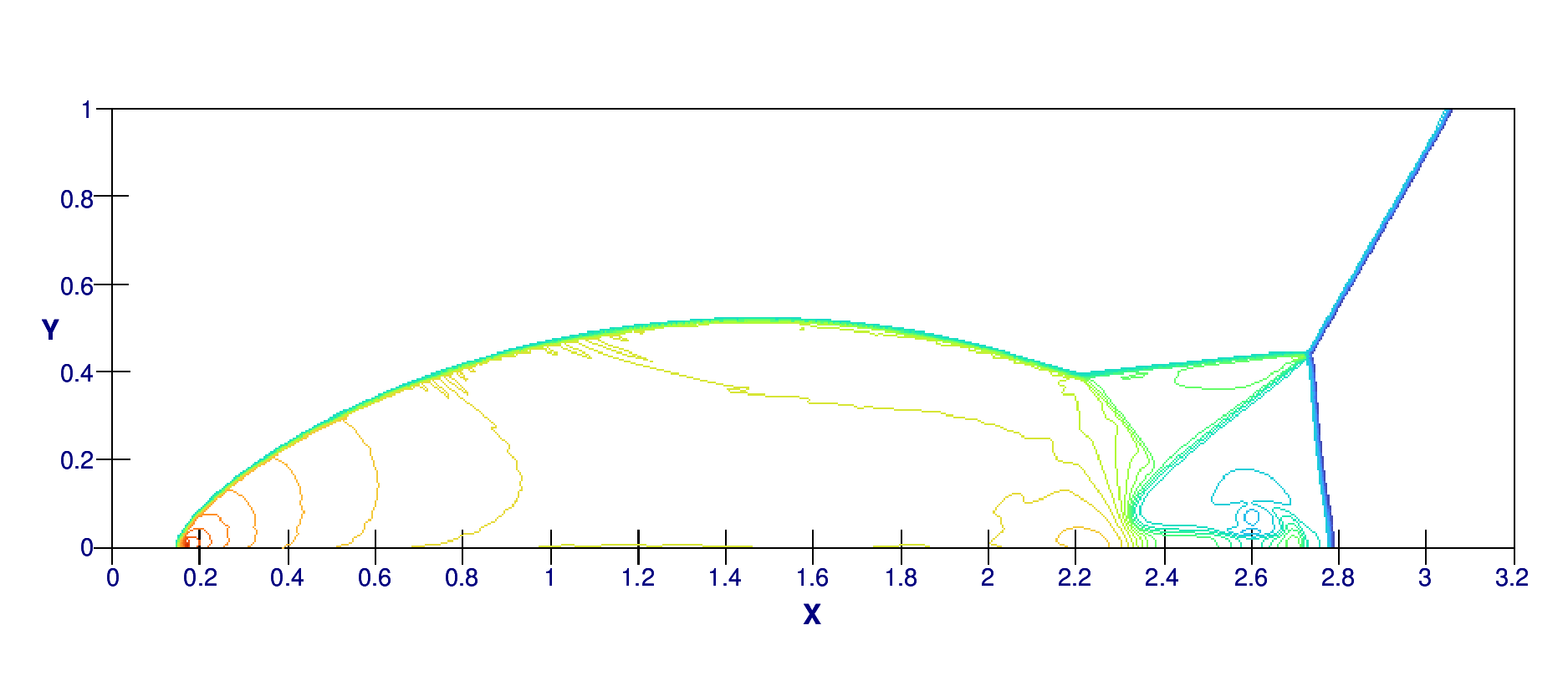}\label{fig:DMR_every_Con}} \\
\caption{Double Mach Reflection. Density contours of the numerical solutions obtained from both limiting approaches. 30 equally spaced contour lines from 1.4 to 23.8.}
\label{fig:DMR_density}
\end{figure}

Figure (\ref{fig:DMR_TC}) presents the troubled-cell region identified by the indicator at the final time step for threshold constants $K = 0.02$ and $K = 0.1$. Figure (\ref{fig:DMR_density}) presents the density contours of the numerical solutions obtained using both limiting approaches.

The results of both unsteady test cases indicate that both the limiting approaches capture the primary features of the solution with comparable accuracy. However, the contour lines for the limiting restricted region approach are less smooth than those observed for the limiting everywhere approach, with the wiggles in the contour lines becoming more pronounced for solutions computed for higher threshold constants.

\section{Conclusions}
\label{sec:Conclusion}
In this paper, we adapted a troubled-cell indicator from the discontinuous Galerkin (DG) framework to the finite volume method (FVM) framework for solving hyperbolic conservation laws. We also introduced a novel monotonicity parameter, $\mu$, that can be used to evaluate the quality of the solution in a neighbourhood containing a shock.

The troubled-cell indicator utilizes only the cell-averaged data of the target cell and its immediate neighbours. Its effectiveness has been demonstrated through various two-dimensional test problems, where it successfully captures shocks for threshold constants in the range of 0.01 to 0.1, particularly for $K = 0.05$.

Once we have information about the troubled-cells i.e., where the limiting is actually needed, we applied the limiter function only in those cells. The results reveal that limiting only in the troubled cells achieves comparable solution accuracy to the traditional limiting everywhere approach, while also improving convergence and reducing computational cost. Pros and cons of this limiting restricted region approach depending on the number of troubled-cells are listed as follows:

Limiting in a larger number of troubled-cells (this also includes limiting everywhere approach): \\

\begin{itemize}[align=left]
 \item[Pros]:
 \begin{itemize}
  \item In the limiting restricted region approach, increasing the number of troubled cells (i.e., reducing the threshold constant) reduces oscillations in the solution (Refer to Figure (\ref{fig:NOS_30_Density}) as an example).
  \item Produces a solution without oscillations and ensures monotonic, as observed in the limiting everywhere approach.
 \end{itemize}
 \item[Cons]:
 \begin{itemize}
  \item Extra computational effort.
  \item In the limiting restricted region approach, higher number of troubled cells (i.e., smaller threshold constant) results in more iterations needed to reach a steady-state solution (Refer to Figure (\ref{fig:NOS_30_RN}) as an example).
  \item Mostly convergence stalls, as seen in the limiting everywhere approach.
 \end{itemize}
\end{itemize}

Limiting in a smaller number of troubled-cells:
\begin{itemize}[align=left]
 \item[Pros]:
 \begin{itemize}
  \item Less computational effort.
  \item Improves convergence compared to the limiting everywhere approach.
  \item In the limiting restricted region approach, fewer troubled cells (i.e., larger threshold constant) lead to faster convergence to a steady-state solution (Refer to Figure (\ref{fig:NOS_30_RN}) as an example).
 \end{itemize}
 \item[Cons]:
 \begin{itemize}
  \item Solutions exhibit oscillations compared to those of the limiting everywhere approach.
  \item In the limiting restricted region approach, fewer number of troubled cells (i.e., larger threshold constant) results in more pronounced oscillations in the solution (Refer to Figure (\ref{fig:NOS_30_Density}) as an example).
  \item Insufficient troubled cells may cause the solution to become unstable, as observed in the case of double Mach reflection.
 \end{itemize}
\end{itemize}

In summary, employing the conventional higher-order MUSCL reconstruction with built-in limiting and no other enhancements results in a monotone solution. However, convergence to steady-state is not guaranteed. By limiting in a restricted region in the vicinity of the shock convergence to steady-state can be achieved and the resulting solution maybe monotone or have minimal oscillations. A systematic study of the effect of the number of troubled-cells on the solver will be presented in a follow up paper.

%
\bibliographystyle{ieeetr}
\bibliography{references}

\begin{thebibliography}{10}

\bibitem{FU2017}
G.~Fu and C.-W. Shu, ``A new troubled-cell indicator for discontinuous
  {G}alerkin methods for hyperbolic conservation laws,'' {\em Journal of
  Computational Physics}, vol.~347, pp.~305--327, 2017.

\bibitem{LeVeque2002}
R.~J. LeVeque, {\em Finite Volume Methods for Hyperbolic Problems}.
\newblock Cambridge Texts in Applied Mathematics, Cambridge University Press,
  2002.

\bibitem{VANLEER1979}
B.~{van Leer}, ``Towards the ultimate conservative difference scheme. {V}. a
  second-order sequel to {G}odunov's method,'' {\em Journal of Computational
  Physics}, vol.~32, no.~1, pp.~101--136, 1979.

\bibitem{VANLEER1974}
B.~{van Leer}, ``Towards the ultimate conservative difference scheme. {II}.
  monotonicity and conservation combined in a second-order scheme,'' {\em
  Journal of Computational Physics}, vol.~14, no.~4, pp.~361--370, 1974.

\bibitem{BISWAS1994}
R.~Biswas, K.~D. Devine, and J.~E. Flaherty, ``Parallel, adaptive finite
  element methods for conservation laws,'' {\em Applied Numerical Mathematics},
  vol.~14, no.~1, pp.~255--283, 1994.

\bibitem{WAN2022}
Y.~Wan and Y.~Xia, ``A hybrid {WENO} scheme for steady-state simulations of
  {E}uler equations,'' {\em Journal of Computational Physics}, vol.~463,
  p.~111292, 2022.

\bibitem{KRIVODONOVA2004}
L.~Krivodonova, J.~Xin, J.-F. Remacle, N.~Chevaugeon, and J.~Flaherty, ``Shock
  detection and limiting with discontinuous {G}alerkin methods for hyperbolic
  conservation laws,'' {\em Applied Numerical Mathematics}, vol.~48, no.~3,
  pp.~323--338, 2004.

\bibitem{SHU2005}
J.~Qiu and C.-W. Shu, ``{R}unge-{K}utta discontinuous {G}alerkin method using
  {WENO} limiters,'' {\em SIAM Journal on Scientific Computing}, vol.~26,
  no.~3, p.~907 – 929, 2005.

\bibitem{HEMKER1988}
P.~Hemker and B.~Koren, ``Multigrid, defect correction and upwind schemes for
  the steady {N}avier-{S}tokes equations,'' in {\em Numerical methods for fluid
  dynamics : international conference : proceedings, 3rd, based on the
  proceedings of a conference held in Oxford in March 1988} (K.~Morton and
  M.~Baines, eds.), The Institute of Mathematics and its Applications
  conference series. New series, pp.~153--170, 1988.

\bibitem{LIOU1996}
M.-S. Liou, ``A sequel to {AUSM}: {AUSM}+,'' {\em Journal of Computational
  Physics}, vol.~129, no.~2, pp.~364--382, 1996.

\bibitem{SHAROV1997}
D.~Sharov, K.~Nakahashi, D.~Sharov, and K.~Nakahashi, {\em Reordering of 3-D
  hybrid unstructured grids for vectorized {LU-SGS} {N}avier-{S}tokes
  computations}.

\bibitem{SHU1988}
C.-W. Shu and S.~Osher, ``Efficient implementation of essentially
  non-oscillatory shock-capturing schemes,'' {\em Journal of Computational
  Physics}, vol.~77, no.~2, pp.~439--471, 1988.

\bibitem{LAX1998}
P.~D. Lax and X.-D. Liu, ``Solution of two-dimensional {R}iemann problems of
  gas dynamics by positive schemes,'' {\em SIAM J. Sci. Comput.}, vol.~19,
  p.~319–340, Mar. 1998.

\bibitem{WOODWARD1984}
P.~Woodward and P.~Colella, ``The numerical simulation of two-dimensional fluid
  flow with strong shocks,'' {\em Journal of Computational Physics}, vol.~54,
  no.~1, pp.~115--173, 1984.

\bibitem{ZHANG2013}
S.-h. Zhang, X.-g. Deng, M.-l. Mao, and C.-W. Shu, ``Improvement of convergence
  to steady state solutions of {E}uler equations with weighted compact
  nonlinear schemes,'' {\em Acta Mathematicae Applicatae Sinica, English
  Series}, vol.~29, pp.~449--464, Jul 2013.

\end{thebibliography}

\end{document}